\RequirePackage{amsmath}
\documentclass[final,authorversion]{jfp-epi}
\usepackage[T1]{fontenc}
\usepackage[utf8]{inputenc}
\usepackage{tikz-cd}
\usetikzlibrary{datavisualization}
\pgfrealjobname{main}
\usepackage{quoting}
\usepackage{hyperref}
\usepackage{cleveref}
\usepackage{enumitem}
\setenumerate[1]{label=\thesection.\arabic*}
\setenumerate[2]{label*=\arabic*}
\usepackage{aliascnt}

\newaliascnt{speceq}{equation}
\crefname{speceq}{eq.}{eqs.}
\Crefname{speceq}{Eq.}{Eqs.}
\newcommand{\speclabel}[1]{\refstepcounter{speceq}\label{#1}}

\makeatletter
\@ifundefined{lhs2tex.lhs2tex.sty.read}%
  {\@namedef{lhs2tex.lhs2tex.sty.read}{}%
   \newcommand\SkipToFmtEnd{}%
   \newcommand\EndFmtInput{}%
   \long\def\SkipToFmtEnd#1\EndFmtInput{}%
  }\SkipToFmtEnd

\newcommand\ReadOnlyOnce[1]{\@ifundefined{#1}{\@namedef{#1}{}}\SkipToFmtEnd}
\usepackage{amstext}
\usepackage{stmaryrd}
\DeclareFontFamily{OT1}{cmtex}{}
\DeclareFontShape{OT1}{cmtex}{m}{n}
  {<5><6><7><8>cmtex8
   <9>cmtex9
   <10><10.95><12><14.4><17.28><20.74><24.88>cmtex10}{}
\DeclareFontShape{OT1}{cmtex}{m}{it}
  {<-> ssub * cmtt/m/it}{}

\DeclareFontShape{OT1}{cmtt}{bx}{n}
  {<5><6><7><8>cmtt8
   <9>cmbtt9
   <10><10.95><12><14.4><17.28><20.74><24.88>cmbtt10}{}
\DeclareFontShape{OT1}{cmtex}{bx}{n}
  {<-> ssub * cmtt/bx/n}{}

\newcommand{\Conid}[1]{\mathit{#1}}
\newcommand{\Varid}[1]{\mathit{#1}}
\newcommand{\anonymous}{\kern0.06em \vbox{\hrule\@width.5em}}

\renewcommand{\leq}{\leqslant}

\usepackage{polytable}

\@ifundefined{mathindent}%
  {\newdimen\mathindent\mathindent\leftmargini}%
  {}%

\def\resethooks{%
  \global\let\SaveRestoreHook\empty
  \global\let\ColumnHook\empty}
\newcommand*{\savecolumns}[1][default]%
  {\g@addto@macro\SaveRestoreHook{\savecolumns[#1]}}
\newcommand*{\restorecolumns}[1][default]%
  {\g@addto@macro\SaveRestoreHook{\restorecolumns[#1]}}
\newcommand*{\aligncolumn}[2]%
  {\g@addto@macro\ColumnHook{\column{#1}{#2}}}

\resethooks

\newcommand{\onelinecommentchars}{\quad-{}- }
\newcommand{\commentbeginchars}{\enskip\{-}
\newcommand{\commentendchars}{-\}\enskip}

\newcommand{\visiblecomments}{%
  \let\onelinecomment=\onelinecommentchars
  \let\commentbegin=\commentbeginchars
  \let\commentend=\commentendchars}

\newcommand{\invisiblecomments}{%
  \let\onelinecomment=\empty
  \let\commentbegin=\empty
  \let\commentend=\empty}

\visiblecomments

\newlength{\blanklineskip}
\newcommand{\hsindent}[1]{\quad}
\let\hspre\empty
\let\hspost\empty

\EndFmtInput
\makeatother
\ReadOnlyOnce{polycode.fmt}%
\makeatletter

\newcommand{\hsnewpar}[1]%
  {{\parskip=0pt\parindent=0pt\par\vskip #1\noindent}}

\newcommand{\hscodestyle}{}

\newcommand{\sethscode}[1]%
  {\expandafter\let\expandafter\hscode\csname #1\endcsname
   \expandafter\let\expandafter\endhscode\csname end#1\endcsname}

  {\par\noindent
   \advance\leftskip\mathindent
   \hscodestyle
   \let\\=\@normalcr
   \let\hspre\(\let\hspost\)%
   \pboxed}%
  {\endpboxed\)%
   \par\noindent
   \ignorespacesafterend}

  {\hsnewpar\abovedisplayskip
   \advance\leftskip\mathindent
   \hscodestyle
   \let\hspre\(\let\hspost\)%
   \pboxed}%
  {\endpboxed%
   \hsnewpar\belowdisplayskip
   \ignorespacesafterend}

  {\hsnewpar\abovedisplayskip
   \advance\leftskip\mathindent
   \hscodestyle
   \let\\=\@normalcr
   \(\pboxed}%
  {\endpboxed\)%
   \hsnewpar\belowdisplayskip
   \ignorespacesafterend}

\newcommand{\plainhs}{\sethscode{plainhscode}}

\plainhs

  {\hsnewpar\abovedisplayskip
   \advance\leftskip\mathindent
   \hscodestyle
   \let\\=\@normalcr
   \(\parray}%
  {\endparray\)%
   \hsnewpar\belowdisplayskip
   \ignorespacesafterend}

  {\parray}{\endparray}

  {\(\parray}{\endparray\)}

\def\codeframewidth{\arrayrulewidth}
\RequirePackage{calc}

  {\parskip=\abovedisplayskip\par\noindent
   \hscodestyle
   \arrayrulewidth=\codeframewidth
   \tabular{@{}|p{\linewidth-2\arraycolsep-2\arrayrulewidth-2pt}|@{}}%
   \hline\framedhslinecorrect\\{-1.5ex}%
   \let\endoflinesave=\\
   \let\\=\@normalcr
   \(\pboxed}%
  {\endpboxed\)%
   \framedhslinecorrect\endoflinesave{.5ex}\hline
   \endtabular
   \parskip=\belowdisplayskip\par\noindent
   \ignorespacesafterend}

\newcommand{\framedhslinecorrect}[2]%
  {#1[#2]}

  {\(\def\column##1##2{}%
   \let\>\undefined\let\<\undefined\let\\\undefined
   \newcommand\>[1][]{}\newcommand\<[1][]{}\newcommand\\[1][]{}%
   \def\fromto##1##2##3{##3}%
   }{\) }%

  {\let\orighscode=\hscode
   \let\origendhscode=\endhscode
   \def\endhscode{\def\hscode{\endgroup\def\@currenvir{hscode}\\}\begingroup}
   \orighscode\def\hscode{\endgroup\def\@currenvir{hscode}}}%
  {\origendhscode
   \global\let\hscode=\orighscode
   \global\let\endhscode=\origendhscode}%

\makeatother
\EndFmtInput
\def\doubleequals{\mathrel{\unitlength 0.01em
  \begin{picture}(78,40)
    \put(7,34){\line(1,0){25}} \put(45,34){\line(1,0){25}}
    \put(7,14){\line(1,0){25}} \put(45,14){\line(1,0){25}}
  \end{picture}}}

\newif\ifhideTODO
\hideTODOfalse

\DeclareMathAlphabet{\mathcal}{OMS}{cmsy}{m}{n}

\newcommand{\Real}{\mathbb{R}}

\newcommand{\Int}{\mathbb{Z}}

\def\commentbegin{\quad\{\ }
\def\commentend{\}}

\providecommand{\DONE}[1]{}

\newcommand{\testresult}[2]{%
\hsnewpar{\abovedisplayskip}
\begin{tabular}{@{}l@{}}
  \hscodestyle
  \ensuremath{\mathtt{\lambda{>}{~~}}} #1
  \\
  \hscodestyle
  \texttt{#2}
\end{tabular}
\hsnewpar{\belowdisplayskip}
\ignorespacesafterend
}

\jfpVolume{NN}
\jfpYear{2026}
\jfpArticle{XX}

\begin{document}

\newcommand{\shorttitle}{Optimization under uncertainty}


\title{Optimization under uncertainty: understanding orders and
  testing programs with specifications}

\author{PATRIK JANSSON}
\orcid{0000-0003-3078-1437}
\affiliation{%
  \institution{Chalmers University of Technology and University of Gothenburg}
  \city{Göteborg}
  \country{Sweden}
  \authoremail{patrikj@chalmers.se}
}

\author{NICOLA BOTTA}
\orcid{0000-0002-8923-2734}
\affiliation{%
  \institution{Potsdam Institute for Climate Impact Research}
  \city{Potsdam}
  \state{Brandenburg}
  \country{Germany}
  \authoremail{botta@pik-potsdam.de}
}

\affiliation{
  \institution{Chalmers University of Technology and University of Gothenburg}
  \city{Göteborg}
  \country{Sweden}
}

\author{TIM RICHTER}
\orcid{0000-0000-0000-0000} 
\affiliation{%
  \institution{Potsdam University}
  \city{Potsdam}
  \state{Brandenburg}
  \country{Germany}
  \authoremail{tim.richter@uni-potsdam.de}
}

\received[Submitted]{2026-03-14} 

\begin{abstract}
  One of the most ubiquitous problems in optimization is that of
  finding all the elements of a finite set at which a function~\ensuremath{\Varid{f}}
  attains its minimum (or maximum).
  When the codomain of~\ensuremath{\Varid{f}} is equipped with a total order, it is easy
  to specify, implement, and verify generic solutions to this problem.
  But what if~\ensuremath{\Varid{f}} is affected by uncertainties?
  What if one seeks values that minimize more than one objective, or
  if~\ensuremath{\Varid{f}} does not return a single result but a set of possible
  results, or even a probability distribution?
  Such situations are common in climate science, economics, and
  engineering.
  Developing trustworthy solution methods for optimization under
  uncertainty requires formulating and answering these questions
  rigorously, including deciding which order relations to apply in
  different cases.
  We show how functional programming can support this task, and apply
  it to specify and test solution methods for cases where optimization
  is affected by two conceptually different kinds of uncertainty:
  \emph{value} and \emph{functorial} uncertainty.
  We analyze the interplay of orders in these contexts,
  demonstrate how standard minimization generalizes to partial
  orders in the multi-objective setting and how it can be lifted via
  monotonicity conditions to handle functorial uncertainty.
\end{abstract}

\maketitle


\setlength{\mathindent}{1em}
\newcommand{\fixlengths}{\setlength{\abovedisplayskip}{6pt plus 1pt minus 1pt}\setlength{\belowdisplayskip}{6pt plus 1pt minus 1pt}}
\renewcommand{\hscodestyle}{\small\fixlengths}
\fixlengths


\section{Introduction}
\label{section:intro}

A common computational problem is finding all the elements of a finite
set at which a function attains its optimum.
Such problems naturally arise in engineering, economics, machine
learning, control theory, integrated assessment, climate policy
advice, etc.

For example, solving sequential decision problems via backward induction
\citep{bellman1957, bertsekas1995, puterman2014markov,
  2017_Botta_Jansson_Ionescu, BREDE_BOTTA_2021} requires computing an \emph{optimal}
 control for every element of a state space and at every
decision step.
A similar computation has to be done at each node of a neural network
when the network is trained via back propagation
\citep{10.5555/541177}.

Integrated assessment models for climate policy advice
\citep{Nordhaus2018, Nordhaus12261} typically seek policies that
minimize a so-called \emph{social welfare} function.
This function lumps together the (typically short-term) economic costs
and political risks of decarbonizing economies with the short- and
long-term impacts of climate change.
These impacts may affect, for example, economic output, human well-being,
and ecological systems \citep{TOL2024113922,
  KotzLevermannWenz2024, NealNewellPitman2025,
  BearparkHoganHsiang2025, Schoetz2025}.

\paragraph*{The totally ordered baseline.}

When the function to be minimized is known with certainty, the problem
can be formulated as that of implementing two functions:
\begin{hscode}\SaveRestoreHook
\column{B}{@{}>{\hspre}l<{\hspost}@{}}%
\column{3}{@{}>{\hspre}l<{\hspost}@{}}%
\column{13}{@{}>{\hspre}c<{\hspost}@{}}%
\column{13E}{@{}l@{}}%
\column{17}{@{}>{\hspre}l<{\hspost}@{}}%
\column{E}{@{}>{\hspre}l<{\hspost}@{}}%
\>[3]{}\Varid{min}{}\<[13]%
\>[13]{}\mathbin{::}{}\<[13E]%
\>[17]{}\Conid{Ord}\;\Varid{b}\Rightarrow (\Varid{a}\,\to\,\Varid{b})\,\to\,\Conid{Set}\;\Varid{a}\,\to\,\Varid{b}{}\<[E]%
\\[\blanklineskip]%
\>[3]{}\Varid{argmin}{}\<[13]%
\>[13]{}\mathbin{::}{}\<[13E]%
\>[17]{}\Conid{Ord}\;\Varid{b}\Rightarrow (\Varid{a}\,\to\,\Varid{b})\,\to\,\Conid{Set}\;\Varid{a}\,\to\,\Conid{Set}\;\Varid{a}{}\<[E]%
\ColumnHook
\end{hscode}\resethooks
The idea is that, for any \ensuremath{\Varid{f}}, \ensuremath{\Varid{as}} of suitable types, \ensuremath{\Varid{min}\;\Varid{f}\;\Varid{as}\in \Varid{map}\;\Varid{f}\;\Varid{as}} and \ensuremath{\Varid{argmin}\;\Varid{f}\;\Varid{as}\subseteq \Varid{as}} shall fulfil three conditions:

\begin{enumerate}
\item \ensuremath{\forall \Varid{a}\in \Varid{as}}. \ \ensuremath{\Varid{min}\;\Varid{f}\;\Varid{as}\mathrel{\leq}\Varid{f}\;\Varid{a}}                               \speclabel{eq1.1}
\item \ensuremath{\forall \Varid{a}\in \Varid{as}}. \ \ensuremath{\Varid{a}\in \Varid{argmin}\;\Varid{f}\;\Varid{as}} \ \ensuremath{\Rightarrow } \ \ensuremath{\Varid{f}\;\Varid{a}\mathrel{=}\Varid{min}\;\Varid{f}\;\Varid{as}}   \speclabel{eq1.2}
\item \ensuremath{\forall \Varid{a}\in \Varid{as}}. \ \ensuremath{\Varid{f}\;\Varid{a}\mathrel{=}\Varid{min}\;\Varid{f}\;\Varid{as}}  \ \ensuremath{\Rightarrow } \ \ensuremath{\Varid{a}\in \Varid{argmin}\;\Varid{f}\;\Varid{as}}  \speclabel{eq1.3}
\end{enumerate}

\noindent
The function \ensuremath{\Varid{f}\ \mathop{:}\ \Conid{A}\,\to\,\Conid{B}} is called the \emph{objective} function.
As mentioned in the abstract, the codomain of \ensuremath{\Varid{f}} is required to be
equipped with a \emph{total} order that we denote with \ensuremath{(\mathrel{\leq})}.
Totality of the order on \ensuremath{\Conid{B}} is important here: it ensures that any
two objective values are comparable, and hence that a unique minimum
value exists for every finite non-empty \ensuremath{\Varid{as}}.
As we will see later, relaxing this assumption, by allowing only a
partial order, fundamentally changes the nature of the optimization
problem and gives rise to sets of incomparable optima rather than a
single best value.

Together with \ensuremath{\Varid{min}\;\Varid{f}\;\Varid{as}\in \Varid{map}\;\Varid{f}\;\Varid{as}} and \ensuremath{\Varid{argmin}\;\Varid{f}\;\Varid{as}\subseteq \Varid{as}},
conditions \cref{eq1.1}-\cref{eq1.3} specify that \ensuremath{\Varid{min}\;\Varid{f}\;\Varid{as}} computes
the minimal objective value attained by \ensuremath{\Varid{f}} on \ensuremath{\Varid{as}}, while \ensuremath{\Varid{argmin}\;\Varid{f}\;\Varid{as}} computes precisely the set of elements in \ensuremath{\Varid{as}} at which this minimum
is attained.
For clarity of reference and for property-based testing, we
deliberately state \cref{eq1.2} and \cref{eq1.3} separately, rather
than combining them into a single biconditional
($\ensuremath{\forall \Varid{a}\in \Varid{as}}. \ \ensuremath{\Varid{a}\in \Varid{argmin}\;\Varid{f}\;\Varid{as}} \ \ensuremath{\Leftrightarrow } \ \ensuremath{\Varid{f}\;\Varid{a}\mathrel{=}\Varid{min}\;\Varid{f}\;\Varid{as}}$).

Note that, compared to more common definitions of \ensuremath{\Varid{min}} and \ensuremath{\Varid{argmin}} in
the literature, we have here an additional second argument of type
\ensuremath{\Conid{Set}\;\Varid{a}} which specifies the finite set on which the function shall be
minimized.

For ease of presentation, we use the type synonym \ensuremath{\mathbf{type}\;\Conid{Set}\;\Varid{a}\mathrel{=}[\mskip1.5mu \Varid{a}\mskip1.5mu]}.
Lists are treated extensionally as finite sets; in particular, order
and multiplicity are ignored, and we will often implicitly assume
non-empty sets.
In Haskell code snippets, type variables are written in lower case,
thus \ensuremath{\Varid{a}} and \ensuremath{\Varid{b}} in \ensuremath{\Conid{Set}\;\Varid{a}\,\to\,\Varid{b}} denote types.
The same lower case name in Haskell can denote both a type variable
and a plain variable.
For example, \ensuremath{\Varid{a}} in \ensuremath{\Varid{f}\;\Varid{a}} is a plain variable, and we picked the name as
a reminder that it is of type \ensuremath{\Varid{a}}.
In the surrounding text, we follow standard mathematical convention
and write type variables as \ensuremath{\Conid{A}} and \ensuremath{\Conid{B}}.

If equality in \ensuremath{\Conid{A}} is decidable, all three requirements are decidable
and polymorphic implementations of \ensuremath{\Varid{min}}, \ensuremath{\Varid{argmin}} return values
in \ensuremath{\Varid{map}\;\Varid{f}\;\Varid{as}} and \ensuremath{\Varid{as}} by parametricity.
For example, these implementations
\begin{hscode}\SaveRestoreHook
\column{B}{@{}>{\hspre}l<{\hspost}@{}}%
\column{3}{@{}>{\hspre}l<{\hspost}@{}}%
\column{18}{@{}>{\hspre}c<{\hspost}@{}}%
\column{18E}{@{}l@{}}%
\column{21}{@{}>{\hspre}l<{\hspost}@{}}%
\column{E}{@{}>{\hspre}l<{\hspost}@{}}%
\>[3]{}\Varid{min}\;\Varid{f}\;\Varid{as}{}\<[18]%
\>[18]{}\mathrel{=}{}\<[18E]%
\>[21]{}\Varid{minimum}\;(\Varid{map}\;\Varid{f}\;\Varid{as}){}\<[E]%
\\[\blanklineskip]%
\>[3]{}\Varid{argmin}\;\Varid{f}\;\Varid{as}{}\<[18]%
\>[18]{}\mathrel{=}{}\<[18E]%
\>[21]{}\mathbf{let}\;\Varid{m}\mathrel{=}\Varid{min}\;\Varid{f}\;\Varid{as}\;\mathbf{in}\;\Varid{filter}\;(\lambda \Varid{a}\,\to\,\Varid{f}\;\Varid{a}\doubleequals\Varid{m})\;\Varid{as}{}\<[E]%
\ColumnHook
\end{hscode}\resethooks
can be tested to fulfil the specification \cref{eq1.1}-\cref{eq1.3}
with QuickCheck \citep{Hughes2010} for arbitrary \ensuremath{\Varid{f}\ \mathop{:}\ \Int\,\to\,\Int}.
We use the helper \ensuremath{\Varid{forAll}} to keep the test and specification close.
\DONE{Perhaps make \ensuremath{\Varid{p}} top-level also here, as in MOO}
\DONE{In the introduction, before value uncertainty, define notions of
  order and discuss the consequences of losing totality for
  \ensuremath{\mathrel{\leq}}. Connect to notions of uncertainty.}
\DONE{Define and type \ensuremath{\Varid{forAll}}.}
\DONE{Nicola: why do we introduce \ensuremath{\Varid{forAll}} here? Is it meant to be an
  example of type variables in Haskell? We already have discussed two
  examples: \ensuremath{\Varid{min}} and \ensuremath{\Varid{argmin}}\dots
  Patrik: It is for readability:
    \ensuremath{\Varid{all}\;(\Varid{some}\;\Varid{long}\;\Varid{predicate})\;\Varid{as}}
  is much more readable if written in this form
    \ensuremath{\Varid{forAll}\;\Varid{as}\;(\Varid{some}\;\Varid{long}\;\Varid{predicate})}
  It is not a major thing, but I like how close the Haskell code gets to the spec when using \ensuremath{\Varid{forAll}}.
  And not many are immediately familiar with \ensuremath{\Varid{all}} anyway.
  I suggest we type it without the \ensuremath{\Conid{Foldable}} constraint: with \ensuremath{\Varid{t}} specialised to \ensuremath{\Conid{Set}\mathrel{=}\Conid{List}}.
}
\begin{hscode}\SaveRestoreHook
\column{B}{@{}>{\hspre}l<{\hspost}@{}}%
\column{3}{@{}>{\hspre}l<{\hspost}@{}}%
\column{11}{@{}>{\hspre}l<{\hspost}@{}}%
\column{14}{@{}>{\hspre}l<{\hspost}@{}}%
\column{17}{@{}>{\hspre}c<{\hspost}@{}}%
\column{17E}{@{}l@{}}%
\column{18}{@{}>{\hspre}l<{\hspost}@{}}%
\column{20}{@{}>{\hspre}l<{\hspost}@{}}%
\column{E}{@{}>{\hspre}l<{\hspost}@{}}%
\>[3]{}\Varid{prop}_{1.1}\;{}\<[11]%
\>[11]{}\Varid{f}\;{}\<[14]%
\>[14]{}\Varid{as}{}\<[18]%
\>[18]{}\mathrel{=}\mathbf{let}\;\Varid{m}\mathrel{=}\Varid{min}\;\Varid{f}\;\Varid{as}\;\mathbf{in}\;\Varid{forAll}\;\Varid{as}\;(\lambda \Varid{a}\,\to\,\Varid{m}\mathrel{\leq}\Varid{f}\;\Varid{a}){}\<[E]%
\\
\>[3]{}\Varid{prop}_{1.2}\;{}\<[11]%
\>[11]{}\Varid{f}\;{}\<[14]%
\>[14]{}\Varid{as}{}\<[18]%
\>[18]{}\mathrel{=}\Varid{forAll}\;\Varid{as}\;(\lambda \Varid{a}\,\to\,(\Varid{a}\in \Varid{argmin}\;\Varid{f}\;\Varid{as})\Rightarrow (\Varid{f}\;\Varid{a}\doubleequals\Varid{min}\;\Varid{f}\;\Varid{as})){}\<[E]%
\\
\>[3]{}\Varid{prop}_{1.3}\;{}\<[11]%
\>[11]{}\Varid{f}\;{}\<[14]%
\>[14]{}\Varid{as}{}\<[18]%
\>[18]{}\mathrel{=}\Varid{forAll}\;\Varid{as}\;(\lambda \Varid{a}\,\to\,(\Varid{f}\;\Varid{a}\doubleequals\Varid{min}\;\Varid{f}\;\Varid{as})\Rightarrow (\Varid{a}\in \Varid{argmin}\;\Varid{f}\;\Varid{as})){}\<[E]%
\\[\blanklineskip]%
\>[3]{}\Varid{forAll}\mathbin{::}\Conid{Set}\;\Varid{a}\,\to\,(\Varid{a}\,\to\,\Conid{Bool})\,\to\,\Conid{Bool}{}\<[E]%
\\
\>[3]{}\Varid{forAll}\mathrel{=}\Varid{flip}\;\Varid{all}{}\<[E]%
\\[\blanklineskip]%
\>[3]{}\Varid{test}_{1.1}\;{}\<[11]%
\>[11]{}\Varid{n}\;{}\<[14]%
\>[14]{}\Varid{f}{}\<[17]%
\>[17]{}\mathrel{=}{}\<[17E]%
\>[20]{}\Varid{check}\;\Varid{n}\;(\Varid{prop}_{1.1}\;\Varid{f}\mathbin{::}\Conid{Set}\;\Int\,\to\,\Conid{Bool}){}\<[E]%
\\
\>[3]{}\Varid{test}_{1.2}\;{}\<[11]%
\>[11]{}\Varid{n}\;{}\<[14]%
\>[14]{}\Varid{f}{}\<[17]%
\>[17]{}\mathrel{=}{}\<[17E]%
\>[20]{}\Varid{check}\;\Varid{n}\;(\Varid{prop}_{1.2}\;\Varid{f}\mathbin{::}\Conid{Set}\;\Int\,\to\,\Conid{Bool}){}\<[E]%
\\
\>[3]{}\Varid{test}_{1.3}\;{}\<[11]%
\>[11]{}\Varid{n}\;{}\<[14]%
\>[14]{}\Varid{f}{}\<[17]%
\>[17]{}\mathrel{=}{}\<[17E]%
\>[20]{}\Varid{check}\;\Varid{n}\;(\Varid{prop}_{1.3}\;\Varid{f}\mathbin{::}\Conid{Set}\;\Int\,\to\,\Conid{Bool}){}\<[E]%
\ColumnHook
\end{hscode}\resethooks
\noindent
For clarity of exposition, we neglect efficiency concerns.
The example implementations of \ensuremath{\Varid{min}}, \ensuremath{\Varid{argmin}} and most of the
programs and tests discussed in the paper can be implemented more
efficiently, e.g. by computing them in one go or by storing
intermediate computations.

In climate science, economics, engineering, and related fields, however,
the function to be minimized is typically \emph{not} known with
certainty.
We will explore two kinds of uncertainty: \emph{value-based}
and \emph{functorial} uncertainty.

\paragraph*{Value-based uncertainty.}
In many decision problems, one cannot meaningfully compare different
decisions using a single \emph{totally ordered} objective value.
In climate policy, for example, it would be all too easy to reduce
greenhouse gas emissions in line with the IPCC recommendations for the
next decades \citep{IPCC2018_SPM} by mass destruction, starvation
or unacceptably low global welfare levels.
In reality, non-trivial decisions almost inevitably entail
trade-offs and compromises between different, often incompatible
objectives.
The outcomes cannot be ranked along a single total order.

For example, sustainable climate policies must keep the earth system
within a ``safe'' state space constrained by a number of planetary
boundaries \citep{ro06010m,
  doi:10.1126/sciadv.adh2458}.
%
Respecting such boundaries is often modeled not as a hard feasibility
constraint, but instead by introducing additional objectives that
quantify proximity to, or violation of, each boundary.
This turns policy design into a genuinely multi-objective problem in
which several, only partially comparable criteria must be minimized
simultaneously.

Even for decision problems involving only two objectives, it is often
difficult \citep{shape2021} to compare values in different dimensions:
how many dollars is it worth being one hour faster on a 6-hour car
trip?
How many human lives should be traded off against avoiding a meltdown
of Greenland's ice sheet?

In economics, the exchange rate between distinct objectives---for example,
the number of human lives equivalent to preserving a unit area of
the Greenland ice sheet---is termed a \emph{price}.
Virtually all integrated assessment models for policy advice in
matters of climate policy \citep{Nordhaus2018, Nordhaus12261} are
based on a system of prices.
By assuming a system of prices and aggregating incomparable measures
into a single number (colloquially called ``cost''), it looks like one
can thus transform a multi-objective problem into a much simpler
single-objective problem.
But that hides the complexity of setting the prices.

For many multi-objective decision problems, it is
better to recognize \citep{pindyck2017} that the codomain of \ensuremath{\Varid{f}} is in
fact a tuple and only supports a partial \emph{product} order instead
of a total order.
In practice, this means recognizing that one is facing a
\emph{multi-objective} optimization (MOO) problem and that the best that
one can do is to provide decision makers with a set of its \emph{Pareto
optimal} solutions \citep{CARLINO202016593}.

Informally, Pareto optimal solutions are those that cannot be
improved: it is not possible to find values in the domain of the
functions to be minimized that obtain a better result in one objective
without sacrificing another.
We define multi-objective optimization problems and Pareto optimality
in \cref{section:moo}.

\paragraph*{Functorial uncertainty.}

The second kind of uncertainty we address is \emph{functorial}
uncertainty.
In many decision problems in engineering, economics, and climate
science, we are faced with decisions whose consequences are not
precisely known.
For example, we do not know exactly how much emission reductions
will result from enforcing a specific carbon tax or a specific trading
scheme for emission rights.
And even for problems in which the consequences of decisions do not
depend on human behavior, our knowledge is most often uncertain
\citep{Moreno-Cruz+Keith2012, Zarnetskee1921854118,
  NealNewellPitman2025}.
Our understanding of the way complex systems respond to inputs is
based on physical experiments, mathematical models, and numerical
approximations.
These are all affected by uncertainties.

At the semantic level, it is useful to distinguish between
\emph{epistemic} and \emph{aleatoric} uncertainty
\citep{shepherd+al2018, shepherd2019}.
In the epistemic case, we do not know precisely the values of
parameters, initial conditions, boundary conditions or other features
of the experiment or simulation that defines the system at stake.
In contrast, in the aleatoric case the system itself is
non-deterministic (even if we would know all input parameters).

In both cases, however, a decision problem must handle functions that
return values in \ensuremath{\Conid{U}\;\Conid{B}}, where \ensuremath{\Conid{U}} is an uncertainty functor, for
example, \ensuremath{\Conid{U}\mathrel{=}\Conid{List}} or \ensuremath{\Conid{U}\mathrel{=}\Conid{SimpleProb}}, see
\citep{10.1017/S0956796805005721, ionescu2009,
  2017_Botta_Jansson_Ionescu}.
Finding all the elements of a finite set at which \ensuremath{\Varid{f}} attains its
minimum requires comparing values in \ensuremath{\Conid{U}\;\Conid{B}} rather than in the base
type \ensuremath{\Conid{B}}, and many different order relations are possible on \ensuremath{\Conid{U}\;\Conid{B}}.

A well-established approach for this is through a \emph{measure}
function that maps values in \ensuremath{\Conid{U}\;\Conid{B}} into a totally ordered set \ensuremath{\Conid{C}} and
fulfills a natural monotonicity condition, see \citep{ionescu2009,
  IONESCU_2016b}.
In many optimization problems, the uncertainty is stochastic, so \ensuremath{\Conid{B}\mathrel{=}\Conid{C}\mathrel{=}\Real}, and a popular measure is the \emph{expected value}.
This is consistent with an attitude towards risk that is neither risk
prone nor risk averse.
Sometimes decision makers want to take decisions that minimize the
risk of worst-case outcomes.
(In the context of stochastic uncertainty, this is often the case when
decision makers are faced with low-probability, high-impact possible
outcomes.)
In these cases the expected value measure is not the right choice.

Measure functions are typically non-injective.
Even for a small set \ensuremath{\Conid{B}}, there can be infinitely many probability
distributions that obtain the same expected value.
This implies that optimization under functorial uncertainty typically
returns a set of values of type \ensuremath{\Conid{U}\;\Conid{B}} rather than a single value of
type \ensuremath{\Conid{B}}.
We discuss how to specify and solve optimization problems under
functorial uncertainty in \cref{section:monadic}.

\paragraph*{Overview and motivation.}

We specify, implement, and test generic methods for solving optimization
problems under value uncertainty and functorial uncertainty in
\cref{section:moo,section:monadic}.
In this manuscript we do not present ``real world'' applications, but
in \cref{subsection:applications} we discuss and demonstrate the
methods of \cref{section:moo} on standard multi-objective optimization
benchmark problems using evolutionary techniques.
We present related work in \cref{section:related} and discuss how to
apply methods for optimization under functorial uncertainty in
\cref{section:conclusion}.

The main motivation for our contribution is pedagogical:
we show how functional programming can support the precise
specification of optimization problems under uncertainty and the
systematic testing of software components against these
specifications.
This provides a principled way of ``understanding orders''---that is,
clarifying which order relations are appropriate in different
uncertainty settings---while also illustrating how specifications can
guide the implementation and validation of correct programs.
In dependently typed systems such as Rocq \citep{CoqProofAssistant},
Agda \citep{norell2007thesis}, and Idris \citep{idrisreference}, these
specifications could further be used to formally verify the
corresponding components, thereby ensuring a higher level of software
correctness \citep{ionescujansson:LIPIcs:2013:3899}.

Equally important, the specifications and tests are intended to help
students and practitioners develop a deeper understanding of
optimization under uncertainty and to approach such problems in a
systematic and disciplined way.
We hope that modelers and domain experts will benefit from the
discussion of monotonicity conditions for optimization under
functorial uncertainty in \cref{subsection:monotonicity}, and that the
examples in \cref{subsection:ufexamples} will help them avoid the
pitfalls of seemingly natural but inconsistent measure functions.

Before we begin, a caveat: as already noted, where appropriate we
trade efficiency for clarity.
While the methods and tests presented in \cref{section:moo,section:monadic}
can in fact be applied to ``real-world'' problems
\citep{Pusztai2023BayesOptMMI, jansson+2025} and to the testing of
more advanced approaches, their primary purpose here is
to make the central ideas clear.

\section{Value-based uncertainty and multi-objective optimization}
\label{section:moo}

As discussed in the introduction, many decision problems cannot be
formulated using a single objective function with values in a totally
ordered codomain.
Instead, decision quality is assessed along several dimensions that
cannot be collapsed into a single total order without loss of
information.
From a formal perspective, this corresponds to replacing a single
objective function \ensuremath{\Varid{f}\ \mathop{:}\ \Conid{A}\,\to\,\Conid{B}}, where \ensuremath{\Conid{B}} is totally ordered, by
\ensuremath{\Varid{n}} functions, or equivalently, returning a vector equipped with a partial order.
For the sake of readability and without loss of generality, we
consider the case in which the \ensuremath{\Varid{n}} objective functions
\ensuremath{f_1,\mathbin{...},f_n} share the same codomain \ensuremath{\Conid{B}}.
They can then be represented by a single function \ensuremath{\Varid{fs}\ \mathop{:}\ \Conid{A}\,\to\,B^n},
with the product order induced by the order on \ensuremath{\Conid{B}}.

In single-objective minimization problems, the image of an element \ensuremath{\Varid{a}\in \Varid{as}} under \ensuremath{\Varid{f}\ \mathop{:}\ \Conid{A}\,\to\,\Conid{B}} is either equal to the minimum of \ensuremath{\Varid{map}\;\Varid{f}\;\Varid{as}}
or is strictly greater (and hence dominated by) that minimum.
The totality of the order ensures that exactly one of these cases
applies.

In multi-objective minimization (MOO) problems, the image of \ensuremath{\Varid{a}} under
\ensuremath{\Varid{fs}} is an element of \ensuremath{B^n}, where the product order is typically only
partial.
Here too, we have to distinguish between two cases: the one in which
\ensuremath{\Varid{fs}\;\Varid{a}} is an element of the \emph{Pareto minimum} front of \ensuremath{\Varid{map}\;\Varid{fs}\;\Varid{as}} and the case in which \ensuremath{\Varid{fs}\;\Varid{a}} is behind (dominated by) that
front.
In the following subsections we make these notions precise and discuss
their computational treatment.

\subsection{Pareto optimality}
\label{subsection:pareto}

We start from a \emph{dominance} (better than) relation \ensuremath{(\mathrel{\prec})} on \ensuremath{B^n}.
Dominance is anti-reflexive, transitive but usually not total.
For example, in \cref{fig:Front} both $q_1$ and $q_2$ dominate $q_3$
(are better than $q_3$ for minimization) but neither dominates the
other: $q_1$ and $q_2$ are mutually indifferent.
A set is called indifferent (or an \emph{antichain}) with respect to \ensuremath{(\mathrel{\prec})}
iff for any two distinct elements \ensuremath{\Varid{x}} and \ensuremath{\Varid{x'}} in the set neither \ensuremath{\Varid{x}\mathrel{\prec}\Varid{x'}} nor \ensuremath{\Varid{x'}\mathrel{\prec}\Varid{x}}.
%
%
\DONE{Tim: \ensuremath{\mathrel{\prec}} is
  not defined. At least two natural definitions would fit the graphics:
  \ensuremath{\Varid{x}\mathrel{\prec}\Varid{y}\Leftrightarrow x_{1}\mathbin{<}y_{1}\wedge x_{2}\mathbin{<}y_{2}} or \ensuremath{\Varid{x}\mathrel{\prec}\Varid{y}\Leftrightarrow x_{1}\mathrel{\leq}y_{1}\wedge x_{2}\mathrel{\leq}y_{2}\wedge (x_{1}\mathbin{<}y_{1}\vee x_{2}\mathbin{<}y_{2})} where \ensuremath{\mathbin{<}} and \ensuremath{\mathrel{\leq}} are the usual relations on
  reals.  Probably the latter is meant. Then it would also be good to
  illustrate in the graphic that e.g. \ensuremath{\Varid{q\char95 1}} could have the same abscissa
  as \ensuremath{p_{2}} and would still be dominated by it. Nicola: That's fine, we
  say ``for \ensuremath{\mathrel{\prec}} in the context, we define \ensuremath{\mathrel{\prec}} in the next section''.}
\definecolor{controlcolor}{RGB}{200,47,68}
\definecolor{expcolor}{RGB}{30,180,200}
\definecolor{objcolor}{RGB}{100,160,60}
\begin{figure}[htp]
\begin{center}
\begin{tikzpicture}
[scale=0.9,axes1/.style={scale=1.2}]
\begin{scope}[axes1]
  \draw[help lines,step=0.5cm,opacity=0.5] (-2,-1.5) grid (5,3);

  \draw[->] (-2.25,0) -- (5.25,0) node[right] {$x_1$} coordinate(x1 axis);
  \draw[->] (0,-1.75) -- (0,3.25) node[above] {$x_2$} coordinate(x2 axis);

  \foreach \x/\xtext in {-2,-1,1,2,3,4}
    \draw[xshift=\x cm] (0pt,2pt) -- (0pt,-2pt) node[below] {$\xtext$};

  \foreach \y/\ytext in {-1,1,2}
    \draw[yshift=\y cm] (2pt,0pt) -- (-2pt,0pt) node[left] {$\ytext$};

  \fill (-1,2.5)    circle [radius=2pt] node [below] {$p_1$};
  \fill [controlcolor,opacity=0.25] (-1,2.5) rectangle (5.0,3.0);

  \fill (1.0,0.75)  circle [radius=2pt] node [below] {$p_2$};
  \fill [controlcolor,opacity=0.25] (1,0.75) rectangle (5.0,3.0);

  \fill (1.5,-0.5)  circle [radius=2pt] node [below] {$p_3$};
  \fill [controlcolor,opacity=0.25] (1.5,-0.5) rectangle (5.0,3.0);

  \fill (3.5,-1)    circle [radius=2pt] node [below] {$p_4$};
  \fill [controlcolor,opacity=0.25] (3.5,-1) rectangle (5.0,3.0);

  \draw (1.0,1.5)   circle [radius=2pt] node [right] {$q_1$};
  \fill [controlcolor,opacity=0.25] (1.0,1.5) rectangle (5.0,3.0);

  \draw (2.0,0.5)   circle [radius=2pt] node [right] {$q_2$};
  \fill [controlcolor,opacity=0.25] (2.0,0.5) rectangle (5.0,3.0);

  \draw (2.5,2.0)   circle [radius=2pt] node [right] {$q_3$};
  \fill [controlcolor,opacity=0.25] (2.5,2.0) rectangle (5.0,3.0);

\end{scope}
\end{tikzpicture}
\end{center}
\caption{Dominance, indifference and Pareto front ($p$-points) in
  $\mathbb{R}^2$: $p_2$ dominates $q_1$ and $q_3$; $p_3$
  dominates $q_2$ and $q_3$; $q_1$ is indifferent to $p_1$ and $p_1$ is
  indifferent to $q_3$ but $q_1$ dominates $q_3$. The shaded area
  associated with each point represents the subset of $\mathbb{R}^2$
  dominated by that point in the squared rectangle.
  The Pareto front of \ensuremath{\{\mskip1.5mu \Varid{p}_{1},\Varid{p}_{2},\Varid{p}_{3},\Varid{p}_{4},q_1,q_2,q_3\mskip1.5mu\}} is \ensuremath{\{\mskip1.5mu \Varid{p}_{1},\Varid{p}_{2},\Varid{p}_{3},\Varid{p}_{4}\mskip1.5mu\}}.
  \label{fig:Front} }
\end{figure}

\noindent
A set \ensuremath{X_P} is then called a Pareto optimal subset (Pareto front)
of \ensuremath{\Conid{X}} w.r.t. \ensuremath{(\mathrel{\prec})} iff
\DONE{Patrik: Talk about ``a'' or ``the'' Pareto front.}
\begin{itemize}
\item \ensuremath{X_P\subseteq \Conid{X}},
\item \ensuremath{X_P} is indifferent with respect to \ensuremath{(\mathrel{\prec})},
\item \ensuremath{X_P} dominates \ensuremath{\Conid{X}\,\backslash\,X_P} that is \ensuremath{\forall \Varid{x'}\in \Conid{X}\,\backslash\,X_P.~\exists \Varid{x}\in X_P.~\Varid{x}\mathrel{\prec}\Varid{x'}}.
\end{itemize}
\DONE{Patrik: Some notation for "\ensuremath{\Conid{A}} dominates \ensuremath{\Conid{B}\mathbin{-}\Conid{A}}"? To make this
  easier to read. I suggest two dots over the element ordering
  relation.}
\DONE{Patrik: Perhaps cite \citep{jansson_jansson_2023} (section
  3.2 Thinning). See Related Work.}
\DONE{Patrik: conditions 1 and 3 together specify a collection of
  dominating subsets (including the full set \ensuremath{\Conid{X}}). This corresponds to
  the relation \ensuremath{\Conid{Thin}\;X_P\;\Conid{X}} in \citep{jansson_jansson_2023}. Condition
  2 requires it to be as small as possible.}
%
\DONE{Tim: as can easily be shown, \ensuremath{X_P} satisfies these conditions iff
  it is the set of minimal elements in \ensuremath{\Conid{X}} w.r.t. \ensuremath{\mathrel{\prec}}, so the rather
  complicated specification here is unnecessary.
  Nicola: it is not clear to me that saying that \ensuremath{X_P} is the set of
  minimal elements in \ensuremath{\Conid{X}} is simpler than saying that \ensuremath{X_P} is the set
  of indifferent elements in \ensuremath{\Conid{X}}. But we should point the reader to
  equivalent formulations of Pareto efficiency if we think that these
  can help to better grasp the notion. In the main text, we should use
  the formulation that is more helpful for explaining/understanding
  the computation of \ensuremath{\Conid{ParetoOpt}} as a fold with \ensuremath{[\mskip1.5mu \mskip1.5mu]} and \ensuremath{\Varid{bump}}.}
%
The third condition in the definition of Pareto optimality expresses
that elements on the front dominate the rest.
Written out directly, this condition simultaneously involves set
difference, universal and existential quantification.
To make its structure explicit and to prepare for later algorithmic
use, we introduce some auxiliary notions that lift the dominance
relation \ensuremath{(\mathrel{\prec})} from elements to sets in two steps.
A set dominates an element if it contains a strictly better one, and
a set dominates another set if it dominates all of its elements
%
\begin{hscode}\SaveRestoreHook
\column{B}{@{}>{\hspre}l<{\hspost}@{}}%
\column{E}{@{}>{\hspre}l<{\hspost}@{}}%
\>[B]{}(\mathrel{\dot{\prec}})\mathbin{::}\Conid{Set}\;\Varid{a}\,\to\,\Varid{a}\,\to\,\Conid{Bool}{}\<[E]%
\\
\>[B]{}\Varid{xs}\mathrel{\dot{\prec}}\Varid{x'}\mathrel{=}\Varid{any}\;(\prec\Varid{x'})\;\Varid{xs}{}\<[E]%
\\[\blanklineskip]%
\>[B]{}(\mathrel{\ddot{\prec}})\mathbin{::}\Conid{Set}\;\Varid{a}\,\to\,\Conid{Set}\;\Varid{a}\,\to\,\Conid{Bool}{}\<[E]%
\\
\>[B]{}\Varid{xs}\mathrel{\ddot{\prec}}\Varid{xs'}\mathrel{=}\Varid{all}\;(\Varid{xs}\mathrel{\dot{\prec}})\;\Varid{xs'}{}\<[E]%
\ColumnHook
\end{hscode}\resethooks
With these in place, the third condition becomes ``\ensuremath{X_P} dominates the
rest of \ensuremath{\Conid{X}}'', that is \ensuremath{X_P\mathrel{\ddot{\prec}}(\Conid{X}\,\backslash\,X_P)}.
It turns out that \ensuremath{X_P} is unique: it is precisely the set of minimal
elements of \ensuremath{\Conid{X}} with respect to \ensuremath{(\prec)}.
We prefer the formulation above because it separates indifference from
coverage of the dominated region, which mirrors the structure of
algorithms that compute Pareto fronts incrementally.

Perhaps not surprisingly, computing the Pareto optimal subset of a
finite set of values of type \ensuremath{B^n} is almost straightforward.
The problem (for \ensuremath{(\prec)} in the context, we define \ensuremath{(\prec)} in the next
section) is to implement a function
\begin{hscode}\SaveRestoreHook
\column{B}{@{}>{\hspre}l<{\hspost}@{}}%
\column{3}{@{}>{\hspre}l<{\hspost}@{}}%
\column{E}{@{}>{\hspre}l<{\hspost}@{}}%
\>[3]{}\Varid{paretoOpt}\mathbin{::}(\Conid{NDim}\;\Varid{n},\Conid{Ord}\;\Varid{b})\Rightarrow \Conid{Set}\;(\Varid{n}\;\Varid{b})\,\to\,\Conid{Set}\;(\Varid{n}\;\Varid{b}){}\<[E]%
\ColumnHook
\end{hscode}\resethooks
that fulfills

\begin{enumerate}
\item \ensuremath{\Varid{paretoOpt}\;\Varid{xs}} \ is a subset of \ \ensuremath{\Varid{xs}},
\item \ensuremath{\Varid{paretoOpt}\;\Varid{xs}} \ is indifferent with respect to \ensuremath{(\prec)},
\item \ensuremath{\Varid{paretoOpt}\;\Varid{xs}} \ dominates \ \ensuremath{\Varid{xs}\mathbin{\char92 \char92 }\Varid{paretoOpt}\;\Varid{xs}}.
\end{enumerate}
%
for any \ensuremath{\Varid{xs}} of suitable type.
The type class \ensuremath{\Conid{NDim}} in the context of \ensuremath{\Varid{paretoOpt}} encodes the notion
of a fixed-length sequence.
Thus, values of type \ensuremath{\Varid{n}\;\Varid{b}} are vectors of \ensuremath{\Varid{b}}-values of the same length.
The difference in the third requirement is set difference.
%
For finite \ensuremath{\Varid{xs}}, all three requirements are decidable and can be
encoded by a generic predicate
\begin{hscode}\SaveRestoreHook
\column{B}{@{}>{\hspre}l<{\hspost}@{}}%
\column{3}{@{}>{\hspre}l<{\hspost}@{}}%
\column{29}{@{}>{\hspre}l<{\hspost}@{}}%
\column{E}{@{}>{\hspre}l<{\hspost}@{}}%
\>[3]{}\Varid{isParetoOptOf}\mathbin{::}(\Conid{NDim}\;\Varid{n},\Conid{Ord}\;\Varid{b})\Rightarrow \Conid{Set}\;(\Varid{n}\;\Varid{b})\,\to\,\Conid{Set}\;(\Varid{n}\;\Varid{b})\,\to\,\Conid{Bool}{}\<[E]%
\\[\blanklineskip]%
\>[3]{}\Varid{pxs}\mathbin{`\Varid{isParetoOptOf}`}\Varid{xs}\mathrel{=}{}\<[29]%
\>[29]{}(\Varid{pxs}\mathbin{`\Varid{isSubset}`}\Varid{xs})\mathrel{\wedge}(\Varid{isIndiff}\;\Varid{pxs})\mathrel{\wedge}(\Varid{pxs}\mathrel{\ddot{\prec}}(\Varid{xs}\mathbin{\char92 \char92 }\Varid{pxs})){}\<[E]%
\ColumnHook
\end{hscode}\resethooks
with straightforward definitions of \ensuremath{\Varid{isSubset}}, \ensuremath{\Varid{isIndiff}} and
\ensuremath{(\mathrel{\ddot{\prec}})}, see the literate Haskell code that generates this
manuscript \citep{MOOcodeRepo2025a}.
As in the computation of the minimum from the introduction, the idea
is to compute the Pareto front as a fold
\begin{hscode}\SaveRestoreHook
\column{B}{@{}>{\hspre}l<{\hspost}@{}}%
\column{3}{@{}>{\hspre}l<{\hspost}@{}}%
\column{18}{@{}>{\hspre}l<{\hspost}@{}}%
\column{23}{@{}>{\hspre}c<{\hspost}@{}}%
\column{23E}{@{}l@{}}%
\column{26}{@{}>{\hspre}l<{\hspost}@{}}%
\column{E}{@{}>{\hspre}l<{\hspost}@{}}%
\>[3]{}\Varid{paretoOpt}\;{}\<[18]%
\>[18]{}[\mskip1.5mu \mskip1.5mu]{}\<[23]%
\>[23]{}\mathrel{=}{}\<[23E]%
\>[26]{}[\mskip1.5mu \mskip1.5mu]{}\<[E]%
\\[\blanklineskip]%
\>[3]{}\Varid{paretoOpt}\;(\Varid{x}\ \mathop{:}\ \Varid{xs}){}\<[23]%
\>[23]{}\mathrel{=}{}\<[23E]%
\>[26]{}\Varid{bump}\;\Varid{x}\;(\Varid{paretoOpt}\;\Varid{xs}){}\<[E]%
\ColumnHook
\end{hscode}\resethooks
with a function \ensuremath{\Varid{bump}} that fulfills three requirements

\begin{enumerate}[resume]
\item \ensuremath{\Varid{bump}\;\Varid{x}\;\Varid{xs}} is a subset of \ensuremath{\Varid{x}\ \mathop{:}\ \Varid{xs}},
\item \ensuremath{\Varid{bump}\;\Varid{x}} preserves indifference: \ensuremath{\Varid{isIndiff}\;\Varid{xs}\Rightarrow \Varid{isIndiff}\;(\Varid{bump}\;\Varid{x}\;\Varid{xs})},
\item If \ensuremath{\Varid{xs}} dominates \ensuremath{\Varid{ys}} then \ensuremath{\Varid{bump}\;\Varid{x}\;\Varid{xs}} dominates \ensuremath{\Varid{ys}}.
\end{enumerate}
\DONE{Patrik: Check with presentation of thinning (compare with Level-p-complexity paper and its citation of Jeremy Gibbons): ``A good (but abstract) reference for thinning is the Algebra of Programming book \cite[Chapter 8]{bird_algebra_1997} and more concrete references are the corresponding developments in Agda \citep{DBLP:journals/jfp/MuKJ09} and Haskell \citep{bird_gibbons_2020}.''}
\DONE{Check if the last condition is correctly stated, or if there should use \ensuremath{\Varid{xp}\mathbin{`\Varid{rel23}`}\Varid{x}\mathrel{=}\Varid{xp}\mathrel{\ddot{\prec}}(\Varid{xp}\mathbin{\char92 \char92 }\Varid{x})}. Basically, I think we could make the right set should also grow to make the property stronger.}

\noindent
for any \ensuremath{\Varid{x}}, \ensuremath{\Varid{xs}}, \ensuremath{\Varid{ys}} of matching types. The implementation of \ensuremath{\Varid{bump}}
is almost straightforward
\DONE{Perhaps use type synonym PSet for ParetoSet to visually indicate the invariant. Patrik: not done to avoid too much explaining.}
\begin{hscode}\SaveRestoreHook
\column{B}{@{}>{\hspre}l<{\hspost}@{}}%
\column{3}{@{}>{\hspre}l<{\hspost}@{}}%
\column{9}{@{}>{\hspre}l<{\hspost}@{}}%
\column{15}{@{}>{\hspre}c<{\hspost}@{}}%
\column{15E}{@{}l@{}}%
\column{18}{@{}>{\hspre}l<{\hspost}@{}}%
\column{21}{@{}>{\hspre}l<{\hspost}@{}}%
\column{53}{@{}>{\hspre}c<{\hspost}@{}}%
\column{53E}{@{}l@{}}%
\column{56}{@{}>{\hspre}l<{\hspost}@{}}%
\column{59}{@{}>{\hspre}l<{\hspost}@{}}%
\column{80}{@{}>{\hspre}l<{\hspost}@{}}%
\column{E}{@{}>{\hspre}l<{\hspost}@{}}%
\>[3]{}\Varid{bump}\mathbin{::}(\Conid{NDim}\;\Varid{n},\Conid{Ord}\;\Varid{b})\Rightarrow \Varid{n}\;\Varid{b}\,\to\,\Conid{Set}\;(\Varid{n}\;\Varid{b})\,\to\,\Conid{Set}\;(\Varid{n}\;\Varid{b}){}\<[E]%
\\[\blanklineskip]%
\>[3]{}\Varid{bump}\;{}\<[9]%
\>[9]{}\Varid{x}\;[\mskip1.5mu \mskip1.5mu]{}\<[15]%
\>[15]{}\mathrel{=}{}\<[15E]%
\>[18]{}\{\mskip1.5mu \Varid{x}\mskip1.5mu\}{}\<[E]%
\\[\blanklineskip]%
\>[3]{}\Varid{bump}\;{}\<[9]%
\>[9]{}\Varid{x}\;(\Varid{p}\ \mathop{:}\ \Varid{ps}){}\<[21]%
\>[21]{}\mid \Varid{p}\doubleequals\Varid{x}\mathrel{\vee}\Varid{p}\prec\Varid{x}\quad {}\<[53]%
\>[53]{}\mathrel{=}{}\<[53E]%
\>[56]{}\Varid{p}{}\<[59]%
\>[59]{}\ \mathop{:}\ \Varid{ps}{}\<[80]%
\>[80]{}\mbox{\onelinecomment  case 1}{}\<[E]%
\\
\>[21]{}\mid \Varid{x}\prec\Varid{p}{}\<[53]%
\>[53]{}\mathrel{=}{}\<[53E]%
\>[56]{}\Varid{x}{}\<[59]%
\>[59]{}\ \mathop{:}\ \Varid{remove}\;(\Varid{x}\prec)\;\Varid{ps}{}\<[80]%
\>[80]{}\mbox{\onelinecomment  case 2}{}\<[E]%
\\
\>[21]{}\mid \Varid{otherwise}{}\<[53]%
\>[53]{}\mathrel{=}{}\<[53E]%
\>[56]{}\Varid{p}{}\<[59]%
\>[59]{}\ \mathop{:}\ \Varid{bump}\;\Varid{x}\;\Varid{ps}{}\<[80]%
\>[80]{}\mbox{\onelinecomment  case 3}{}\<[E]%
\ColumnHook
\end{hscode}\resethooks
\DONE{Perhaps combine \ensuremath{\Varid{p}\doubleequals\Varid{x}\mid \Varid{p}\prec\Varid{x}} into an infix ``less than or equal'' operator above. (Not done.)}
\DONE{Then the three cases below are correct (now, part of case 1 is ``missing'').}
\DONE{Perhaps add notation for indifference - perhaps ||. (Not done.)}
\noindent
Beside the trivial case \ensuremath{\Varid{p}\doubleequals\Varid{x}}, the induction step distinguishes
between three (sub-)cases: \ensuremath{\Varid{p}\prec\Varid{x}}, \ensuremath{\Varid{x}\prec\Varid{p}}, or neither dominates the
other.
\DONE{Patrik: Perhaps use other (more set-like) notation for \ensuremath{\{\mskip1.5mu \Varid{x}\mskip1.5mu\}}.}
In case 1, the element \ensuremath{\Varid{p}} on the front dominates the new element \ensuremath{\Varid{x}},
thus the front is unchanged.
In case 2, the new element \ensuremath{\Varid{x}} replaces \ensuremath{\Varid{p}}, and is also given the
chance to remove more elements from the front.
Finally, in case 3, we keep \ensuremath{\Varid{p}} and recurse to compare \ensuremath{\Varid{x}} to the rest
of the front.

As for the problem (of finding all the elements of a finite set at which
a function attains its minimum) discussed in the introduction,
2.4-2.6 are decidable properties and one can test that \ensuremath{\Varid{bump}}
fulfills these requirements.
We implement first the polymorphic properties, then the monomorphic
tests:
\begin{hscode}\SaveRestoreHook
\column{B}{@{}>{\hspre}l<{\hspost}@{}}%
\column{3}{@{}>{\hspre}l<{\hspost}@{}}%
\column{11}{@{}>{\hspre}l<{\hspost}@{}}%
\column{14}{@{}>{\hspre}c<{\hspost}@{}}%
\column{14E}{@{}l@{}}%
\column{15}{@{}>{\hspre}l<{\hspost}@{}}%
\column{17}{@{}>{\hspre}l<{\hspost}@{}}%
\column{18}{@{}>{\hspre}l<{\hspost}@{}}%
\column{22}{@{}>{\hspre}l<{\hspost}@{}}%
\column{34}{@{}>{\hspre}l<{\hspost}@{}}%
\column{38}{@{}>{\hspre}l<{\hspost}@{}}%
\column{44}{@{}>{\hspre}l<{\hspost}@{}}%
\column{E}{@{}>{\hspre}l<{\hspost}@{}}%
\>[3]{}\Varid{prop}_{2.4}\;{}\<[11]%
\>[11]{}\Varid{xs}\;{}\<[15]%
\>[15]{}\Varid{x}{}\<[22]%
\>[22]{}\mathrel{=}\Varid{bump}\;\Varid{x}\;\Varid{xs}\mathbin{`\Varid{isSubset}`}(\Varid{x}\ \mathop{:}\ \Varid{xs}){}\<[E]%
\\
\>[3]{}\Varid{prop}_{2.5}\;{}\<[11]%
\>[11]{}\Varid{xs}\;{}\<[15]%
\>[15]{}\Varid{x}{}\<[22]%
\>[22]{}\mathrel{=}\Varid{isIndiff}\;\Varid{xs}{}\<[38]%
\>[38]{}\Rightarrow \Varid{isIndiff}\;(\Varid{bump}\;\Varid{x}\;\Varid{xs}){}\<[E]%
\\
\>[3]{}\Varid{prop}_{2.6}\;{}\<[11]%
\>[11]{}\Varid{xs}\;{}\<[15]%
\>[15]{}\Varid{x}\;{}\<[18]%
\>[18]{}\Varid{ys}{}\<[22]%
\>[22]{}\mathrel{=}\Varid{xs}\mathrel{\ddot{\prec}}\Varid{ys}{}\<[44]%
\>[44]{}\Rightarrow (\Varid{bump}\;\Varid{x}\;\Varid{xs})\mathrel{\ddot{\prec}}\Varid{ys}{}\<[E]%
\\[\blanklineskip]%
\>[3]{}\Varid{test}_{2.4}\;{}\<[11]%
\>[11]{}\Varid{n}{}\<[14]%
\>[14]{}\mathrel{=}{}\<[14E]%
\>[17]{}\Varid{check}\;\Varid{n}\;(\Varid{prop}_{2.4}{}\<[34]%
\>[34]{}\mathbin{::}\Conid{Set}\;\Real^2\,\to\,\Real^2\,\to\,\Conid{Bool}){}\<[E]%
\\
\>[3]{}\Varid{test}_{2.5}\;{}\<[11]%
\>[11]{}\Varid{n}{}\<[14]%
\>[14]{}\mathrel{=}{}\<[14E]%
\>[17]{}\Varid{check}\;\Varid{n}\;(\Varid{prop}_{2.5}{}\<[34]%
\>[34]{}\mathbin{::}\Conid{Set}\;\Real^2\,\to\,\Real^2\,\to\,\Conid{Bool}){}\<[E]%
\\
\>[3]{}\Varid{test}_{2.6}\;{}\<[11]%
\>[11]{}\Varid{n}{}\<[14]%
\>[14]{}\mathrel{=}{}\<[14E]%
\>[17]{}\Varid{check}\;\Varid{n}\;(\Varid{prop}_{2.6}{}\<[34]%
\>[34]{}\mathbin{::}\Conid{Set}\;\Real^2\,\to\,\Real^2\,\to\,\Conid{Set}\;\Real^2\,\to\,\Conid{Bool}){}\<[E]%
\ColumnHook
\end{hscode}\resethooks
\DONE{Talk about uniqueness somewhere (``the'' or ``a'' Pareto optimal subset).}
With this \ensuremath{\Varid{bump}}, \ensuremath{\Varid{paretoOpt}} does compute the Pareto optimal subset
of a finite set:
\begin{hscode}\SaveRestoreHook
\column{B}{@{}>{\hspre}l<{\hspost}@{}}%
\column{3}{@{}>{\hspre}l<{\hspost}@{}}%
\column{E}{@{}>{\hspre}l<{\hspost}@{}}%
\>[3]{}\Varid{propParetoOpt}\;\Varid{xs}\mathrel{=}\Varid{paretoOpt}\;\Varid{xs}\mathbin{`\Varid{isParetoOptOf}`}\Varid{xs}{}\<[E]%
\\[\blanklineskip]%
\>[3]{}\Varid{testParetoOpt}\;\Varid{n}\mathrel{=}\Varid{check}\;\Varid{n}\;(\Varid{propParetoOpt}\mathbin{::}\Conid{Set}\;\Real^2\,\to\,\Conid{Bool}){}\<[E]%
\ColumnHook
\end{hscode}\resethooks
We conclude this section with three observations.

The first one is that, while our definition of \ensuremath{\Varid{paretoOpt}} uses a
sequential fold for simplicity, the underlying structure supports a
divide-and-conquer strategy.
One can define a \ensuremath{\Varid{merge}} operation such that \ensuremath{\Varid{paretoOpt}\;(\Varid{xs}\mathbin{+\!\!+}\Varid{ys})\mathrel{=}\Varid{merge}\;(\Varid{paretoOpt}\;\Varid{xs})\;(\Varid{paretoOpt}\;\Varid{ys})}.
This property paves the way for efficient parallel implementations,
where local Pareto fronts are computed on subsets of data and then
merged hierarchically.

The second observation is that indifference (with respect to
dominance) is a \emph{tolerance} relation: it is reflexive and
symmetric, but not transitive:
\begin{hscode}\SaveRestoreHook
\column{B}{@{}>{\hspre}l<{\hspost}@{}}%
\column{3}{@{}>{\hspre}l<{\hspost}@{}}%
\column{36}{@{}>{\hspre}l<{\hspost}@{}}%
\column{E}{@{}>{\hspre}l<{\hspost}@{}}%
\>[3]{}\Varid{testNotTransitive}\;\Varid{n}\mathrel{=}\Varid{falsify}\;\Varid{n}\;{}\<[36]%
\>[36]{}(\Varid{isTransitive}\;(\Varid{indiff}\;(\prec))\mathbin{::}\Int^2\,\to\,\Int^2\,\to\,\Int^2\,\to\,\Conid{Bool}){}\<[E]%
\ColumnHook
\end{hscode}\resethooks
The lack of transitivity necessitates the recursion in case 3.
Even though \ensuremath{\Varid{x}} is indifferent to the head \ensuremath{\Varid{p}}, it is not necessarily
indifferent to the elements in the tail \ensuremath{\Varid{ps}}.
Thus, simply prepending \ensuremath{\Varid{x}} would risk violating the invariant of
mutual indifference; we must instead traverse the rest of the list.
\DONE{Discuss \ensuremath{\Varid{bump}}, non-transitivity of indifference (test!)  and
  relate to thinning methods}
\DONE{See log.org: ``*** Preorder \ensuremath{\preccurlyeq}'' for some some notation + properties. (Not done - perhaps later.)}

The third observation is that \ensuremath{\Varid{bump}} plays a crucial role in
computational methods for approximating the Pareto front when \ensuremath{\Conid{B}\mathrel{=}\Real} and \ensuremath{\Varid{xs}} is a genuine subset of \ensuremath{\Real^n}.
We come back to this point in \cref{subsection:applications}.

\subsection{Multi-objective optimization}
\label{subsection:moo}

For \emph{minimization} problems, a tuple \ensuremath{\Varid{x}\in B^n} is better than \ensuremath{\Varid{x'}\in B^n} iff \ensuremath{\Varid{x}} is \emph{smaller or equal} to \ensuremath{\Varid{x'}} in all components and
\emph{strictly smaller} in at least one:
\DONE{See log.org: ** TODO Detail: alternative definition of the dominance relation. (Not done. Not important.)}
\DONE{Patrik: explain \ensuremath{\Varid{comps}}}
\begin{hscode}\SaveRestoreHook
\column{B}{@{}>{\hspre}l<{\hspost}@{}}%
\column{3}{@{}>{\hspre}l<{\hspost}@{}}%
\column{E}{@{}>{\hspre}l<{\hspost}@{}}%
\>[3]{}(\prec)\mathbin{::}(\Conid{NDim}\;\Varid{n},\Conid{Ord}\;\Varid{b})\Rightarrow \Varid{n}\;\Varid{b}\,\to\,\Varid{n}\;\Varid{b}\,\to\,\Conid{Bool}{}\<[E]%
\\[\blanklineskip]%
\>[3]{}\Varid{x}\prec\Varid{x'}\mathrel{=}\Varid{comps}\;\Varid{all}\;(\mathrel{\leq})\;\Varid{x}\;\Varid{x'}~\mathrel{\wedge}~\Varid{comps}\;\Varid{any}\;(\mathbin{<})\;\Varid{x}\;\Varid{x'}{}\<[E]%
\ColumnHook
\end{hscode}\resethooks
\noindent
The helper \ensuremath{\Varid{comps}} takes a list quantifier \ensuremath{\Varid{test}} (like \ensuremath{\Varid{all}} or
\ensuremath{\Varid{any}}), a binary relation \ensuremath{\Varid{op}} 
and two sequences \ensuremath{\Varid{xs}} and \ensuremath{\Varid{ys}} of values of the same length.
It applies the quantifier after zipping the two sequences with \ensuremath{\Varid{op}}:
\begin{hscode}\SaveRestoreHook
\column{B}{@{}>{\hspre}l<{\hspost}@{}}%
\column{3}{@{}>{\hspre}l<{\hspost}@{}}%
\column{E}{@{}>{\hspre}l<{\hspost}@{}}%
\>[3]{}\Varid{comps}\;\Varid{test}\;\Varid{op}\;\Varid{xs}\;\Varid{ys}\mathrel{=}\Varid{test}\;\Varid{id}\;(\Varid{toList}\;(\Varid{ndimZipWith}\;\Varid{op}\;\Varid{xs}\;\Varid{ys})){}\<[E]%
\ColumnHook
\end{hscode}\resethooks
We test that dominance is anti-reflexive, transitive and not total:
\begin{hscode}\SaveRestoreHook
\column{B}{@{}>{\hspre}l<{\hspost}@{}}%
\column{3}{@{}>{\hspre}l<{\hspost}@{}}%
\column{24}{@{}>{\hspre}l<{\hspost}@{}}%
\column{37}{@{}>{\hspre}l<{\hspost}@{}}%
\column{E}{@{}>{\hspre}l<{\hspost}@{}}%
\>[3]{}\Varid{testAntiReflexive}\;\Varid{n}{}\<[24]%
\>[24]{}\mathrel{=}\Varid{check}\;\Varid{n}\;{}\<[37]%
\>[37]{}(\Varid{isAntiReflexive}\;(\prec)\mathbin{::}\Int^2\,\to\,\Conid{Bool}){}\<[E]%
\\[\blanklineskip]%
\>[3]{}\Varid{testTransitive}\;\Varid{n}{}\<[24]%
\>[24]{}\mathrel{=}\Varid{check}\;\Varid{n}\;{}\<[37]%
\>[37]{}(\Varid{isTransitive}\;(\prec)\mathbin{::}\Int^2\,\to\,\Int^2\,\to\,\Int^2\,\to\,\Conid{Bool}){}\<[E]%
\\[\blanklineskip]%
\>[3]{}\Varid{testNotTotal}\;\Varid{n}{}\<[24]%
\>[24]{}\mathrel{=}\Varid{falsify}\;\Varid{n}\;{}\<[37]%
\>[37]{}(\Varid{isTotal}\;(\prec)\mathbin{::}\Int^2\,\to\,\Int^2\,\to\,\Conid{Bool}){}\<[E]%
\ColumnHook
\end{hscode}\resethooks
With \ensuremath{(\prec)} and \ensuremath{\Varid{paretoOpt}} in place, it is easy to implement a generic
solver for multi-objective minimization as done in the introduction for
the single-objective case:
\DONE{Patrik: why has the type changed from \ensuremath{\Varid{fs}\ \mathop{:}\ \Conid{A}\,\to\,B^n} in the
  first paragraph of this section, to \ensuremath{\Varid{fs}\ \mathop{:}\ \Varid{n}\;(\Conid{A}\,\to\,\Conid{B})} here?}
\begin{hscode}\SaveRestoreHook
\column{B}{@{}>{\hspre}l<{\hspost}@{}}%
\column{3}{@{}>{\hspre}l<{\hspost}@{}}%
\column{E}{@{}>{\hspre}l<{\hspost}@{}}%
\>[3]{}\Varid{paretoMin}\mathbin{::}(\Conid{NDim}\;\Varid{n},\Conid{Eq}\;\Varid{a},\Conid{Ord}\;\Varid{b})\Rightarrow (\Varid{a}\,\to\,\Varid{n}\;\Varid{b})\,\to\,\Conid{Set}\;\Varid{a}\,\to\,\Conid{Set}\;(\Varid{n}\;\Varid{b}){}\<[E]%
\\[\blanklineskip]%
\>[3]{}\Varid{paretoMin}\;\Varid{fs}\;\Varid{as}\mathrel{=}\Varid{paretoOpt}\;(\Varid{map}\;\Varid{fs}\;\Varid{as}){}\<[E]%
\ColumnHook
\end{hscode}\resethooks
\begin{hscode}\SaveRestoreHook
\column{B}{@{}>{\hspre}l<{\hspost}@{}}%
\column{3}{@{}>{\hspre}l<{\hspost}@{}}%
\column{27}{@{}>{\hspre}l<{\hspost}@{}}%
\column{E}{@{}>{\hspre}l<{\hspost}@{}}%
\>[3]{}\Varid{argparetoMin}\mathbin{::}(\Conid{NDim}\;\Varid{n},\Conid{Eq}\;\Varid{a},\Conid{Ord}\;\Varid{b})\Rightarrow (\Varid{a}\,\to\,\Varid{n}\;\Varid{b})\,\to\,\Conid{Set}\;\Varid{a}\,\to\,\Conid{Set}\;\Varid{a}{}\<[E]%
\\[\blanklineskip]%
\>[3]{}\Varid{argparetoMin}\;\Varid{fs}\;\Varid{as}\mathrel{=}{}\<[27]%
\>[27]{}\mathbf{let}\;\Varid{ps}\mathrel{=}\Varid{paretoMin}\;\Varid{fs}\;\Varid{as}\;\mathbf{in}\;\Varid{filter}\;(\lambda \Varid{a}\,\to\,\Varid{fs}\;\Varid{a}\in \Varid{ps})\;\Varid{as}{}\<[E]%
\ColumnHook
\end{hscode}\resethooks
%
%

Again, \ensuremath{\Varid{paretoMin}\;\Varid{fs}\;\Varid{as}} and \ensuremath{\Varid{argparetoMin}\;\Varid{fs}\;\Varid{as}} are in \ensuremath{\Varid{map}\;\Varid{fs}\;\Varid{as}}
and \ensuremath{\Varid{as}} respectively and we can take advantage of Quickcheck to test
that the two functions fulfil the specification (remember the
specification of \ensuremath{\Varid{min}}, \ensuremath{\Varid{argmin}} from \cref{section:intro})

\DONE{These properties would be useful to ``name'' in the code because
  they are reused inside the tests below. A bit like \ensuremath{\Varid{isTransitive}}
  etc. earlier: Haskell code which can be used separately from the
  specific tests. I have started with \ensuremath{\Varid{prop}_{2.7}}, and I think 1-2
  examples + explanation is enough, we don't need to strictly define
  all of them. For example, if a few of the earlier ones are explain,
  all these (\ensuremath{\Varid{prop}_{2.7}}, \ensuremath{\Varid{prop}_{2.8}}, \ensuremath{\Varid{prop}_{2.9}}) could be elided and just
  the \ensuremath{\Varid{test}_{2.7}}, etc.\ kept.}
\DONE{Perhaps let QuickCheck generate random functions. (Not done - perhaps later.)}
\DONE{Perhaps use non-injective functions (instead of \ensuremath{\Varid{id}} / \ensuremath{\Varid{fst}} / \ensuremath{\Varid{snd}}) - for example \ensuremath{(\mathbin{\Varid{`div`}}\mathrm{2})}. Not done. TODO later.}
\begin{enumerate}[resume]
\item \ensuremath{\Varid{paretoMin}\;\Varid{fs}\;\Varid{as}\mathbin{`\Varid{isParetoOptOf}`}\Varid{map}\;\Varid{fs}\;\Varid{as}},
\item \ensuremath{\forall \Varid{a}\in \Varid{as}}. \ \ensuremath{\Varid{a}\in \Varid{argparetoMin}\;\Varid{fs}\;\Varid{as}} \ \ensuremath{\Rightarrow } \ \ensuremath{\Varid{fs}\;\Varid{a}\in \Varid{paretoMin}\;\Varid{fs}\;\Varid{as}},
\item \ensuremath{\forall \Varid{a}\in \Varid{as}}. \ \ensuremath{\Varid{fs}\;\Varid{a}\in \Varid{paretoMin}\;\Varid{fs}\;\Varid{as}} \ \ensuremath{\Rightarrow } \ \ensuremath{\Varid{a}\in \Varid{argparetoMin}\;\Varid{fs}\;\Varid{as}}.
\end{enumerate}

\noindent
for all \ensuremath{\Varid{fs}}, \ensuremath{\Varid{as}} of suitable type.
\DONE{Check why the font size of the enumeration equation is different from that of the text code. (We use small in the haskell code - see [[file:main.lhs::\renewcommand{\hscodestyle}{\small\fixlengths}]].)}
We elide the definitions of \ensuremath{\Varid{prop}_{2.7}}, \ensuremath{\Varid{prop}_{2.8}}, and \ensuremath{\Varid{prop}_{2.9}} as they
are very similar to the specification clauses listed above but we
include the tests:

\DONE{Patrik: Explain or hide \ensuremath{\Conid{C}_2}.}
\DONE{Patrik: Explain the clever use of \ensuremath{\Varid{fst}} and \ensuremath{\Varid{snd}} for testing.}

\begin{hscode}\SaveRestoreHook
\column{B}{@{}>{\hspre}l<{\hspost}@{}}%
\column{3}{@{}>{\hspre}l<{\hspost}@{}}%
\column{11}{@{}>{\hspre}l<{\hspost}@{}}%
\column{14}{@{}>{\hspre}c<{\hspost}@{}}%
\column{14E}{@{}l@{}}%
\column{17}{@{}>{\hspre}l<{\hspost}@{}}%
\column{34}{@{}>{\hspre}l<{\hspost}@{}}%
\column{E}{@{}>{\hspre}l<{\hspost}@{}}%
\>[3]{}\Varid{test}_{2.7}\;{}\<[11]%
\>[11]{}\Varid{n}{}\<[14]%
\>[14]{}\mathrel{=}{}\<[14E]%
\>[17]{}\Varid{check}\;\Varid{n}\;(\Varid{prop}_{2.7}\;{}\<[34]%
\>[34]{}\Varid{id}\mathbin{::}\Conid{Set}\;(\Conid{Coord2D}\;\Int)\,\to\,\Conid{Bool}){}\<[E]%
\\
\>[3]{}\Varid{test}_{2.8}\;{}\<[11]%
\>[11]{}\Varid{n}{}\<[14]%
\>[14]{}\mathrel{=}{}\<[14E]%
\>[17]{}\Varid{check}\;\Varid{n}\;(\Varid{prop}_{2.8}\;{}\<[34]%
\>[34]{}\Varid{id}\mathbin{::}\Conid{Set}\;(\Conid{Coord2D}\;\Int)\,\to\,\Conid{Bool}){}\<[E]%
\\
\>[3]{}\Varid{test}_{2.9}\;{}\<[11]%
\>[11]{}\Varid{n}{}\<[14]%
\>[14]{}\mathrel{=}{}\<[14E]%
\>[17]{}\Varid{check}\;\Varid{n}\;(\Varid{prop}_{2.9}\;{}\<[34]%
\>[34]{}\Varid{id}\mathbin{::}\Conid{Set}\;(\Conid{Coord2D}\;\Int)\,\to\,\Conid{Bool}){}\<[E]%
\ColumnHook
\end{hscode}\resethooks
In the tests, we simulate a general objective function \ensuremath{\Varid{fs}\ \mathop{:}\ \Conid{A}\,\to\,B^n}
(with \ensuremath{\Varid{n}\mathrel{=}\mathrm{2}}) by generating random pairs directly and using \ensuremath{\Varid{fs}\mathrel{=}\Varid{id}\ \mathop{:}\ B^n\,\to\,B^n}.
This is equivalent to using the two functions \ensuremath{\Varid{fst}} and \ensuremath{\Varid{snd}} and we
are computing the Pareto minimum front of a (QuickCheck-generated) set
of random points in \ensuremath{\Int^2}.

In the definition of \ensuremath{\Varid{test}_{2.7-9}} and later throughout this paper, we use
a collection of instances of \ensuremath{\Conid{NDim}}: \ensuremath{\Conid{Coord1D}}, \ensuremath{\Conid{Coord2D}}, \ldots,
together with their respective data constructors \ensuremath{\Conid{C}_1}, \ensuremath{\Conid{C}_2}, \ldots, see
\citep{MOOcodeRepo2025a}.
\DONE{We should probably introduce all these CoordND in a box/fig somewhere.}

\subsection{Multi-objective and single-objective optimization}
\label{subsection:mosoo}


The functions \ensuremath{\Varid{paretoMin}} and \ensuremath{\Varid{argparetoMin}} generalize \ensuremath{\Varid{min}} and
\ensuremath{\Varid{argmin}} from the introduction. Thus, for a single objective function
\ensuremath{\Varid{fs}\mathrel{=}\Conid{C}_1\mathbin{\circ}\Varid{f}},
they should fulfil
\DONE{Patrik: I changed the formatting of \text{\ttfamily \char61{}\char61{}\char126{}} to be as \text{\ttfamily \char61{}\char61{}}, but perhaps something it should be more visible (stack a tilde on top?).}
\begin{enumerate}[resume]
\item \ensuremath{\Varid{paretoMin}\;(\Conid{C}_1\mathbin{\circ}\Varid{f})\;\Varid{as}\,\doubleequals\,\{\mskip1.5mu \Conid{C}_1\;(\Varid{min}\;\Varid{f}\;\Varid{as})\mskip1.5mu\}},
\item \ensuremath{\Varid{argparetoMin}\;(\Conid{C}_1\mathbin{\circ}\Varid{f})\;\Varid{as}\,\doubleequals\,\Varid{argmin}\;\Varid{f}\;\Varid{as}}.
\end{enumerate}

\noindent
for arbitrary, non empty \ensuremath{\Varid{as}}.
We can test these properties with
\DONE{I've hidden \ensuremath{\Varid{prop}_{2.10}}, \ensuremath{\Varid{prop}_{2.11}}. No comment about nonempty sets
  here because we have a top-level comment about it earlier.}

\begin{hscode}\SaveRestoreHook
\column{B}{@{}>{\hspre}l<{\hspost}@{}}%
\column{3}{@{}>{\hspre}l<{\hspost}@{}}%
\column{12}{@{}>{\hspre}l<{\hspost}@{}}%
\column{34}{@{}>{\hspre}l<{\hspost}@{}}%
\column{E}{@{}>{\hspre}l<{\hspost}@{}}%
\>[3]{}\Varid{testfun}\;(\Varid{x},\Varid{y})\mathrel{=}\Varid{x}\mathbin{+}\Varid{y}\mathbin{*}\Varid{x}{}\<[E]%
\\
\>[3]{}\Varid{test}_{2.10}\;{}\<[12]%
\>[12]{}\Varid{n}\mathrel{=}\Varid{check}\;\Varid{n}\;(\Varid{prop}_{2.10}\;{}\<[34]%
\>[34]{}\Varid{testfun}\mathbin{::}\Conid{Set}\;(\Int,\Int)\,\to\,\Conid{Bool}){}\<[E]%
\\
\>[3]{}\Varid{test}_{2.11}\;{}\<[12]%
\>[12]{}\Varid{n}\mathrel{=}\Varid{check}\;\Varid{n}\;(\Varid{prop}_{2.11}\;{}\<[34]%
\>[34]{}\Varid{testfun}\mathbin{::}\Conid{Set}\;(\Int,\Int)\,\to\,\Conid{Bool}){}\<[E]%
\ColumnHook
\end{hscode}\resethooks
A perhaps less obvious relationship between single-objective and
multi-objective optimization is that for \ensuremath{\Varid{f}\in \Varid{fs}}, the minimum of \ensuremath{\Varid{f}} on
\ensuremath{\Varid{as}} is in the image of \ensuremath{\Varid{argparetoMin}\;\Varid{fs}\;\Varid{as}} through \ensuremath{\Varid{f}} for any \ensuremath{\Varid{as}}

\begin{enumerate}[resume]
\item \ensuremath{\Varid{f}\in \Varid{fs}} \ \ensuremath{\Rightarrow } \ \ensuremath{\Varid{min}\;\Varid{f}\;\Varid{as}\in \Varid{map}\;\Varid{f}\;(\Varid{argparetoMin}\;\Varid{fs}\;\Varid{as})}
\end{enumerate}

\noindent
\DONE{Improve the \ensuremath{\Varid{fs}}, \ensuremath{\Varid{fs'}} usage + explain}
We can test this property with

\begin{hscode}\SaveRestoreHook
\column{B}{@{}>{\hspre}l<{\hspost}@{}}%
\column{3}{@{}>{\hspre}l<{\hspost}@{}}%
\column{16}{@{}>{\hspre}l<{\hspost}@{}}%
\column{21}{@{}>{\hspre}l<{\hspost}@{}}%
\column{28}{@{}>{\hspre}l<{\hspost}@{}}%
\column{E}{@{}>{\hspre}l<{\hspost}@{}}%
\>[3]{}\Varid{testfun3}\mathrel{=}{}\<[16]%
\>[16]{}\mathbf{let}\;{}\<[21]%
\>[21]{}\Varid{f}_{1}\mathrel{=}\Varid{const}\;\mathrm{0};\Varid{f}_{2}\mathrel{=}\Varid{fst};\Varid{f}_{3}\mathrel{=}\Varid{snd};{}\<[E]%
\\
\>[16]{}\mathbf{in}{}\<[21]%
\>[21]{}\lambda \Varid{x}\,\to\,{}\<[28]%
\>[28]{}\Conid{C3}\;(\Varid{f}_{1}\;\Varid{x},\Varid{f}_{2}\;\Varid{x},\Varid{f}_{3}\;\Varid{x}){}\<[E]%
\\
\>[3]{}\Varid{test}_{2.12}\;\Varid{n}\mathrel{=}\Varid{check}\;\Varid{n}\;(\Varid{prop}_{2.12}\;\Varid{testfun3}\mathbin{::}\Conid{Set}\;(\Int,\Int)\,\to\,\Conid{Bool}){}\<[E]%
\ColumnHook
\end{hscode}\resethooks
\noindent
The fact that property 2.12 holds and that
\DONE{This format statement should be in main.fmt, but it did not seem to work.}
\begin{hscode}\SaveRestoreHook
\column{B}{@{}>{\hspre}l<{\hspost}@{}}%
\column{3}{@{}>{\hspre}l<{\hspost}@{}}%
\column{E}{@{}>{\hspre}l<{\hspost}@{}}%
\>[3]{}\Varid{min}\;\Varid{f}\;\Varid{as}\in \Varid{map}\;\Varid{f}\;(\Varid{argparetoMin}\;\Varid{fs}\;\Varid{as})\phantom{x} \Leftrightarrow\phantom{x} \{\mskip1.5mu \Varid{min}\;\Varid{f}\;\Varid{as}\mskip1.5mu\}\subseteq\Varid{map}\;\Varid{f}\;(\Varid{argparetoMin}\;\Varid{fs}\;\Varid{as}){}\<[E]%
\ColumnHook
\end{hscode}\resethooks
\noindent
and that (because of \cref{eq1.2})
\DONE{the cref. rendered as ``Item 1.2'' which was a bit confusing. Now fixed.}
\begin{hscode}\SaveRestoreHook
\column{B}{@{}>{\hspre}l<{\hspost}@{}}%
\column{3}{@{}>{\hspre}l<{\hspost}@{}}%
\column{E}{@{}>{\hspre}l<{\hspost}@{}}%
\>[3]{}\{\mskip1.5mu \Varid{min}\;\Varid{f}\;\Varid{as}\mskip1.5mu\}\phantom{x} \,\doubleequals\,\phantom{x} \Varid{map}\;\Varid{f}\;(\Varid{argmin}\;\Varid{f}\;\Varid{as}){}\<[E]%
\ColumnHook
\end{hscode}\resethooks
\noindent
and hence
\begin{hscode}\SaveRestoreHook
\column{B}{@{}>{\hspre}l<{\hspost}@{}}%
\column{3}{@{}>{\hspre}l<{\hspost}@{}}%
\column{E}{@{}>{\hspre}l<{\hspost}@{}}%
\>[3]{}\Varid{min}\;\Varid{f}\;\Varid{as}\in \Varid{map}\;\Varid{f}\;(\Varid{argparetoMin}\;\Varid{fs}\;\Varid{as})\phantom{x} \Leftrightarrow\phantom{x} \Varid{map}\;\Varid{f}\;(\Varid{argmin}\;\Varid{f}\;\Varid{as})\subseteq\Varid{map}\;\Varid{f}\;(\Varid{argparetoMin}\;\Varid{fs}\;\Varid{as}){}\<[E]%
\ColumnHook
\end{hscode}\resethooks
\noindent
raises the question of whether 2.12 isn't perhaps a consequence of the
stronger property
\begin{hscode}\SaveRestoreHook
\column{B}{@{}>{\hspre}l<{\hspost}@{}}%
\column{3}{@{}>{\hspre}l<{\hspost}@{}}%
\column{E}{@{}>{\hspre}l<{\hspost}@{}}%
\>[3]{}\Varid{f}\in \Varid{fs}\phantom{x} \Rightarrow \phantom{x} \Varid{argmin}\;\Varid{f}\;\Varid{xs}\subseteq\Varid{argparetoMin}\;\Varid{fs}\;\Varid{xs}{}\<[E]%
\ColumnHook
\end{hscode}\resethooks
\noindent
or (remember that \ensuremath{\Varid{argparetoMin}\;(\Conid{C}_1\mathbin{\circ}\Varid{f})\;\Varid{as}\,\doubleequals\,\Varid{argmin}\;\Varid{f}\;\Varid{as}}) of

\begin{enumerate}[resume]
\item \ensuremath{\Varid{f}\in \Varid{fs}} \ \ensuremath{\Rightarrow } \ \ensuremath{\Varid{argparetoMin}\;(\Conid{C}_1\mathbin{\circ}\Varid{f})\;\Varid{as}\subseteq\Varid{argparetoMin}\;\Varid{fs}\;\Varid{as}}
\end{enumerate}

\noindent
This would be very useful in applications but unfortunately it turns out to be
false. We have
\begin{hscode}\SaveRestoreHook
\column{B}{@{}>{\hspre}l<{\hspost}@{}}%
\column{3}{@{}>{\hspre}l<{\hspost}@{}}%
\column{E}{@{}>{\hspre}l<{\hspost}@{}}%
\>[3]{}\Varid{argparetoMin}\;(\Conid{C}_1\mathbin{\circ}\Varid{const}\;\mathrm{0})\;\Varid{as}\,\doubleequals\,\Varid{as}{}\<[E]%
\ColumnHook
\end{hscode}\resethooks
but combining \ensuremath{\Varid{const}\;\mathrm{0}} with \ensuremath{\Varid{id}} we always get a singleton:
\begin{hscode}\SaveRestoreHook
\column{B}{@{}>{\hspre}l<{\hspost}@{}}%
\column{3}{@{}>{\hspre}l<{\hspost}@{}}%
\column{E}{@{}>{\hspre}l<{\hspost}@{}}%
\>[3]{}\Varid{isSingletonSet}\;(\Varid{argparetoMin}\;(\Conid{C}_2\mathbin{\circ}\lambda \Varid{x}\,\to\,(\mathrm{0},\Varid{x}))\;\Varid{as}){}\<[E]%
\ColumnHook
\end{hscode}\resethooks
thus the property does not hold in general. The first function does not
have to be constant everywhere, only on the actual inputs:

\begin{hscode}\SaveRestoreHook
\column{B}{@{}>{\hspre}l<{\hspost}@{}}%
\column{3}{@{}>{\hspre}l<{\hspost}@{}}%
\column{E}{@{}>{\hspre}l<{\hspost}@{}}%
\>[3]{}\Varid{test}_{2.13}\;\Varid{n}\mathrel{=}\mathbf{do}\;\Varid{falsify}\;\Varid{n}\;(\Varid{prop}_{2.13}\;\Varid{testfun3}\mathbin{::}\Conid{Set}\;(\Int,\Int)\,\to\,\Conid{Bool}){}\<[E]%
\ColumnHook
\end{hscode}\resethooks
\subsection{Applications}
\label{subsection:applications}

In many important applications, the domain of the objective function is
a subset \ensuremath{\Conid{Y}} of \ensuremath{\Real^m}, typically a hyper-rectangle or, in other
words, the Cartesian product of \ensuremath{\Varid{m}} intervals and \ensuremath{\Conid{B}\mathrel{=}\Real}. Thus \ensuremath{\Varid{fs}\ \mathop{:}\ \Real^m\,\to\,\Real^n}.

For example, in \citet{Pusztai2023BayesOptMMI}, \ensuremath{\Varid{m}\mathrel{=}\mathrm{2}} or \ensuremath{\Varid{m}\mathrel{=}\mathrm{3}} and
points in \ensuremath{\Conid{Y}} represents the concentrations of neon and deuterium (and,
for \ensuremath{\Varid{m}\mathrel{=}\mathrm{3}} their distribution in space) which have to be
injected\footnote{Virtually, through numerical simulation of disruption
events obtained, e.g., with DREAM \citep{HOPPE2021108098}.} into the
ITER\footnote{ITER = International Thermonuclear Experimental Reactor,
\href{https://www.iter.org}{www.iter.org}.} tokamak fusion reactor to
avoid dangerous \emph{runaway current}, a phenomenon where electrons
accelerate to relativistic speeds and escape the magnetic confinement of
the reactor.

In climate policy, \ensuremath{\Conid{Y}} may represent a combination of greenhouse gas
(GHG) emission reduction and geo-engineering measures as in
\citep{esd-10-453-2019, Moreno-Cruz+Keith2012, CARLINO202016593} or
perhaps a combination of short-term (a few decades) and long-term (a few
centuries) emission reduction measures, see \citet{10.1093/oxfclm/kgae004}.
In these applications, values in \ensuremath{\Conid{Y}} are often called \emph{controls},
the image of \ensuremath{\Conid{Y}} under \ensuremath{\Varid{fs}}, a subset of \ensuremath{\Real^n}, is called the
\emph{operational} space and one seeks Pareto optimal controls that
obtain outcomes in a \emph{safe} subset of the operational space.

For example, controls that obtain manageable values of the runaway
current and of the thermal quench \citep{Pusztai2023BayesOptMMI} or, in
climate policy, anthropogenic forcings that keep the earth system within
safe planetary boundaries \citep{ro06010m, doi:10.1126/science.1259855,
  doi:10.1126/sciadv.adh2458} or ensure that economic damages from
climate change are still manageable \citep{TOL2024113922,
  KotzLevermannWenz2024, NealNewellPitman2025, BearparkHoganHsiang2025, Schoetz2025, Nordhaus12261}.

When the control space is infinite (as in these examples),
the Pareto optimal controls and
the ideal Pareto front can only be approximated.

\paragraph*{Sampling.}

A straightforward way to approximate the ideal Pareto front and Pareto
optimal controls for multi-objective optimization problems in which
the control space \ensuremath{\Conid{Y}} is an infinite subset of \ensuremath{\Real^m} is to apply
\ensuremath{\Varid{paretoMin}} and \ensuremath{\Varid{argparetoMin}} to a finite sample of \ensuremath{\Conid{Y}}.
For each sample set, we can compute its Pareto front, and with
reasonable assumptions about the objective function, these fronts
converge to the ideal front as the sample size grows.

\begin{figure}[htbp]
  \includegraphics[width=0.49\textwidth]{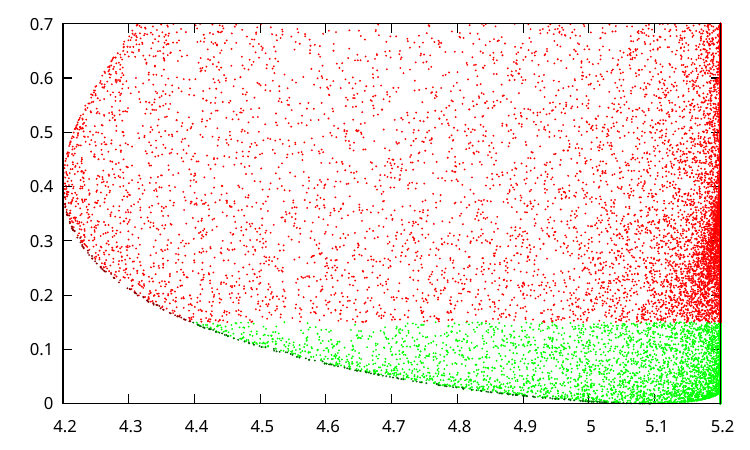}
  \includegraphics[width=0.49\textwidth]{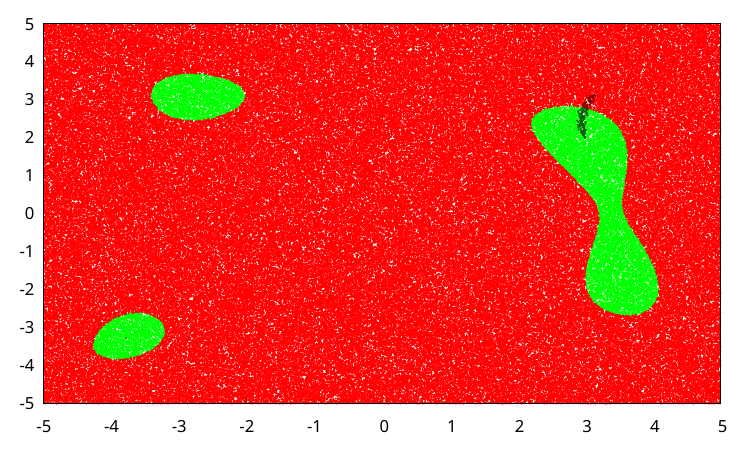}
  \caption{Benchmark MOO problem, filtering a 250k sample: operational
    space and Pareto front (left, \ensuremath{\Varid{f}_{1}} and \ensuremath{\Varid{f}_{2}} values on the \ensuremath{\Varid{x}} and
    \ensuremath{\Varid{y}} axes, respectively) and control space (right)
    partitioning: safe (green), unsafe (red), safe and unsafe Pareto
    optimal (black) points.}
  \label{fig:MOObenchmark250k}
\end{figure}

The sample can be obtained with a regular grid or by selecting a
finite number of points in \ensuremath{\Conid{Y}} randomly.
QuickCheck can also be applied to generate random samples of Cartesian
products of intervals and sampling can be combined with memoised
evaluation \citep{jansson+2025} of the objective functions, for
example, to increase the number of sampling points in regions where
the objective functions vary strongly. We do not discuss sampling
methods here, but the approach can be summarized as
\begin{hscode}\SaveRestoreHook
\column{B}{@{}>{\hspre}l<{\hspost}@{}}%
\column{3}{@{}>{\hspre}l<{\hspost}@{}}%
\column{8}{@{}>{\hspre}c<{\hspost}@{}}%
\column{8E}{@{}l@{}}%
\column{11}{@{}>{\hspre}l<{\hspost}@{}}%
\column{E}{@{}>{\hspre}l<{\hspost}@{}}%
\>[3]{}\Varid{pys}{}\<[8]%
\>[8]{}\mathrel{=}{}\<[8E]%
\>[11]{}\Varid{argparetoMin}\;\Varid{fs}\;(\Varid{randomGrid2D}\;\Varid{n}\;\Varid{ry}\;\Varid{g}){}\<[E]%
\ColumnHook
\end{hscode}\resethooks
\noindent
where \ensuremath{\Varid{n}} denotes the size of the sample, \ensuremath{\Varid{ry}} a rectangle in \ensuremath{\Real^2}
and \ensuremath{\Varid{g}} is a random number generator, see \citep{MOOcodeRepo2025a}.
\DONE{Nicola: add public repository}
\Cref{fig:MOObenchmark250k} shows the results obtained by sampling \ensuremath{\Conid{Y}\mathrel{=}[\mskip1.5mu \mathbin{-}\mathrm{5},\mathrm{5}\mskip1.5mu]\;\times\;[\mskip1.5mu \mathbin{-}\mathrm{5},\mathrm{5}\mskip1.5mu]\subseteq \Real^2} with 250k points for a standard MOO
benchmark problem with two objective functions:
\DONE{The colors in
  \cref{fig:MOObenchmark250k} are a impossible to see in a grey-scale
  print. Not sure what to do about that, though.}
\DONE{Nicola: we
  can print pictures in color, most journals welcome color
  pictures. Or we can try to make the pictures in Grey, DarkGrey and
  Black, see color palette at
  https://hackage.haskell.org/package/easyplot-1.0/docs/Graphics-EasyPlot.html}
\begin{equation*}
  \begin{split}
  f_1 (x, y) & = 5.2 - \cos(x) * \cos(y) * \exp(-((x - \pi)^2 + (y - \pi)^2)) \\
  f_2 (x, y) & = ((x^2 + y - 11)^2 + (x + y^2 -7)^2) / 100
  \end{split}
\end{equation*}
%
\DONE{Make sure the figures correspond to the seed 137
  computations. Nicola: the computations detailed the README do in fact
  correspond to the seed 137. The problem is that the black dort in the
  middle of the bottom left bubble is only visible in the corresponding
  PDF file, see updated README.}
\begin{figure}[htbp]
  \includegraphics[width=0.49\textwidth]{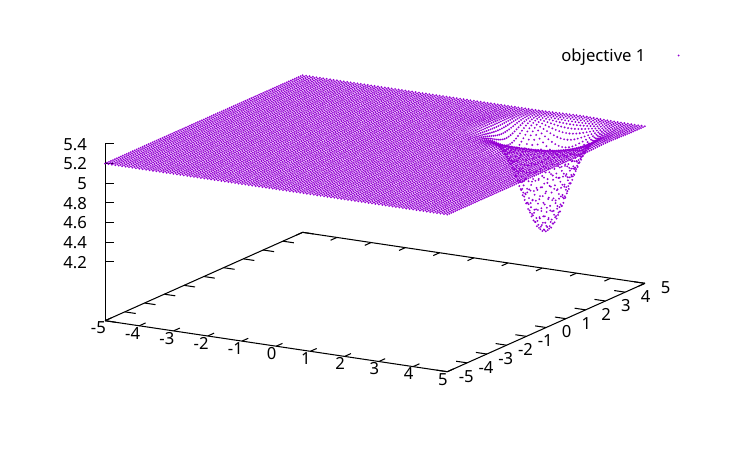}
  \includegraphics[width=0.49\textwidth]{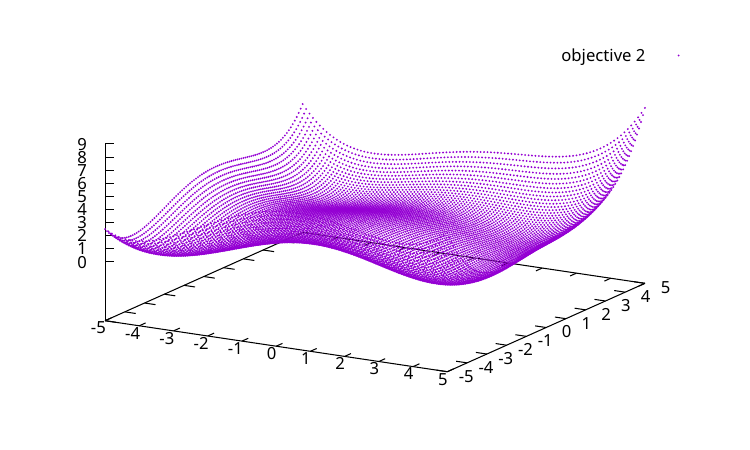}
  \caption{Benchmark MOO problem: objective functions.}
  \label{fig:BenchmarkMOOfunctions}
\end{figure}%
%
The ``ground truth'' is that \ensuremath{\Varid{f}_{1}} has a single global minimum (of 4.2)
at \ensuremath{(\pi,\pi)} while \ensuremath{\Varid{f}_{2}} has one local mininum in each quadrant (all
with value 0) at exactly \ensuremath{(\mathrm{3},\mathrm{2})} and roughly (-3.779, -3.283),
(-2.805, 3.131), and (3.584, -1.848).
Function \ensuremath{\Varid{f}_{1}} is close to its minimum 4.2 only in the first quadrant, near \ensuremath{(\pi,\pi)}.

The computed minima of \ensuremath{\Varid{f}_{1}} and \ensuremath{\Varid{f}_{2}} on the sample are at about these
values (all results rounded to three decimal places):
\begin{hscode}\SaveRestoreHook
\column{B}{@{}>{\hspre}l<{\hspost}@{}}%
\column{3}{@{}>{\hspre}l<{\hspost}@{}}%
\column{21}{@{}>{\hspre}l<{\hspost}@{}}%
\column{31}{@{}>{\hspre}l<{\hspost}@{}}%
\column{42}{@{}>{\hspre}l<{\hspost}@{}}%
\column{E}{@{}>{\hspre}l<{\hspost}@{}}%
\>[3]{}\Varid{mygrid}\mathrel{=}\Varid{randomGrid2D}\;\mathrm{250000}\;\Varid{ry}\;(\Varid{mkStdGen}\;\mathrm{137}){}\<[E]%
\\
\>[3]{}\Varid{f1min}\mathrel{=}\Varid{min}\;\Varid{f}_{1}\;{}\<[21]%
\>[21]{}\Varid{mygrid}\mathrel{=}{}\<[31]%
\>[31]{}\mathrm{4.200}{}\<[42]%
\>[42]{}\mbox{\onelinecomment  close to exact 4.2}{}\<[E]%
\\
\>[3]{}\Varid{f2min}\mathrel{=}\Varid{min}\;\Varid{f}_{2}\;{}\<[21]%
\>[21]{}\Varid{mygrid}\mathrel{=}{}\<[31]%
\>[31]{}\mathrm{6.326e{-}6}{}\<[42]%
\>[42]{}\mbox{\onelinecomment  close to exact 0}{}\<[E]%
\ColumnHook
\end{hscode}\resethooks
The corresponding points in the control space are at
\begin{hscode}\SaveRestoreHook
\column{B}{@{}>{\hspre}l<{\hspost}@{}}%
\column{3}{@{}>{\hspre}l<{\hspost}@{}}%
\column{16}{@{}>{\hspre}l<{\hspost}@{}}%
\column{38}{@{}>{\hspre}l<{\hspost}@{}}%
\column{47}{@{}>{\hspre}l<{\hspost}@{}}%
\column{E}{@{}>{\hspre}l<{\hspost}@{}}%
\>[3]{}\Varid{argmin}\;\Varid{f}_{1}\;{}\<[16]%
\>[16]{}\Varid{mygrid}\mathrel{=}[\mskip1.5mu \Conid{C}_2\;(\mathrm{3.132},{}\<[38]%
\>[38]{}\mathrm{3.130})\mskip1.5mu]{}\<[47]%
\>[47]{}\mbox{\onelinecomment  close to exact (\ensuremath{\pi}, \ensuremath{\pi})}{}\<[E]%
\\
\>[3]{}\Varid{argmin}\;\Varid{f}_{2}\;{}\<[16]%
\>[16]{}\Varid{mygrid}\mathrel{=}[\mskip1.5mu \Conid{C}_2\;(\mathrm{2.998},{}\<[38]%
\>[38]{}\mathrm{1.995})\mskip1.5mu]{}\<[47]%
\>[47]{}\mbox{\onelinecomment  close to exact (3, 2)}{}\<[E]%
\ColumnHook
\end{hscode}\resethooks
\noindent
At the first point of the Pareto front (\cref{fig:MOObenchmark250k} left) the
value of \ensuremath{\Varid{f}_{2}} is about \ensuremath{\mathrm{0.4}}
\begin{hscode}\SaveRestoreHook
\column{B}{@{}>{\hspre}l<{\hspost}@{}}%
\column{3}{@{}>{\hspre}l<{\hspost}@{}}%
\column{E}{@{}>{\hspre}l<{\hspost}@{}}%
\>[3]{}\Varid{f}_{2}\;(\Conid{C}_2\;(\mathrm{3.132},\mathrm{3.130}))\mathrel{=}\mathrm{0.389}{}\<[E]%
\ColumnHook
\end{hscode}\resethooks
\DONE{It looks like the the minimum of \ensuremath{\Varid{f}_{2}} is 0 and then \ensuremath{\Varid{f}_{1}} is
  between 5.0 and 5.1. Is this true? It seems to be between 5.0
  and 5.1 in the figure. By 5.2 the point cloud has "curved up" a
  bit. Nicola: Yes, f1 (C2 (-3.782,-3.283)) = 5.2 and argminOf f2
  (randomGrid2D 250000 ry g) = [(-3.782,-3.283)] as shown. The
  min. of f2 on the sample is 3.401e-6 and it is attained at the point
  in the middle of the bottom left green bubble. The point cloud curves
  up between 5.1 and 5.2 but there is a (Pareto) point at (5.2,
  3.401e-6).}%
\noindent
The value of \ensuremath{\Varid{f}_{1}} at the minimum of \ensuremath{\Varid{f}_{2}} is around \ensuremath{\mathrm{5.1}} as one
would expect from \cref{fig:BenchmarkMOOfunctions}
\begin{hscode}\SaveRestoreHook
\column{B}{@{}>{\hspre}l<{\hspost}@{}}%
\column{3}{@{}>{\hspre}l<{\hspost}@{}}%
\column{E}{@{}>{\hspre}l<{\hspost}@{}}%
\>[3]{}\Varid{f}_{1}\;(\Conid{C}_2\;(\mathrm{2.998},\mathrm{1.995}))\mathrel{=}\mathrm{5.093}{}\<[E]%
\ColumnHook
\end{hscode}\resethooks
Note that \ensuremath{\Varid{f}_{2}} attains its minimum on the sample at a point near \ensuremath{(\pi,\pi)} in the first quadrant of \cref{fig:MOObenchmark250k} but with
other random seeds other quadrants show up.
%
%
%
%
In this example the objective space is equipped with a safety
predicate (\ensuremath{(x_1,x_2)} is safe if \ensuremath{x_2\mathrel{\leq}\mathrm{0.15}}, see
\cref{fig:MOObenchmark250k} left) and the objective and the control
space are partitioned into safe (green) and unsafe (red)
points. Pareto optimal points are shown in black.
\DONE{Tim: I can see no dark green or dark red points. Perhaps make a
  ``zoomed in'' cutout of the main region of the (inputs to) the
  Pareto front. And change color scheme!}
\DONE{Nicola: Mention the ``finite area'' aspect of the set of control
  points mapping to the Pareto front and its relation to
  uncertainty. (Perhaps two different reasons: numerical /
  experimental inaccuracy, and perhaps also the model uncertainty.)
  Nicola: not sure it fits here and how to do this... if you have an
  idea, just go on!}

\DONE{Start of a formulation: "Even when the ideal Pareto front is a
  1D curve, the set of control points mapping to it under evaluation
  forms a region of finite area, due to numerical noise and model
  uncertainty." But it is not quite accurate.}

Even when the ideal Pareto front in objective space is a
one-dimensional curve (as in this two-objective case), the
corresponding set of control points mapping to it under evaluation
typically occupies a region of non-zero area in control space (see
\cref{fig:MOObenchmark250k} right).
In the present example, this region has a characteristic
crescent-like shape, but its precise geometry depends on the objective
functions; what is robust is its two-dimensional nature.
This apparent increase in dimensionality is not an artefact of the
visualization but reflects a structural property of multi-objective
optimization under sampling.
Near the ideal Pareto front, the objective functions often vary only
weakly along certain directions in control space: small changes in the
control may lead to negligible or compensating changes in the
objective values.
As a consequence, many nearby control points are mutually
non-dominating and are all mapped to points close to the same segment
of the Pareto front.
When the front is approximated from a finite sample, this local
flatness manifests itself as a two-dimensional region of Pareto-optimal
(or near-Pareto-optimal) controls rather than a one-dimensional curve.
From the perspective of decision making under uncertainty, this
finite-area structure is relevant: it indicates robustness, in the
sense that a range of controls yield essentially equivalent trade-offs
between objectives.

\paragraph*{Evolutionary methods, robust control.}

In many applications the objective functions are evaluated by running
computationally expensive numerical simulations.
For example, a single run of DREAM (Disruption and Runaway Electron
Analysis Model by \citet{HOPPE2021108098}) at the lowest level of
fidelity takes about 6 minutes on an off-the-shelf 4-core CPU.
Similarly, running even a reduced climate simulator like SURFER
\citep{gmd-15-8059-2022} (not to mention an earth-system model of
intermediate complexity or a general circulation model
\citep{Easterbrook_2023}) 250k times takes significant time.
\DONE{Too complex sentence: split.}
While these computations are fully parallelizable, brute-force sampling
is in practice often too expensive.

One can take advantage of the computational structure of \ensuremath{\Varid{paretoOptOf}}
(remember that \ensuremath{\Varid{paretoOpt}\mathrel{=}\Varid{fold}\;\Varid{bump}\;[\mskip1.5mu \mskip1.5mu]}) and of the fact that, in most
practical applications, one is interested in controls that are
\emph{robust}: small deviations from a selected safe, optimal control
should not yield very unsafe or very sub-optimal outcomes\footnote{There
are at least two reasons why control shall be robust: 1) in most
practical cases, decisions can be imposed only up to a certain accuracy
and 2) the results obtained through physical experiments and
computer-based simulations are typically uncertain. We discuss how to
tackle this problem in \cref{section:monadic}.}.
This means that one can neglect outliers and focus on regions of the
control space in which Pareto optimal controls are reasonably
\emph{dense}.

The observation suggests an alternative approach to sampling: start with
a small, possibly empty set of Pareto optimal controls and incrementally
grow this set. The method can be summarized in two steps:

\begin{itemize}
\item[1] Start with a small seed \ensuremath{\Varid{pys},\Varid{npys}} of controls with
\begin{hscode}\SaveRestoreHook
\column{B}{@{}>{\hspre}l<{\hspost}@{}}%
\column{3}{@{}>{\hspre}l<{\hspost}@{}}%
\column{E}{@{}>{\hspre}l<{\hspost}@{}}%
\>[3]{}\Varid{pys}\,\doubleequals\,\Varid{argparetoMin}\;\Varid{fs}\;(\Varid{npys}\mathbin{+\!\!+}\Varid{pys})\mbox{\onelinecomment  (*)}{}\<[E]%
\ColumnHook
\end{hscode}\resethooks
\item[2] Repeat a number of times:
  \begin{itemize}
  \item[2.1] Pick up a new control \ensuremath{\Varid{y}} near \ensuremath{\Varid{pys}} or explore \ensuremath{\Conid{Y}},
  \item[2.2] Add \ensuremath{\Varid{y}} to \ensuremath{\Varid{pys}} or \ensuremath{\Varid{npys}} as to preserve (*)
  \end{itemize}
\end{itemize}

\noindent
In the above program sketch \ensuremath{\Varid{pys}} and \ensuremath{\Varid{npys}} represent disjoint Pareto
optimal and non Pareto optimal controls, respectively.
It is easy to construct a seed that fulfills (*) by applying combinators
that we have already implemented. Remember that \ensuremath{\Varid{argparetoMin}\;\Varid{fs}\;\Varid{ys}}
is a subset of~\ensuremath{\Varid{ys}}.
Therefore, for any \ensuremath{\Varid{ys}}, \ensuremath{\Varid{pys}\mathrel{=}\Varid{argparetoMin}\;\Varid{fs}\;\Varid{ys}} implies \ensuremath{\Varid{pys}\,\doubleequals\,\Varid{argparetoMin}\;\Varid{fs}\;(\Varid{npys}\mathbin{+\!\!+}\Varid{pys})} with \ensuremath{\Varid{npys}\mathrel{=}\Varid{ys}\mathbin{\char92 \char92 }\Varid{pys}}. One can start
with a very coarse (or just empty) sample of controls~\ensuremath{\Varid{ys}}.

In Step 2.1 one needs a rule to pick up a new control near \ensuremath{\Varid{pys}} or to
explore \ensuremath{\Conid{Y}}, perhaps for sub-optimal but safe controls. This requires
striking a compromise between \emph{exploring} less known regions of the
control space and \emph{exploiting} the knowledge which is already
available about Pareto optimal controls.

\begin{figure}[h]
  \includegraphics[width=0.49\textwidth]{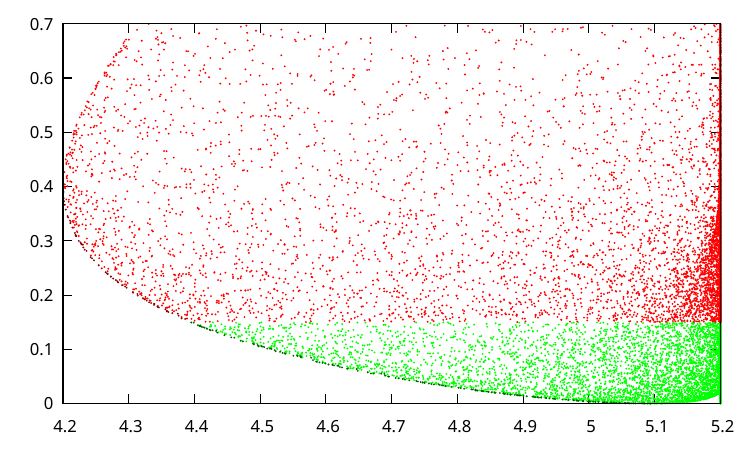}
  \includegraphics[width=0.49\textwidth]{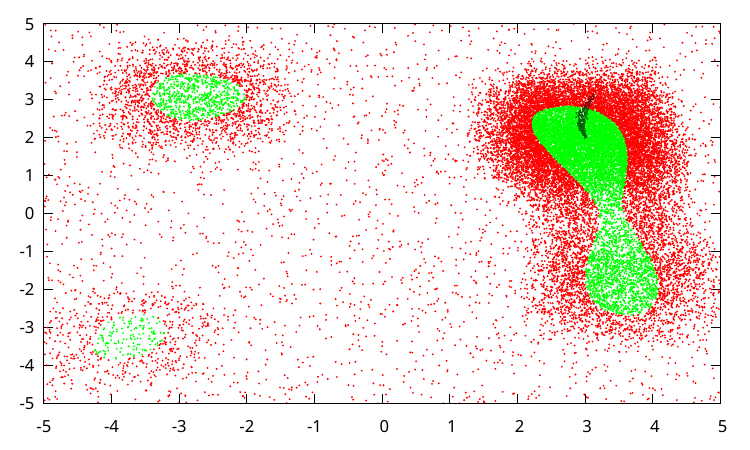}
  \caption{Benchmark MOO problem, evolving a 2500 points sample after
    22500 iterations (1/10 function evaluations as compared to
    \cref{fig:MOObenchmark250k}): operational space and Pareto front
    (left, \ensuremath{\Varid{f}_{1}} and \ensuremath{\Varid{f}_{2}} values on the \ensuremath{\Varid{x}} and \ensuremath{\Varid{y}} axes, respectively)
    and control space (right) partitioning: safe (green), unsafe (red),
    safe and unsafe Pareto optimal (black) points.}
  \label{BenchmarkMOOevolution}
\end{figure}

Step 2.2 can be implemented in terms of a modified \ensuremath{\Varid{bump}} function. This
is where \emph{evolutionary} algorithms take advantage of the
computational structure of \ensuremath{\Varid{paretoOptOf}}.
While a full treatment of evolutionary algorithms is beyond the scope
of the paper, they typically reduce the number of function evaluations
that are needed to achieve a comparable description of the Pareto
front by a factor of about 10.
\Cref{BenchmarkMOOevolution} illustrates this situation for 22500
iterations starting from a 50 $\times$ 50 sample and thus in total
25000 function evaluations, compare with \cref{fig:MOObenchmark250k}
with 10 times more function evaluations.
Results obtained starting from an empty sample are very similar to
those of \cref{BenchmarkMOOevolution}.

\section{Functorial uncertainty}
\label{section:monadic}

As discussed in the introduction, most decision problems in engineering,
economics and climate science are affected by epistemic uncertainty,
aleatoric uncertainty or a combination of both.
These can be modelled with functions that return values in \ensuremath{\Conid{U}\;\Conid{B}} (rather
than \ensuremath{\Conid{B}}) where \ensuremath{\Conid{U}} is an uncertainty functor.
We discuss uncertainty functors in
\cref{subsection:ufexamples,subsection:ufgeneral}, typical examples are
\ensuremath{\Conid{Id}}, \ensuremath{\Conid{Maybe}}, \ensuremath{\Conid{NDim}} instances (lists of values of a fixed
length), \ensuremath{\Conid{Set}} and the finite probability distributions discussed in
\citep{10.1017/S0956796805005721, ionescu2009, BREDE_BOTTA_2021}.

We specify and solve the problem of \emph{minimizing} an arbitrary
function \ensuremath{\Varid{f}\ \mathop{:}\ \Conid{A}\,\to\,\Conid{U}\;\Conid{B}} according to a measure function \ensuremath{\mu\ \mathop{:}\ \Conid{U}\;\Conid{B}\,\to\,\Conid{C}} over a finite set of values of type \ensuremath{\Conid{A}} with \ensuremath{\Conid{B}} and \ensuremath{\Conid{C}} in \ensuremath{\Conid{Ord}}.
In many applications in climate science and engineering, uncertainty is
treated implicitly: \ensuremath{\Conid{B}} and \ensuremath{\Conid{C}} are assumed to be equal to \ensuremath{\Real},
\ensuremath{\mu} is assumed to be the expected value function and optimization
under uncertainty is reduced to the application of \ensuremath{\Varid{min}} and
\ensuremath{\Varid{argmin}} to the composition \ensuremath{\mu\mathbin{\circ}\Varid{f}\ \mathop{:}\ \Conid{A}\,\to\,\Conid{C}}.
This approach is straightforward but implicitly assumes that decision
makers are risk-neutral, and it does not distinguish between different
values in the inverse image of \ensuremath{\Varid{min}\;(\mu\mathbin{\circ}\Varid{f})\;\Varid{as}} through \ensuremath{\mu}.

Here we avoid making these implicit assumptions.
Instead, we generalize \ensuremath{\Varid{min}} and \ensuremath{\Varid{argmin}} to functions that take a
measure function as an additional argument and return sets of values
in \ensuremath{\Conid{U}\;\Conid{B}} and \ensuremath{\Conid{A}}, respectively:
\begin{hscode}\SaveRestoreHook
\column{B}{@{}>{\hspre}l<{\hspost}@{}}%
\column{3}{@{}>{\hspre}l<{\hspost}@{}}%
\column{14}{@{}>{\hspre}c<{\hspost}@{}}%
\column{14E}{@{}l@{}}%
\column{18}{@{}>{\hspre}l<{\hspost}@{}}%
\column{47}{@{}>{\hspre}l<{\hspost}@{}}%
\column{E}{@{}>{\hspre}l<{\hspost}@{}}%
\>[3]{}\Varid{min}_\Varid{u}{}\<[14]%
\>[14]{}\mathbin{::}{}\<[14E]%
\>[18]{}(\Conid{Ord}\;\Varid{b},\Conid{Ord}\;\Varid{c})\Rightarrow {}\<[47]%
\>[47]{}(\Varid{u}\;\Varid{b}\,\to\,\Varid{c})\,\to\,(\Varid{a}\,\to\,\Varid{u}\;\Varid{b})\,\to\,\Conid{Set}\;\Varid{a}\,\to\,\Conid{Set}\;(\Varid{u}\;\Varid{b}){}\<[E]%
\\[\blanklineskip]%
\>[3]{}\Varid{argmin}_\Varid{u}{}\<[14]%
\>[14]{}\mathbin{::}{}\<[14E]%
\>[18]{}(\Conid{Ord}\;\Varid{b},\Conid{Ord}\;\Varid{c},\Conid{Eq}\;(\Varid{u}\;\Varid{b}))\Rightarrow {}\<[47]%
\>[47]{}(\Varid{u}\;\Varid{b}\,\to\,\Varid{c})\,\to\,(\Varid{a}\,\to\,\Varid{u}\;\Varid{b})\,\to\,\Conid{Set}\;\Varid{a}\,\to\,\Conid{Set}\;\Varid{a}{}\<[E]%
\ColumnHook
\end{hscode}\resethooks
\noindent
Notice the additional constraint \ensuremath{\Conid{Eq}\;(\Varid{u}\;\Varid{b})} in \ensuremath{\Varid{argmin}_\Varid{u}}.
In the implementation of \ensuremath{\Varid{argmin}} from \cref{section:intro}, \ensuremath{\Conid{Ord}\;\Conid{B}}
was sufficient for values in \ensuremath{\Conid{B}} to be comparable for equality.
In order to implement \ensuremath{\Varid{argmin}_\Varid{u}}, we also need to compare values in
\ensuremath{\Conid{U}\;\Conid{B}} for equality.
With \ensuremath{\Conid{U}\;\Conid{B}} in \ensuremath{\Conid{Eq}}, it is straightforward to implement a function that
computes the inverse image of \ensuremath{\Varid{min}\;(\mu\mathbin{\circ}\Varid{f})\;\Varid{as}} through \ensuremath{\mu}
generically in terms of \ensuremath{\Varid{argmin}}:
\begin{hscode}\SaveRestoreHook
\column{B}{@{}>{\hspre}l<{\hspost}@{}}%
\column{3}{@{}>{\hspre}l<{\hspost}@{}}%
\column{E}{@{}>{\hspre}l<{\hspost}@{}}%
\>[3]{}\Varid{min}_\Varid{u}\;\mu\;\Varid{f}\mathrel{=}\Varid{argmin}\;\mu\mathbin{\circ}\Varid{map}\;\Varid{f}{}\<[E]%
\ColumnHook
\end{hscode}\resethooks
With \ensuremath{\Varid{min}_\Varid{u}} in place, \ensuremath{\Varid{argmin}_\Varid{u}} can be defined through the same
(inefficient but obviously correct) pattern used for \ensuremath{\Varid{argmin}} in the
introduction:
\begin{hscode}\SaveRestoreHook
\column{B}{@{}>{\hspre}l<{\hspost}@{}}%
\column{3}{@{}>{\hspre}l<{\hspost}@{}}%
\column{24}{@{}>{\hspre}c<{\hspost}@{}}%
\column{24E}{@{}l@{}}%
\column{27}{@{}>{\hspre}l<{\hspost}@{}}%
\column{E}{@{}>{\hspre}l<{\hspost}@{}}%
\>[3]{}\Varid{argmin}_\Varid{u}\;\mu\;\Varid{f}\;\Varid{as}{}\<[24]%
\>[24]{}\mathrel{=}{}\<[24E]%
\>[27]{}\mathbf{let}\;\Varid{ms}\mathrel{=}\Varid{min}_\Varid{u}\;\mu\;\Varid{f}\;\Varid{as}\;\mathbf{in}\;\Varid{filter}\;(\lambda \Varid{a}\,\to\,\Varid{f}\;\Varid{a}\in \Varid{ms})\;\Varid{as}{}\<[E]%
\ColumnHook
\end{hscode}\resethooks
\paragraph*{Structural Comparison of Uncertain Values.}
Before we can specify the properties of these functions, we have to
clarify what it means to compare uncertain values.
In understanding orders under value uncertainty in \cref{section:moo},
we introduced the partial relation \ensuremath{(\prec)} to compare values in \ensuremath{B^n}.
Here, we seek a relation \ensuremath{(\mathrel{\prec_{u}})} that allows us to compare values in \ensuremath{\Conid{U}\;\Conid{B}}
(where \ensuremath{\Conid{U}} can be \ensuremath{\Conid{Id}}, \ensuremath{\Conid{Maybe}}, \ensuremath{\Conid{NDim}} instances, \ensuremath{\Conid{Set}}, etc.)
consistently with \ensuremath{\Varid{ub}\in \Varid{min}_\Varid{u}\;\mu\;\Varid{f}\;\Varid{as}} being a solution of a
minimization problem.

A crucial observation is that \ensuremath{\Conid{Id}}, \ensuremath{\Conid{Maybe}}, \ensuremath{\Conid{NDim}} instances, etc.
support quantifiers and a membership relation.
Remember that the values in \ensuremath{\Varid{map}\;\Varid{f}\;\Varid{as}} encode (non-empty) sets of
``possible outcomes''%
\footnote{For example, of a numerical experiment obtained by applying
  a deterministic model to \ensuremath{\Varid{as}} for different values of some model
  parameter (epistemic uncertainty) or by integrating a genuinely
  non-deterministic model (aleatoric uncertainty)}, perhaps associated
with specific weights or probabilities.
This suggests defining \ensuremath{(\mathrel{\prec_{u}})} (``universally strictly better'') to ensure
that \emph{every} outcome on the left is smaller than \emph{every} outcome
on the right:
\DONE{Check the logic related to the unrelated cases.}
\begin{hscode}\SaveRestoreHook
\column{B}{@{}>{\hspre}l<{\hspost}@{}}%
\column{3}{@{}>{\hspre}l<{\hspost}@{}}%
\column{E}{@{}>{\hspre}l<{\hspost}@{}}%
\>[3]{}ub_1\mathrel{\prec_{u}}ub_2\phantom{x} \mathrel{=}\phantom{x} \forall\ \Varid{x}\in ub_1.\ \ \forall\ \Varid{y}\in ub_2.\ \ \Varid{x}\mathbin{<}\Varid{y}{}\<[E]%
\ColumnHook
\end{hscode}\resethooks
\noindent
This relation is anti-reflexive, transitive and, just like the dominance
relation from \cref{subsection:pareto}, partial.
It is notably coarse: a single point of overlap makes any two sets incomparable.
Notice that, with \ensuremath{\Conid{B}} in \ensuremath{\Conid{Ord}}, \ensuremath{(\mathrel{\prec_{u}})} can be defined
straightforwardly for \ensuremath{\Conid{U}\mathrel{=}\Conid{Set}}.
We discuss \ensuremath{(\mathrel{\prec_{u}})} for other examples of uncertainty functors in
\cref{subsection:ufexamples} and, more in general, in
\cref{subsection:ufgeneral}.

\DONE{Check formulation (Nicola 26-02-06)}
From this strict partial order, we derive the negation \ensuremath{ub_1\mathrel{\precsim_{u}}ub_2},
representing that \ensuremath{ub_1} is ``not strictly dominated by'' \ensuremath{ub_2}:
\begin{hscode}\SaveRestoreHook
\column{B}{@{}>{\hspre}l<{\hspost}@{}}%
\column{3}{@{}>{\hspre}l<{\hspost}@{}}%
\column{16}{@{}>{\hspre}c<{\hspost}@{}}%
\column{16E}{@{}l@{}}%
\column{19}{@{}>{\hspre}l<{\hspost}@{}}%
\column{E}{@{}>{\hspre}l<{\hspost}@{}}%
\>[3]{}ub_1\mathrel{\precsim_{u}}ub_2{}\<[16]%
\>[16]{}\mathrel{=}{}\<[16E]%
\>[19]{}\neg \;(ub_2\mathrel{\prec_{u}}ub_1){}\<[E]%
\ColumnHook
\end{hscode}\resethooks
\paragraph*{Properties of functorial minimization.}
We can now state the soundness requirements for our minimization functions.
To rule out trivial solutions, we first require that non-empty inputs
yield non-empty outputs that are subsets of the image:

\begin{enumerate}
  \setcounter{enumi}{-1}
  \setcounter{speceq}{-1}
\item \ensuremath{\Varid{as}\not\doubleequals\{\mskip1.5mu \mskip1.5mu\}} \ \ensuremath{\Rightarrow } \ \ensuremath{\{\mskip1.5mu \mskip1.5mu\}\not\doubleequals\Varid{min}_\Varid{u}\;\mu\;\Varid{f}\;\Varid{as}\subseteq\Varid{map}\;\Varid{f}\;\Varid{as}}                      \speclabel{eq3.0}
\end{enumerate}

\noindent
Next, we adapt the properties from the pure case (\cref{section:intro}) to
this context.
There we had:
\begin{enumerate}
\item[1.1] \ensuremath{\forall \Varid{a}\in \Varid{as}}. \ \ensuremath{\Varid{min}\;\Varid{f}\;\Varid{as}\mathrel{\leq}\Varid{f}\;\Varid{a}}
\item[1.2] \ensuremath{\forall \Varid{a}\in \Varid{as}}. \ \ensuremath{\Varid{a}\in \Varid{argmin}\;\Varid{f}\;\Varid{as}} \ \ensuremath{\Rightarrow } \ \ensuremath{\Varid{f}\;\Varid{a}\mathrel{=}\Varid{min}\;\Varid{f}\;\Varid{as}}
\item[1.3] \ensuremath{\forall \Varid{a}\in \Varid{as}}. \ \ensuremath{\Varid{f}\;\Varid{a}\mathrel{=}\Varid{min}\;\Varid{f}\;\Varid{as}}  \ \ensuremath{\Rightarrow } \ \ensuremath{\Varid{a}\in \Varid{argmin}\;\Varid{f}\;\Varid{as}}
\end{enumerate}
\noindent
For consistency, we would like \ensuremath{\Varid{min}_\Varid{u}}, \ensuremath{\Varid{argmin}_\Varid{u}} to fulfil these
conditions when \ensuremath{\Conid{U}} is the identity functor and \ensuremath{\mu\mathrel{=}\Varid{id}}.
Our definitions satisfy this ``modulo typing'': \ensuremath{\Varid{min}_\Varid{u}\;\Varid{id}\;\Varid{f}\;\Varid{as}\doubleequals\{\mskip1.5mu \Varid{min}\;\Varid{f}\;\Varid{as}\mskip1.5mu\}}.

In the presence of uncertainty,
%
a sound minimization strategy should never return a result \ensuremath{\Varid{ub}\in \Varid{min}_\Varid{u}\;\mu\;\Varid{f}\;\Varid{as}} that is \emph{strictly dominated} by some other available
option in \ensuremath{\Varid{map}\;\Varid{f}\;\Varid{as}}.
%
This leads to the specification:

\begin{enumerate}
\item \ensuremath{\forall \Varid{a}\in \Varid{as}}. \ \ensuremath{\forall \Varid{ub}\in \Varid{min}_\Varid{u}\;\mu\;\Varid{f}\;\Varid{as}}. \ \ensuremath{\Varid{ub}\mathrel{\precsim_{u}}\Varid{f}\;\Varid{a}}                         \speclabel{eq3.1}
\item \ensuremath{\forall \Varid{a}\in \Varid{as}}. \ \ensuremath{\Varid{a}\in \Varid{argmin}_\Varid{u}\;\mu\;\Varid{f}\;\Varid{as}} \ \ensuremath{\Rightarrow } \ \ensuremath{\Varid{f}\;\Varid{a}\in \Varid{min}_\Varid{u}\;\mu\;\Varid{f}\;\Varid{as}}       \speclabel{eq3.2}
\item \ensuremath{\forall \Varid{a}\in \Varid{as}}. \ \ensuremath{\Varid{f}\;\Varid{a}\in \Varid{min}_\Varid{u}\;\mu\;\Varid{f}\;\Varid{as}}  \ \ensuremath{\Rightarrow } \ \ensuremath{\Varid{a}\in \Varid{argmin}_\Varid{u}\;\mu\;\Varid{f}\;\Varid{as}}      \speclabel{eq3.3}
\end{enumerate}

\noindent
%
Notice that \cref{eq3.1} replaces the pure condition \ensuremath{\Varid{min}\;\Varid{f}\;\Varid{as}\mathrel{\leq}\Varid{f}\;\Varid{a}}
with the new condition \ensuremath{\Varid{ub}\mathrel{\precsim_{u}}\Varid{f}\;\Varid{a}} using the uncertainty order.
The specification does not mention the measure function \ensuremath{\mu}, but the implementation uses it.
This raises the central question of this section.

\paragraph*{What are suitable measure functions?}
We consider a measure function \ensuremath{\mu} to be suitable if it guarantees
that the operationally defined \ensuremath{\Varid{min}_\Varid{u}\;\mu} satisfies the structural
specification \crefrange{eq3.0}{eq3.3}.

As mentioned earlier, the expected value is popular but assumes risk neutrality.
Constant functions are unsuitable as they fail to discriminate, causing
\ensuremath{\Varid{min}_\Varid{u}} to return the entire set \ensuremath{\Varid{map}\;\Varid{f}\;\Varid{as}}
.
\DONE{Nicola: why does \ensuremath{\Varid{map}\;\Varid{f}\;\Varid{as}} technically satisfy \cref{eq3.1}? I
  do not understand the remark ``which technically satisfies
  \cref{eq3.1} due to the coarseness of \ensuremath{\mathrel{\prec_{u}}}, but is practically
  useless'', later on we do in fact show that constant measures fail to
  satisfy \cref{eq3.1}!}
We follow the approach originally proposed by \citet{ionescu2009,
 IONESCU_2016b} and formally derive the specifications \ensuremath{\mu} must fulfill
in \cref{subsection:monotonicity}.
\subsection{Monotonicity conditions}
\label{subsection:monotonicity}

With \cref{eq3.1} we have required that a sound minimization strategy
shall not return results that are structurally strictly dominated (with
respect to \ensuremath{(\mathrel{\prec_{u}})}) by values of \ensuremath{\Varid{map}\;\Varid{f}\;\Varid{as}}.
We now investigate which properties the measure function \ensuremath{\mu} must possess
to fulfil this specification.

\paragraph*{Pointwise Monotonicity (M1).}
As discussed in \citet{ionescu2009, IONESCU_2016b}, it is natural to expect
measure functions to fulfill a \emph{pointwise} monotonicity condition.
If we increase the value of some elements inside a structure, the
measure should not decrease.
We denote this pointwise comparison as \ensuremath{ub_1\mathrel{\leq\!\cdot}ub_2}, meaning \ensuremath{ub_2} has
the same shape as \ensuremath{ub_1} but with element-wise greater-or-equal values.
This leads to condition (M1):

\begin{itemize}
\item[(M1)] \ensuremath{ub_1\mathrel{\leq\!\cdot}ub_2\Rightarrow \mu\;ub_1\mathrel{\leq}\mu\;ub_2} \speclabel{eqM1}
\end{itemize}

\noindent
While (M1) is a sensible property that any measure function should
fulfil, it is \emph{insufficient} to guarantee the soundness of
\ensuremath{\Varid{min}_\Varid{u}}.
%
The reason is that (M1) is fulfilled by measure functions that do not
depend on the values of the elements of their arguments, for example
constant functions.
\DONE{Nicola: I am not sure this is the reason why (M1) is not
  sufficient to guarantee the soundness of \ensuremath{\Varid{min}_\Varid{u}}. I think the reason
  is that (M1) is fulfilled by measure functions that do not depend on
  the values of the elements of their arguments, for example constant
  functions.}

Consider the case where \ensuremath{\mu} is a constant function (e.g., \ensuremath{\Varid{const}\;\mathrm{42}}) or
simple geometric properties like \ensuremath{\Varid{length}}. These satisfy (M1).
\DONE{Nicola: why often? \ensuremath{\Varid{const}\;\mathrm{42}} and \ensuremath{\Varid{length}} do satisfy (M1).}
However, take \ensuremath{\Conid{B}\mathrel{=}\Int}, \ensuremath{\Conid{U}\mathrel{=}\Conid{List}}, \ensuremath{\Varid{as}\mathrel{=}} \ensuremath{\{\mskip1.5mu \mathrm{0},\mathrm{1}\mskip1.5mu\}}, \ensuremath{\Varid{f}\;\mathrm{0}\mathrel{=}[\mskip1.5mu \mathrm{3},\mathrm{2},\mathrm{2}\mskip1.5mu]} and \ensuremath{\Varid{f}\;\mathrm{1}\mathrel{=}[\mskip1.5mu \mathrm{1},\mathrm{1},\mathrm{0}\mskip1.5mu]}.
Structurally, \ensuremath{[\mskip1.5mu \mathrm{3},\mathrm{2},\mathrm{2}\mskip1.5mu]} is strictly dominated by \ensuremath{[\mskip1.5mu \mathrm{1},\mathrm{1},\mathrm{0}\mskip1.5mu]} (since $3>1, 2>1, 2>0$).
Yet, if \ensuremath{\mu} is constant:
1. \ensuremath{\Varid{min}\;\mu\;\{\mskip1.5mu \Varid{f}\;\mathrm{0},\Varid{f}\;\mathrm{1}\mskip1.5mu\}} finds both values minimal (equal measure);
2. \ensuremath{\Varid{min}_\Varid{u}} returns both \ensuremath{\{\mskip1.5mu [\mskip1.5mu \mathrm{3},\mathrm{2},\mathrm{2}\mskip1.5mu],[\mskip1.5mu \mathrm{1},\mathrm{1},\mathrm{0}\mskip1.5mu]\mskip1.5mu\}};
3. this violates \cref{eq3.1}, because the result contains \ensuremath{[\mskip1.5mu \mathrm{3},\mathrm{2},\mathrm{2}\mskip1.5mu]}
   which is strictly dominated by the reachable option \ensuremath{\Varid{f}\;\mathrm{1}}.

It follows that (M1) is too weak to enforce soundness.

\paragraph*{Structural Monotonicity (M2).}
To ensure consistency with the structural order \ensuremath{(\mathrel{\prec_{u}})} used in our specification,
we propose a stronger condition that links \ensuremath{\mu} directly to strict dominance:

\begin{itemize}
\item[(M2)] \ensuremath{ub_1\mathrel{\prec_{u}}ub_2\Rightarrow \mu\;ub_1\mathbin{<}\mu\;ub_2} \speclabel{eqM2}
\end{itemize}

\noindent
This condition effectively requires that if \ensuremath{ub_1} is universally strictly
better than \ensuremath{ub_2}, the measure \emph{must} reflect this by assigning a
strictly lower value.

Unlike (M1), this condition is violated by constant functions (since \ensuremath{\Varid{m}\mathbin{<}\Varid{m}} is false).
In the example above, because \ensuremath{[\mskip1.5mu \mathrm{1},\mathrm{1},\mathrm{0}\mskip1.5mu]\mathrel{\prec_{u}}[\mskip1.5mu \mathrm{3},\mathrm{2},\mathrm{2}\mskip1.5mu]}, (M2) requires
\ensuremath{\mu\;[\mskip1.5mu \mathrm{1},\mathrm{1},\mathrm{0}\mskip1.5mu]\mathbin{<}\mu\;[\mskip1.5mu \mathrm{3},\mathrm{2},\mathrm{2}\mskip1.5mu]}.
Consequently, (M2) rules out any measure that depends solely on the
structure (e.g., \ensuremath{\Varid{length}}), as such measures cannot distinguish between
strictly dominated options of the same shape.

\paragraph*{Relating M1 and M2.}
It is important to note that these two conditions are logically distinct:
\begin{itemize}
\item {(M1) does not imply (M2):} As seen above, \ensuremath{\Varid{const}\;\Varid{c}} satisfies
  (M1) but fails (M2).
\item {(M2) does not imply (M1):} There are measure that respect
  structural dominance but violate pointwise monotonicity in specific
  edge cases, we discuss an example in the last paragraph of \cref{subsection:ufexamples}.
\end{itemize}
\DONE{Nicola: I do not see that we discuss an example of a measure
  function that fulfils (M2) but not (M1) in
  \cref{subsection:ufexamples}.}
In \cref{subsection:testing,subsection:ufexamples}, we systematically test
standard measures against these conditions.
Our results support the conjecture that (M2) is the necessary and sufficient
condition for \ensuremath{\Varid{min}_\Varid{u}} to fulfill the soundness specification \cref{eq3.1}.

\DONE{Check the property (due to possible changes of order relation in
  M2). (Nicola 26-02-06)}
\DONE{Perhaps add semi formal proofs, see ``if False - endif'' below?}

\subsection{Testing the theory: monotonicity conditions and \ensuremath{\Varid{min}_\Varid{u}} specification}
\label{subsection:testing}

In this section, we test the theory outlined in the previous sections
for \ensuremath{\Conid{U}\mathrel{=}\Conid{List}}.
Specifically, we check which measure functions fulfil the monotonicity
conditions (M1) and (M2) from \cref{subsection:monotonicity} and confirm
that only those satisfying (M2) ensure that \ensuremath{\Varid{min}_\Varid{u}} and \ensuremath{\Varid{argmin}_\Varid{u}}
satisfy the structural specifications \crefrange{eq3.1}{eq3.3}.

It is easy to see that, for non-empty lists, constant functions and
\ensuremath{\Varid{length}} fulfil the pointwise monotonicity condition (M1):

\testresult{\ensuremath{\Varid{testM1}_{[]}\;\mathrm{10000}\;\Varid{length}}}{+++ OK, passed 10000 tests}
\hsnewpar{-2\belowdisplayskip}
\testresult{\ensuremath{\Varid{testM1}_{[]}\;\mathrm{10000}\;(\Varid{const}\;\mathrm{3})}}{+++ OK, passed 10000 tests}%
\noindent
with
\begin{hscode}\SaveRestoreHook
\column{B}{@{}>{\hspre}l<{\hspost}@{}}%
\column{3}{@{}>{\hspre}l<{\hspost}@{}}%
\column{5}{@{}>{\hspre}l<{\hspost}@{}}%
\column{13}{@{}>{\hspre}c<{\hspost}@{}}%
\column{13E}{@{}l@{}}%
\column{16}{@{}>{\hspre}l<{\hspost}@{}}%
\column{21}{@{}>{\hspre}l<{\hspost}@{}}%
\column{26}{@{}>{\hspre}l<{\hspost}@{}}%
\column{E}{@{}>{\hspre}l<{\hspost}@{}}%
\>[3]{}\Varid{testM1}_{[]}\;\Varid{n}\;\mu\mathrel{=}\mathbf{do}\;\Varid{check}\;\Varid{n}\;(\Varid{p}\mathbin{::}\Int\,\to\,[\mskip1.5mu \Int\mskip1.5mu]\,\to\,\Conid{Bool})\;\mathbf{where}{}\<[E]%
\\
\>[3]{}\hsindent{2}{}\<[5]%
\>[5]{}\Varid{p}\;\Varid{u}\;\Varid{ub}{}\<[13]%
\>[13]{}\mathrel{=}{}\<[13E]%
\>[16]{}\mathbf{let}\;{}\<[21]%
\>[21]{}\Varid{ub}_1{}\<[26]%
\>[26]{}\mathrel{=}\Varid{u}\ \mathop{:}\ \Varid{ub}{}\<[E]%
\\
\>[21]{}\Varid{ub}_2{}\<[26]%
\>[26]{}\mathrel{=}\Varid{fmap}\;(\mathrm{1}\mathbin{+})\;\Varid{ub}_1{}\<[E]%
\\
\>[16]{}\mathbf{in}\;{}\<[21]%
\>[21]{}(\Varid{ub}_1\mathrel{\leq\!\cdot}\Varid{ub}_2)\Rightarrow (\mu\;\Varid{ub}_1\mathrel{\leq}\mu\;\Varid{ub}_2){}\<[E]%
\ColumnHook
\end{hscode}\resethooks
\noindent
However, these same measures fail the structural monotonicity condition (M2):

\testresult{\ensuremath{\Varid{testM2}_{[]}\;\mathrm{10000}\;\Varid{length}}}{*** Failed, Falsified (after 1 test)}%
\noindent
with
\begin{hscode}\SaveRestoreHook
\column{B}{@{}>{\hspre}l<{\hspost}@{}}%
\column{3}{@{}>{\hspre}l<{\hspost}@{}}%
\column{5}{@{}>{\hspre}l<{\hspost}@{}}%
\column{13}{@{}>{\hspre}c<{\hspost}@{}}%
\column{13E}{@{}l@{}}%
\column{16}{@{}>{\hspre}l<{\hspost}@{}}%
\column{21}{@{}>{\hspre}l<{\hspost}@{}}%
\column{26}{@{}>{\hspre}l<{\hspost}@{}}%
\column{E}{@{}>{\hspre}l<{\hspost}@{}}%
\>[3]{}\Varid{testM2}_{[]}\;\Varid{n}\;\mu\mathrel{=}\mathbf{do}\;\Varid{check}\;\Varid{n}\;(\Varid{p}\mathbin{::}\Int\,\to\,[\mskip1.5mu \Int\mskip1.5mu]\,\to\,\Conid{Bool})\;\mathbf{where}{}\<[E]%
\\
\>[3]{}\hsindent{2}{}\<[5]%
\>[5]{}\Varid{p}\;\Varid{u}\;\Varid{ub}{}\<[13]%
\>[13]{}\mathrel{=}{}\<[13E]%
\>[16]{}\mathbf{let}\;{}\<[21]%
\>[21]{}\Varid{ub}_1{}\<[26]%
\>[26]{}\mathrel{=}\Varid{u}\ \mathop{:}\ \Varid{ub}{}\<[E]%
\\
\>[21]{}\Varid{shift}\mathrel{=}\mathrm{1}\mathbin{+}(\Varid{maximum}\;\Varid{ub}_1\mathbin{-}\Varid{minimum}\;\Varid{ub}_1){}\<[E]%
\\
\>[21]{}\Varid{ub}_2{}\<[26]%
\>[26]{}\mathrel{=}\Varid{fmap}\;(\Varid{shift}\mathbin{+})\;\Varid{ub}_1{}\<[E]%
\\
\>[16]{}\mathbf{in}\;{}\<[21]%
\>[21]{}(\Varid{ub}_1\mathrel{\prec_{u}}\Varid{ub}_2)\Rightarrow (\mu\;\Varid{ub}_1\mathbin{<}\mu\;\Varid{ub}_2){}\<[E]%
\ColumnHook
\end{hscode}\resethooks
\noindent
This failure to satisfy (M2) directly correlates with the violation of
the soundness specification \cref{eq3.1}.
As shown below, using such measures causes \ensuremath{\Varid{min}_\Varid{u}} to return solutions
that are structurally strictly dominated by other available options:

\testresult{\ensuremath{\Varid{test}_{3.1}\;\mathrm{10000}\;\,\Varid{f}_{[]}\,\;\Varid{length}}}{*** Failed, Falsified (after 9 tests and 5 shrinks)}
\hsnewpar{-2\belowdisplayskip}
\testresult{\ensuremath{\Varid{test}_{3.1}\;\mathrm{10000}\;\,\Varid{f}_{[]}\,\;(\Varid{const}\;\mathrm{3})}}{*** Failed, Falsified (after 9 tests and 3 shrinks)}%
\noindent
with
\begin{hscode}\SaveRestoreHook
\column{B}{@{}>{\hspre}l<{\hspost}@{}}%
\column{3}{@{}>{\hspre}l<{\hspost}@{}}%
\column{5}{@{}>{\hspre}l<{\hspost}@{}}%
\column{11}{@{}>{\hspre}c<{\hspost}@{}}%
\column{11E}{@{}l@{}}%
\column{14}{@{}>{\hspre}l<{\hspost}@{}}%
\column{18}{@{}>{\hspre}l<{\hspost}@{}}%
\column{E}{@{}>{\hspre}l<{\hspost}@{}}%
\>[3]{}\Varid{test}_{3.1}\;\Varid{n}\;\Varid{f}\;\mu\mathrel{=}\mathbf{do}\;\Varid{check}\;\Varid{n}\;\Varid{p}\;\mathbf{where}{}\<[E]%
\\
\>[3]{}\hsindent{2}{}\<[5]%
\>[5]{}\Varid{p}\;\Varid{as}{}\<[11]%
\>[11]{}\mathrel{=}{}\<[11E]%
\>[14]{}\mathbf{let}\;\Varid{ubs}\mathrel{=}\Varid{min}_\Varid{u}\;\mu\;\Varid{f}\;\Varid{as}{}\<[E]%
\\
\>[14]{}\mathbf{in}\;{}\<[18]%
\>[18]{}\Varid{forAll}\;\Varid{as}\;(\lambda \Varid{a}\,\to\,\Varid{forAll}\;\Varid{ubs}\;(\lambda \Varid{ub}\,\to\,\Varid{ub}\mathrel{\precsim_{u}}\Varid{f}\;\Varid{a})){}\<[E]%
\ColumnHook
\end{hscode}\resethooks
\noindent
and the counter-example function derived in \cref{subsection:monotonicity}:
\begin{hscode}\SaveRestoreHook
\column{B}{@{}>{\hspre}l<{\hspost}@{}}%
\column{3}{@{}>{\hspre}l<{\hspost}@{}}%
\column{E}{@{}>{\hspre}l<{\hspost}@{}}%
\>[3]{}\,\Varid{f}_{[]}\,\mathbin{::}\Int\,\to\,[\mskip1.5mu \Int\mskip1.5mu]{}\<[E]%
\\
\>[3]{}\,\Varid{f}_{[]}\,\;\mathrm{0}\mathrel{=}[\mskip1.5mu \mathrm{3},\mathrm{2},\mathrm{2}\mskip1.5mu];\,\Varid{f}_{[]}\,\;\mathrm{1}\mathrel{=}[\mskip1.5mu \mathrm{1},\mathrm{1},\mathrm{0}\mskip1.5mu];\,\Varid{f}_{[]}\,\;\Varid{n}\mathrel{=}[\mskip1.5mu \mathrm{3},\mathrm{1},\mathrm{0}\mskip1.5mu]{}\<[E]%
\ColumnHook
\end{hscode}\resethooks
\noindent
By contrast, measures that satisfy (M2), such as \ensuremath{\Varid{sum}_{[]}}, \ensuremath{\Varid{average}_{[]}},
\ensuremath{\Varid{head}_{[]}}, \ensuremath{\Varid{best}_{[]}} and \ensuremath{\Varid{worst}_{[]}}, allow \ensuremath{\Varid{min}_\Varid{u}} to satisfy \cref{eq3.1}.

Finally, it is worth noting that \cref{eq3.2,eq3.3} are satisfied by any
measure function by the very definition of \ensuremath{\Varid{min}_\Varid{u}} and \ensuremath{\Varid{argmin}_\Varid{u}}.
This is confirmed by:
\begin{hscode}\SaveRestoreHook
\column{B}{@{}>{\hspre}l<{\hspost}@{}}%
\column{3}{@{}>{\hspre}l<{\hspost}@{}}%
\column{5}{@{}>{\hspre}l<{\hspost}@{}}%
\column{E}{@{}>{\hspre}l<{\hspost}@{}}%
\>[3]{}\Varid{test}_{3.2}\;\Varid{n}\;\Varid{f}\;\mu\mathrel{=}\mathbf{do}\;\Varid{check}\;\Varid{n}\;\Varid{p}\;\mathbf{where}{}\<[E]%
\\
\>[3]{}\hsindent{2}{}\<[5]%
\>[5]{}\Varid{p}\;\Varid{as}\mathrel{=}\Varid{forAll}\;\Varid{as}\;(\lambda \Varid{a}\,\to\,(\Varid{a}\in \Varid{argmin}_\Varid{u}\;\mu\;\Varid{f}\;\Varid{as})\Rightarrow (\Varid{f}\;\Varid{a}\in \Varid{min}_\Varid{u}\;\mu\;\Varid{f}\;\Varid{as})){}\<[E]%
\\[\blanklineskip]%
\>[3]{}\Varid{test}_{3.3}\;\Varid{n}\;\Varid{f}\;\mu\mathrel{=}\mathbf{do}\;\Varid{check}\;\Varid{n}\;\Varid{p}\;\mathbf{where}{}\<[E]%
\\
\>[3]{}\hsindent{2}{}\<[5]%
\>[5]{}\Varid{p}\;\Varid{as}\mathrel{=}\Varid{forAll}\;\Varid{as}\;(\lambda \Varid{a}\,\to\,(\Varid{f}\;\Varid{a}\in \Varid{min}_\Varid{u}\;\mu\;\Varid{f}\;\Varid{as})\Rightarrow (\Varid{a}\in \Varid{argmin}_\Varid{u}\;\mu\;\Varid{f}\;\Varid{as})){}\<[E]%
\ColumnHook
\end{hscode}\resethooks
As mentioned in the introduction, we neglect efficiency concerns in
the implementation of the tests for the sake of readability.

\subsection{Uncertainty functors: examples}
\label{subsection:ufexamples}

We discuss specific examples of uncertainty functors and, as done in
the previous section for \ensuremath{\Conid{U}\mathrel{=}\Conid{List}}, test the theory.
The tests suggest that measure functions fulfilling the monotonicity
condition (M2) from \cref{subsection:monotonicity} yield
implementations of \ensuremath{\Varid{min}_\Varid{u}} and \ensuremath{\Varid{argmin}_\Varid{u}} that satisfy
\crefrange{eq3.1}{eq3.3}.

\subsubsection{Identity functor}
\label{subsubsection:id}

The most obvious example of an uncertainty functor is the identity.
As already mentioned, we expect \ensuremath{\Varid{min}_\Varid{u}} and \ensuremath{\Varid{argmin}_\Varid{u}} to reduce
to \ensuremath{\Varid{min}} and \ensuremath{\Varid{argmin}} when \ensuremath{\Conid{U}\mathrel{=}\Conid{Id}}.
More formally, we require:

\begin{enumerate}[resume]
\item \ensuremath{\Varid{isSingletonSet}\;(\Varid{min}_\Varid{u}\;\mu\;\Varid{f}\;\Varid{as})},
\item \ensuremath{\Varid{min}_\Varid{u}\;\mu\;\Varid{f}\;\Varid{as}\,\doubleequals\,[\mskip1.5mu \Varid{min}\;\Varid{f}\;\Varid{as}\mskip1.5mu]},
\item \ensuremath{\Varid{argmin}_\Varid{u}\;\mu\;\Varid{f}\;\Varid{as}\,\doubleequals\,\Varid{argmin}\;\Varid{f}\;\Varid{as}}
\end{enumerate}

\noindent
for any \ensuremath{\Varid{f}}, non-empty \ensuremath{\Varid{as}}, and measure functions \ensuremath{\mu} that fulfil
(M2).
The equality in (3.5) and (3.6) is set equality.
First, we observe that \ensuremath{\Varid{const}\;\mathrm{7}} violates the singleton requirement:
\begin{hscode}\SaveRestoreHook
\column{B}{@{}>{\hspre}l<{\hspost}@{}}%
\column{3}{@{}>{\hspre}l<{\hspost}@{}}%
\column{5}{@{}>{\hspre}l<{\hspost}@{}}%
\column{E}{@{}>{\hspre}l<{\hspost}@{}}%
\>[3]{}\Varid{test}_{3.4}\;\Varid{n}\;\Varid{f}\;\mu\mathrel{=}\mathbf{do}\;\Varid{check}\;\Varid{n}\;(\Varid{p}\mathbin{::}[\mskip1.5mu \Int\mskip1.5mu]\,\to\,\Conid{Bool})\;\mathbf{where}{}\<[E]%
\\
\>[3]{}\hsindent{2}{}\<[5]%
\>[5]{}\Varid{p}\;\Varid{as}\mathrel{=}\neg \;(\Varid{isEmpty}\;\Varid{as})\Rightarrow \Varid{isSingletonSet}\;(\Varid{min}_\Varid{u}\;\mu\;\Varid{f}\;\Varid{as}){}\<[E]%
\ColumnHook
\end{hscode}\resethooks
\noindent
(e.g., for \ensuremath{\Varid{f}\;\Varid{x}\mathrel{=}\Conid{Id}\;(\Varid{x}\mathbin{*}\Varid{x})}).
Unsurprisingly, this measure also violates (M2) for the identity functor:
\begin{hscode}\SaveRestoreHook
\column{B}{@{}>{\hspre}l<{\hspost}@{}}%
\column{3}{@{}>{\hspre}l<{\hspost}@{}}%
\column{5}{@{}>{\hspre}l<{\hspost}@{}}%
\column{12}{@{}>{\hspre}c<{\hspost}@{}}%
\column{12E}{@{}l@{}}%
\column{15}{@{}>{\hspre}l<{\hspost}@{}}%
\column{20}{@{}>{\hspre}l<{\hspost}@{}}%
\column{E}{@{}>{\hspre}l<{\hspost}@{}}%
\>[3]{}\Varid{testM2}_{\Conid{Id}}\;\Varid{n}\;\mu\mathrel{=}\mathbf{do}\;\Varid{check}\;\Varid{n}\;(\Varid{p}\mathbin{::}\Conid{Id}\;\Int\,\to\,\Conid{Bool})\;\mathbf{where}{}\<[E]%
\\
\>[3]{}\hsindent{2}{}\<[5]%
\>[5]{}\Varid{p}\;\Varid{ub}_1{}\<[12]%
\>[12]{}\mathrel{=}{}\<[12E]%
\>[15]{}\mathbf{let}\;{}\<[20]%
\>[20]{}\Varid{ub}_2\mathrel{=}\Varid{fmap}\;(\mathrm{1}\mathbin{+})\;\Varid{ub}_1{}\<[E]%
\\
\>[15]{}\mathbf{in}\;{}\<[20]%
\>[20]{}(\Varid{ub}_1\mathrel{\prec_{u}}\Varid{ub}_2)\Rightarrow (\mu\;\Varid{ub}_1\mathbin{<}\mu\;\Varid{ub}_2){}\<[E]%
\ColumnHook
\end{hscode}\resethooks
By contrast, \ensuremath{\Varid{unwrap}} (with \ensuremath{\Varid{unwrap}\;(\Conid{Id}\;\Varid{x})\mathrel{=}\Varid{x}}) fulfils (M2) and
passes \ensuremath{\Varid{test}_{3.4}} (and similarly (3.5) and (3.6)).

\subsubsection{Probability functors}

Many decision problems involving stochastic uncertainty can be modeled
using probability functors.
The simplest case is when the set of possible outcomes is finite.
However, even when the theoretical outcome space \ensuremath{\Omega} is infinite
(e.g., a subset of \ensuremath{\Real^n}), the set of possible outcomes with non-zero
probability \ensuremath{\Varid{f}\;\Varid{a}} in a specific minimization problem is often finite.

\paragraph*{Simple probability.}

In this case, probability distributions can be represented by lists of
value-probability pairs \citep{10.1017/S0956796805005721,ionescu2009}
\begin{hscode}\SaveRestoreHook
\column{B}{@{}>{\hspre}l<{\hspost}@{}}%
\column{3}{@{}>{\hspre}l<{\hspost}@{}}%
\column{E}{@{}>{\hspre}l<{\hspost}@{}}%
\>[3]{}\mathbf{data}\;\Conid{SP}\;\Varid{b}\mathrel{=}\Conid{SP}\;[\mskip1.5mu (\Varid{b},\Real)\mskip1.5mu]{}\<[E]%
\ColumnHook
\end{hscode}\resethooks
and \ensuremath{\Conid{SP}} is a functor.
We can verify that standard measures like \ensuremath{\Varid{expVal}_{S\!P}}, \ensuremath{\Varid{best}_{S\!P}}, and \ensuremath{\Varid{worst}_{S\!P}}
fulfill the monotonicity condition (M2).
Note that to ensure valid test cases where \ensuremath{\Varid{ub}_1\mathrel{\prec_{u}}\Varid{ub}_2}, we must shift the values
by the spread of the support, similar to the list case:
\begin{hscode}\SaveRestoreHook
\column{B}{@{}>{\hspre}l<{\hspost}@{}}%
\column{3}{@{}>{\hspre}l<{\hspost}@{}}%
\column{5}{@{}>{\hspre}l<{\hspost}@{}}%
\column{14}{@{}>{\hspre}l<{\hspost}@{}}%
\column{19}{@{}>{\hspre}l<{\hspost}@{}}%
\column{26}{@{}>{\hspre}l<{\hspost}@{}}%
\column{E}{@{}>{\hspre}l<{\hspost}@{}}%
\>[3]{}\Varid{testM2}_{S\!P}\;\Varid{n}\;\mu\mathrel{=}\mathbf{do}\;\Varid{check}\;\Varid{n}\;(\Varid{p}\mathbin{::}\Conid{SP}\;\Int\,\to\,\Conid{Bool})\;\mathbf{where}{}\<[E]%
\\
\>[3]{}\hsindent{2}{}\<[5]%
\>[5]{}\Varid{p}\;\Varid{ub}_1\mathrel{=}{}\<[14]%
\>[14]{}\mathbf{let}\;{}\<[19]%
\>[19]{}\Varid{vals}{}\<[26]%
\>[26]{}\mathrel{=}\Varid{elemsSP}\;\Varid{ub}_1{}\<[E]%
\\
\>[19]{}\Varid{shift}{}\<[26]%
\>[26]{}\mathrel{=}\mathrm{1}\mathbin{+}(\Varid{maximum}\;\Varid{vals}\mathbin{-}\Varid{minimum}\;\Varid{vals}){}\<[E]%
\\
\>[19]{}\Varid{ub}_2{}\<[26]%
\>[26]{}\mathrel{=}\Varid{fmap}\;(\Varid{shift}\mathbin{+})\;\Varid{ub}_1{}\<[E]%
\\
\>[14]{}\mathbf{in}\;{}\<[19]%
\>[19]{}(\Varid{isSP}\;\Varid{ub}_1\mathrel{\wedge}\Varid{isSP}\;\Varid{ub}_2)\Rightarrow (\Varid{ub}_1\mathrel{\prec_{u}}\Varid{ub}_2\Rightarrow (\mu\;\Varid{ub}_1\mathbin{<}\mu\;\Varid{ub}_2)){}\<[E]%
\ColumnHook
\end{hscode}\resethooks
With these measures, \ensuremath{\Varid{min}_\Varid{u}} passes the structural specifications
(\ensuremath{\Varid{test}_{3.1}}, \ensuremath{\Varid{test}_{3.2}} and \ensuremath{\Varid{test}_{3.3}}).
For example, with:
\begin{hscode}\SaveRestoreHook
\column{B}{@{}>{\hspre}l<{\hspost}@{}}%
\column{3}{@{}>{\hspre}l<{\hspost}@{}}%
\column{14}{@{}>{\hspre}c<{\hspost}@{}}%
\column{14E}{@{}l@{}}%
\column{17}{@{}>{\hspre}l<{\hspost}@{}}%
\column{E}{@{}>{\hspre}l<{\hspost}@{}}%
\>[3]{}\Varid{f}_{S\!P}\mathbin{::}\Conid{Bool}\,\to\,\Conid{SP}\;\Int{}\<[E]%
\\
\>[3]{}\Varid{f}_{S\!P}\;\Conid{False}{}\<[14]%
\>[14]{}\mathrel{=}{}\<[14E]%
\>[17]{}\Conid{SP}\;[\mskip1.5mu (\mathrm{2},\mathrm{0.4}),(\mathrm{0},\mathrm{0.3}),(\mathrm{1},\mathrm{0.3})\mskip1.5mu]{}\<[E]%
\\
\>[3]{}\Varid{f}_{S\!P}\;\Conid{True}{}\<[14]%
\>[14]{}\mathrel{=}{}\<[14E]%
\>[17]{}\Conid{SP}\;[\mskip1.5mu (\mathrm{3},\mathrm{0.4}),(\mathrm{4},\mathrm{0.3}),(\mathrm{4},\mathrm{0.3})\mskip1.5mu]{}\<[E]%
\ColumnHook
\end{hscode}\resethooks
By contrast, constant measure functions violate (M2) and fail \cref{eq3.1}.

\paragraph*{Intervals of real numbers.}

In many physical applications, \ensuremath{\Varid{f}} maps inputs to intervals of real numbers:
\begin{hscode}\SaveRestoreHook
\column{B}{@{}>{\hspre}l<{\hspost}@{}}%
\column{3}{@{}>{\hspre}l<{\hspost}@{}}%
\column{E}{@{}>{\hspre}l<{\hspost}@{}}%
\>[3]{}\mathbf{data}\;\Conid{I}\;\Varid{b}\mathrel{=}\Conid{I}\;(\Varid{b},\Varid{b}){}\<[E]%
\ColumnHook
\end{hscode}\resethooks
\noindent
The interpretation of \ensuremath{\Varid{f}\;\Varid{a}\mathrel{=}\Conid{I}\;(x_s,x_e)} is that all values between
\ensuremath{x_s} and \ensuremath{x_e} are possible.
We define \ensuremath{(\mathrel{\prec_{u}})} and \ensuremath{(\mathrel{\leq\!\cdot})} for intervals via their endpoints:
\begin{hscode}\SaveRestoreHook
\column{B}{@{}>{\hspre}l<{\hspost}@{}}%
\column{3}{@{}>{\hspre}l<{\hspost}@{}}%
\column{15}{@{}>{\hspre}c<{\hspost}@{}}%
\column{15E}{@{}l@{}}%
\column{16}{@{}>{\hspre}c<{\hspost}@{}}%
\column{16E}{@{}l@{}}%
\column{20}{@{}>{\hspre}l<{\hspost}@{}}%
\column{34}{@{}>{\hspre}c<{\hspost}@{}}%
\column{34E}{@{}l@{}}%
\column{37}{@{}>{\hspre}l<{\hspost}@{}}%
\column{47}{@{}>{\hspre}c<{\hspost}@{}}%
\column{47E}{@{}l@{}}%
\column{48}{@{}>{\hspre}c<{\hspost}@{}}%
\column{48E}{@{}l@{}}%
\column{52}{@{}>{\hspre}l<{\hspost}@{}}%
\column{E}{@{}>{\hspre}l<{\hspost}@{}}%
\>[3]{}\Conid{I}\;(x_s,x_e){}\<[16]%
\>[16]{}\mathrel{\prec_{u}}{}\<[16E]%
\>[20]{}\Conid{I}\;(\Varid{x'}\!_s,\Varid{x'}\!_e){}\<[34]%
\>[34]{}\mathrel{=}{}\<[34E]%
\>[37]{}[\mskip1.5mu x_s,x_e\mskip1.5mu]{}\<[48]%
\>[48]{}\mathrel{\prec_{u}}{}\<[48E]%
\>[52]{}[\mskip1.5mu \Varid{x'}\!_s,\Varid{x'}\!_e\mskip1.5mu]{}\<[E]%
\\
\>[3]{}\Conid{I}\;(x_s,x_e){}\<[15]%
\>[15]{}\mathrel{\leq\!\cdot}{}\<[15E]%
\>[20]{}\Conid{I}\;(\Varid{x'}\!_s,\Varid{x'}\!_e){}\<[34]%
\>[34]{}\mathrel{=}{}\<[34E]%
\>[37]{}[\mskip1.5mu x_s,x_e\mskip1.5mu]{}\<[47]%
\>[47]{}\mathrel{\leq\!\cdot}{}\<[47E]%
\>[52]{}[\mskip1.5mu \Varid{x'}\!_s,\Varid{x'}\!_e\mskip1.5mu]{}\<[E]%
\ColumnHook
\end{hscode}\resethooks
\noindent
Notice that \ensuremath{i_1\mathrel{\prec_{u}}i_2} implies that the two intervals are disjoint
while \ensuremath{i_1\mathrel{\leq\!\cdot}i_2} allows for \ensuremath{i_1} and \ensuremath{i_2} to overlap.
Measures like \ensuremath{\Varid{sum}_{I}} (sum of endpoints), \ensuremath{\Varid{average}_{I}} (midpoint),
\ensuremath{\Varid{best}_{I}} (lower end), and \ensuremath{\Varid{worst}_{I}} (upper end of the interval) fulfil
the monotonicity condition (M2).
We verify this using a shift logic ensuring disjointness:
\begin{hscode}\SaveRestoreHook
\column{B}{@{}>{\hspre}l<{\hspost}@{}}%
\column{3}{@{}>{\hspre}l<{\hspost}@{}}%
\column{5}{@{}>{\hspre}l<{\hspost}@{}}%
\column{13}{@{}>{\hspre}l<{\hspost}@{}}%
\column{17}{@{}>{\hspre}l<{\hspost}@{}}%
\column{24}{@{}>{\hspre}l<{\hspost}@{}}%
\column{E}{@{}>{\hspre}l<{\hspost}@{}}%
\>[3]{}\Varid{testM2}_{I}\;\Varid{n}\;\mu\mathrel{=}\mathbf{do}\;\Varid{check}\;\Varid{n}\;(\Varid{p}\mathbin{::}\Conid{I}\;\Real\,\to\,\Conid{Bool})\;\mathbf{where}{}\<[E]%
\\
\>[3]{}\hsindent{2}{}\<[5]%
\>[5]{}\Varid{p}\;\Varid{ub}_1\mathrel{=}\mathbf{let}\;\Varid{shift}{}\<[24]%
\>[24]{}\mathrel{=}\mathrm{1}\mathbin{+}(\Varid{end}\;\Varid{ub}_1\mathbin{-}\Varid{start}\;\Varid{ub}_1){}\<[E]%
\\
\>[5]{}\hsindent{12}{}\<[17]%
\>[17]{}\Varid{ub}_2{}\<[24]%
\>[24]{}\mathrel{=}\Varid{fmap}\;(\Varid{shift}\mathbin{+})\;\Varid{ub}_1{}\<[E]%
\\
\>[5]{}\hsindent{8}{}\<[13]%
\>[13]{}\mathbf{in}\;(\Varid{is}_{I}\;\Varid{ub}_1)\Rightarrow (\Varid{ub}_1\mathrel{\prec_{u}}\Varid{ub}_2\Rightarrow (\mu\;\Varid{ub}_1\mathbin{<}\mu\;\Varid{ub}_2)){}\<[E]%
\ColumnHook
\end{hscode}\resethooks
However, measures like \ensuremath{\Varid{width}_{I}} (the size of the uncertainty) or \ensuremath{\Varid{const}} violate (M2) and cause \ensuremath{\Varid{min}_\Varid{u}} to fail.
For example:
\begin{hscode}\SaveRestoreHook
\column{B}{@{}>{\hspre}l<{\hspost}@{}}%
\column{3}{@{}>{\hspre}l<{\hspost}@{}}%
\column{13}{@{}>{\hspre}c<{\hspost}@{}}%
\column{13E}{@{}l@{}}%
\column{16}{@{}>{\hspre}l<{\hspost}@{}}%
\column{31}{@{}>{\hspre}l<{\hspost}@{}}%
\column{E}{@{}>{\hspre}l<{\hspost}@{}}%
\>[3]{}\Varid{f}_{I}\mathbin{::}\Conid{Bool}\,\to\,\Conid{I}\;\Real{}\<[E]%
\\
\>[3]{}\Varid{f}_{I}\;\Conid{False}{}\<[13]%
\>[13]{}\mathrel{=}{}\<[13E]%
\>[16]{}\Conid{I}\;(\mathrm{1.1},\mathrm{1.2}){}\<[31]%
\>[31]{}\mbox{\onelinecomment  Width 0.1}{}\<[E]%
\\
\>[3]{}\Varid{f}_{I}\;\Conid{True}{}\<[13]%
\>[13]{}\mathrel{=}{}\<[13E]%
\>[16]{}\Conid{I}\;(\mathrm{0.0},\mathrm{1.0}){}\<[31]%
\>[31]{}\mbox{\onelinecomment  Width 1.0}{}\<[E]%
\ColumnHook
\end{hscode}\resethooks
If minimizing \ensuremath{\Varid{width}_{I}}, \ensuremath{\Conid{False}} is preferred, even though \ensuremath{\Conid{True}} offers strictly lower values.
This violates \cref{eq3.1}.

\paragraph*{Probability density functions (Histogram approximation).}

A common way to model continuous uncertainty is via a \emph{probability density function} (PDF).
Here, we consider a discrete approximation where the PDF is represented by a finite number of equally sized ``bins'' (a histogram) over a support interval:
\begin{hscode}\SaveRestoreHook
\column{B}{@{}>{\hspre}l<{\hspost}@{}}%
\column{3}{@{}>{\hspre}l<{\hspost}@{}}%
\column{E}{@{}>{\hspre}l<{\hspost}@{}}%
\>[3]{}\mathbf{data}\;\Conid{PDF}\;\Varid{b}\mathrel{=}\Conid{PDF}\;(\Conid{I}\;\Varid{b},[\mskip1.5mu \Real\mskip1.5mu]){}\<[E]%
\ColumnHook
\end{hscode}\resethooks
\DONE{This definition seems very far off from the earlier cases and is
  only one (rather unusual) way to model pdfs. Nicola: why far off? I
  thought that approximating PDFs with a finite number of bins was
  quite standard. It also seems a natural way of applying finite
  probability distributions to the continuous case. Patrik: I tried to
  make it just one possibility now.}
\noindent
The structural comparison of PDFs is inherited strictly from their
support intervals (ignoring weights for dominance purposes):
\ensuremath{\Varid{pdf1}\mathrel{\prec_{u}}\Varid{pdf2}} iff \ensuremath{\Varid{interval}_{P\!D\!F}\;\Varid{pdf1}\mathrel{\prec_{u}}\Varid{interval}_{P\!D\!F}\;\Varid{pdf2}}.
We verify that \ensuremath{\Varid{expVal}_{P\!D\!F}} fulfills (M2):
\begin{hscode}\SaveRestoreHook
\column{B}{@{}>{\hspre}l<{\hspost}@{}}%
\column{3}{@{}>{\hspre}l<{\hspost}@{}}%
\column{5}{@{}>{\hspre}l<{\hspost}@{}}%
\column{14}{@{}>{\hspre}l<{\hspost}@{}}%
\column{19}{@{}>{\hspre}l<{\hspost}@{}}%
\column{25}{@{}>{\hspre}l<{\hspost}@{}}%
\column{E}{@{}>{\hspre}l<{\hspost}@{}}%
\>[3]{}\Varid{testM2}_{P\!D\!F}\;\Varid{n}\;\mu\mathrel{=}\mathbf{do}\;\Varid{check}\;\Varid{n}\;(\Varid{p}\mathbin{::}\Conid{PDF}\;\Real\,\to\,\Conid{Bool})\;\mathbf{where}{}\<[E]%
\\
\>[3]{}\hsindent{2}{}\<[5]%
\>[5]{}\Varid{p}\;\Varid{ub}_1\mathrel{=}{}\<[14]%
\>[14]{}\mathbf{let}\;{}\<[19]%
\>[19]{}\Varid{shift}\mathrel{=}\mathrm{1}\mathbin{+}\Varid{width}_{I}\;(\Varid{interval}_{P\!D\!F}\;\Varid{ub}_1){}\<[E]%
\\
\>[19]{}\Varid{ub}_2{}\<[25]%
\>[25]{}\mathrel{=}\Varid{fmap}\;(\Varid{shift}\mathbin{+})\;\Varid{ub}_1{}\<[E]%
\\
\>[14]{}\mathbf{in}\;{}\<[19]%
\>[19]{}(\Varid{is}_{P\!D\!F}\;\Varid{ub}_1)\Rightarrow (\Varid{ub}_1\mathrel{\prec_{u}}\Varid{ub}_2\Rightarrow (\mu\;\Varid{ub}_1\mathbin{<}\mu\;\Varid{ub}_2)){}\<[E]%
\ColumnHook
\end{hscode}\resethooks
\begin{hscode}\SaveRestoreHook
\column{B}{@{}>{\hspre}l<{\hspost}@{}}%
\column{3}{@{}>{\hspre}l<{\hspost}@{}}%
\column{15}{@{}>{\hspre}c<{\hspost}@{}}%
\column{15E}{@{}l@{}}%
\column{18}{@{}>{\hspre}l<{\hspost}@{}}%
\column{E}{@{}>{\hspre}l<{\hspost}@{}}%
\>[3]{}\Varid{f}_{P\!D\!F}\mathbin{::}\Conid{Bool}\,\to\,\Conid{PDF}\;\Real{}\<[E]%
\\
\>[3]{}\Varid{f}_{P\!D\!F}\;\Conid{False}{}\<[15]%
\>[15]{}\mathrel{=}{}\<[15E]%
\>[18]{}\Conid{PDF}\;(\Conid{I}\;(\mathrm{1.0},\mathrm{2.0}),[\mskip1.5mu \mathrm{0.6},\mathrm{0.4}\mskip1.5mu]){}\<[E]%
\\
\>[3]{}\Varid{f}_{P\!D\!F}\;\Conid{True}{}\<[15]%
\>[15]{}\mathrel{=}{}\<[15E]%
\>[18]{}\Conid{PDF}\;(\Conid{I}\;(\mathrm{0.0},\mathrm{0.1}),[\mskip1.5mu \mathrm{0.6},\mathrm{0.4}\mskip1.5mu]){}\<[E]%
\ColumnHook
\end{hscode}\resethooks
\noindent
As usual, constant measure functions are a valid counterexample.
Interestingly, \ensuremath{\Varid{mostLikely}_{P\!D\!F}} (the midpoint of the highest bin)
fulfills (M1), unlike in the simple probability case where
\ensuremath{\Varid{mostLikely}_{S\!P}} violates (M1), see \citet{IONESCU_2016b}.

It is worth noticing that \ensuremath{\Varid{mostLikely}_{S\!P}} fulfils (M2) (if every element
of \ensuremath{ub_1} is strictly smaller than every element of \ensuremath{ub_2}, then the most
likely element of \ensuremath{ub_1} is strictly smaller than the most likely element
of \ensuremath{ub_2}) and hence (M2) does not imply (M1) as argued in
\cref{subsection:monotonicity}.
%

%
%
The fact that popular measures (like \ensuremath{\Varid{mostLikely}}) satisfy
monotonicity for some functors (\ensuremath{\Conid{PDF}}) but fail for others (\ensuremath{\Conid{SP}})
underscores the need for the generic theory presented here.

\subsection{Uncertainty functors in general}
\label{subsection:ufgeneral}

So far, we have discussed inequalities in \ensuremath{\Conid{U}\;\Conid{B}} and monotonicity
conditions for optimization under functorial (epistemic, aleatoric and
mixed) uncertainty.
In \cref{subsection:ufexamples}, we discussed examples of uncertainty
functors relevant for applications and argued that measure functions
fulfilling the monotonicity condition (M2) yield implementations of
\ensuremath{\Varid{min}_\Varid{u}} and \ensuremath{\Varid{argmin}_\Varid{u}} that satisfy \crefrange{eq3.1}{eq3.3}.

In doing so, we implicitly characterized the notion of an uncertainty
functor via an equality test for shape and two inequalities:
\begin{hscode}\SaveRestoreHook
\column{B}{@{}>{\hspre}l<{\hspost}@{}}%
\column{3}{@{}>{\hspre}l<{\hspost}@{}}%
\column{5}{@{}>{\hspre}l<{\hspost}@{}}%
\column{12}{@{}>{\hspre}c<{\hspost}@{}}%
\column{12E}{@{}l@{}}%
\column{16}{@{}>{\hspre}l<{\hspost}@{}}%
\column{E}{@{}>{\hspre}l<{\hspost}@{}}%
\>[3]{}\mathbf{class}\;\Conid{Functor}\;\Varid{u}\Rightarrow \Conid{UncertaintyFunctor}\;\Varid{u}\;\mathbf{where}{}\<[E]%
\\
\>[3]{}\hsindent{2}{}\<[5]%
\>[5]{}(\thickapprox){}\<[12]%
\>[12]{}\mathbin{::}{}\<[12E]%
\>[16]{}\Varid{u}\;\Varid{b}\,\to\,\Varid{u}\;\Varid{b}\,\to\,\Conid{Bool}{}\<[E]%
\\
\>[3]{}\hsindent{2}{}\<[5]%
\>[5]{}(\mathrel{\prec_{u}}){}\<[12]%
\>[12]{}\mathbin{::}{}\<[12E]%
\>[16]{}\Conid{Ord}\;\Varid{b}\Rightarrow \Varid{u}\;\Varid{b}\,\to\,\Varid{u}\;\Varid{b}\,\to\,\Conid{Bool}{}\<[E]%
\\
\>[3]{}\hsindent{2}{}\<[5]%
\>[5]{}(\mathrel{\leq\!\cdot}){}\<[12]%
\>[12]{}\mathbin{::}{}\<[12E]%
\>[16]{}\Conid{Ord}\;\Varid{b}\Rightarrow \Varid{u}\;\Varid{b}\,\to\,\Varid{u}\;\Varid{b}\,\to\,\Conid{Bool}{}\<[E]%
\ColumnHook
\end{hscode}\resethooks
Recall from \cref{subsection:monotonicity} that \ensuremath{(\mathrel{\leq\!\cdot})} was only
defined for structures of the same shape (and for \ensuremath{\Conid{U}\mathrel{=}\Conid{List}}, shape
corresponds to length).
However, we have not yet discussed the minimal requirements that allow
one to define uncertainty functors generically, nor how this notion
relates to well-understood abstractions in functional programming.

It is worth noting that the definitions of \ensuremath{(\mathrel{\prec_{u}})} and \ensuremath{(\mathrel{\leq\!\cdot})}
discussed above can be formulated entirely in terms of membership
tests, quantifiers, and (for \ensuremath{(\mathrel{\leq\!\cdot})}) a \ensuremath{\Varid{zip}} combinator.
This suggests that a fruitful way of understanding uncertainty
functors is as \emph{functors that preserve decidable equality}.
Intuitively, if we can decide equality for elements of type \ensuremath{\Conid{A}}, we
should be able to decide equality for the structure \ensuremath{\Conid{F}\;\Conid{A}}.
In the context of constructive type theory, this property is formulated
precisely as a predicate on type constructors:
\begin{hscode}\SaveRestoreHook
\column{B}{@{}>{\hspre}l<{\hspost}@{}}%
\column{3}{@{}>{\hspre}l<{\hspost}@{}}%
\column{16}{@{}>{\hspre}c<{\hspost}@{}}%
\column{16E}{@{}l@{}}%
\column{19}{@{}>{\hspre}l<{\hspost}@{}}%
\column{E}{@{}>{\hspre}l<{\hspost}@{}}%
\>[3]{}\Conid{PresDec}{\equiv}\,\;\Conid{F}{}\<[16]%
\>[16]{}\mathrel{=}{}\<[16E]%
\>[19]{}(\Conid{A}\ \mathop{:}\ \Conid{Set})\,\to\,\Conid{Dec}{\equiv}\,\;\Conid{A}\,\to\,\Conid{Dec}{\equiv}\,\;(\Conid{F}\;\Conid{A}){}\<[E]%
\ColumnHook
\end{hscode}\resethooks
\noindent
where \ensuremath{\Conid{Dec}{\equiv}\,\;\Conid{A}} witnesses that equality on \ensuremath{\Conid{A}} is decidable.
A dependently typed formalization of this concept in Agda
\citep{norell2007thesis} is available in the accompanying repository
\citep{MOOcodeRepo2025b}.
Additionally, we require a \ensuremath{\Varid{zip}} operation
\begin{hscode}\SaveRestoreHook
\column{B}{@{}>{\hspre}l<{\hspost}@{}}%
\column{3}{@{}>{\hspre}l<{\hspost}@{}}%
\column{E}{@{}>{\hspre}l<{\hspost}@{}}%
\>[3]{}\Varid{zip}\mathbin{::}(\Conid{F}\;\Conid{A},\Conid{F}\;\Conid{B})\,\to\,\Conid{F}\;(\Conid{A},\Conid{B}){}\<[E]%
\ColumnHook
\end{hscode}\resethooks


\noindent
Perhaps not surprisingly, \ensuremath{\Conid{List}} preserves decidable equality (and has a
\ensuremath{\Varid{zip}} operation), whereas \ensuremath{\Conid{Reader}} (over infinite domains) does not.
For functors that preserve decidable equality, one can easily define
an equality test for shapes, a membership test, and canonical
quantifiers generically.
For example, using \ensuremath{\Conid{List}} syntax to illustrate the generic concept, we
can test the equality of shapes:
\begin{hscode}\SaveRestoreHook
\column{B}{@{}>{\hspre}l<{\hspost}@{}}%
\column{3}{@{}>{\hspre}l<{\hspost}@{}}%
\column{E}{@{}>{\hspre}l<{\hspost}@{}}%
\>[3]{}\Varid{shapeList}\mathbin{::}[\mskip1.5mu \Varid{a}\mskip1.5mu]\,\to\,[\mskip1.5mu ()\mskip1.5mu]{}\<[E]%
\\
\>[3]{}\Varid{shapeList}\mathrel{=}\Varid{fmap}\;(\Varid{const}\;()){}\<[E]%
\\[\blanklineskip]%
\>[3]{}\Varid{sameShapeList}\mathbin{::}[\mskip1.5mu \Varid{a}\mskip1.5mu]\,\to\,[\mskip1.5mu \Varid{a}\mskip1.5mu]\,\to\,\Conid{Bool}{}\<[E]%
\\
\>[3]{}\Varid{sameShapeList}\;\Varid{xs}\;\Varid{ys}\mathrel{=}\Varid{shapeList}\;\Varid{xs}\doubleequals\Varid{shapeList}\;\Varid{ys}{}\<[E]%
\ColumnHook
\end{hscode}\resethooks
\noindent
without invoking \ensuremath{\Varid{length}} or the constructors of \ensuremath{\Conid{List}}.
Similarly, we can define membership:
\begin{hscode}\SaveRestoreHook
\column{B}{@{}>{\hspre}l<{\hspost}@{}}%
\column{3}{@{}>{\hspre}l<{\hspost}@{}}%
\column{E}{@{}>{\hspre}l<{\hspost}@{}}%
\>[3]{}\Varid{elemList}\mathbin{::}\Conid{Eq}\;\Varid{a}\Rightarrow \Varid{a}\,\to\,[\mskip1.5mu \Varid{a}\mskip1.5mu]\,\to\,\Conid{Bool}{}\<[E]%
\\
\>[3]{}\Varid{elemList}\;\Varid{a}\;\Varid{as}\mathrel{=}\neg \;(\Varid{fmap}\;(\Varid{a}\doubleequals)\;\Varid{as}\doubleequals\Varid{fmap}\;(\Varid{const}\;\Conid{False})\;\Varid{as}){}\<[E]%
\ColumnHook
\end{hscode}\resethooks
\noindent
and quantifiers:
\begin{hscode}\SaveRestoreHook
\column{B}{@{}>{\hspre}l<{\hspost}@{}}%
\column{3}{@{}>{\hspre}l<{\hspost}@{}}%
\column{E}{@{}>{\hspre}l<{\hspost}@{}}%
\>[3]{}\Varid{allList}\mathbin{::}[\mskip1.5mu \Conid{Bool}\mskip1.5mu]\,\to\,\Conid{Bool}{}\<[E]%
\\
\>[3]{}\Varid{allList}\;\Varid{as}\mathrel{=}\Varid{as}\doubleequals\Varid{fmap}\;(\Varid{const}\;\Conid{True})\;\Varid{as}{}\<[E]%
\\[\blanklineskip]%
\>[3]{}\Varid{anyList}\mathbin{::}[\mskip1.5mu \Conid{Bool}\mskip1.5mu]\,\to\,\Conid{Bool}{}\<[E]%
\\
\>[3]{}\Varid{anyList}\;\Varid{as}\mathrel{=}\Varid{as}\not\doubleequals\Varid{fmap}\;(\Varid{const}\;\Conid{False})\;\Varid{as}{}\<[E]%
\ColumnHook
\end{hscode}\resethooks
\noindent
without resorting to implementation-specific pattern matching.
These constructions hold for \emph{any} functor that preserves
decidable equality.
Consequently, for such functors (provided they also have \ensuremath{\Varid{zip}}),
one can define the operations that characterize uncertainty functors
generically; see \citep{MOOcodeRepo2025b} for the full derivation.

We therefore propose that uncertainty functors are best understood as
functors that preserve decidable equality and possess a suitable
\ensuremath{\Varid{zip}}.
We also observe that all our finite examples of uncertainty functors
are \emph{Traversable} \citep{mcbride2008applicative}.
While the precise relationship between traversability and the
preservation of decidable equality is beyond the scope of this paper,
we consider it a promising avenue for future work.

\section{Related work}
\label{section:related}

\paragraph*{Classical multi-objective optimization.}

Optimization under uncertainty is a mature research area and numerical
libraries for multi-objective optimization are available in, among
others, C, C++, Python, Julia, R, C\#, Matlab and Java, see Table 1 in
\citep{Emmerich+Deutz}.
Classical approaches include scalarization, $\varepsilon$-constraint
methods, and evolutionary algorithms \citep{NSGAII, MOEAD}.

Most libraries are focused on the minimization of functions of real
variables and MOO algorithms are often presented and implemented in
imperative form. The lack of formal specifications limits the
possibility of program testing and verification
\citep{ionescujansson:LIPIcs:2013:3899} to benchmarking against selected
problems like, for example, the one sketched in
\cref{fig:MOObenchmark250k,fig:BenchmarkMOOfunctions} of
\cref{subsection:applications}.

\paragraph*{Economic foundations.}

In microeconomic terms \citep{varian1992}, the strict dominance relation \ensuremath{(\prec)}
used in \cref{section:moo} corresponds to \emph{strong Pareto efficiency}.
Consequently, the computational task of finding the Pareto front is
equivalent to identifying the set of feasible allocations that are
strong Pareto efficient.

\paragraph*{Functional programming and explicit specifications.}

By contrast, the functional programming approach emphasizes correctness
and systematic property-based testing~\citep{Hughes2010} and the methods
discussed in \cref{section:moo} are a first step towards overcoming
implicit specifications and standard black-box heuristics in MOO.
In particular, our operational definition of the Pareto front via the
\ensuremath{\Varid{bump}} function corresponds to the \emph{thinning} strategy described
by \citet{bird_algebra_1997} and later refined by
\citet{mu2009algebra} and used in \citep{jansson_jansson_2023}.
This connects our approach to the broader field of algorithm
derivation, where optimization problems are solved by filtering
(thinning) results to maintain only those that are optimal according
to a preorder.
By contrast, the functional programming approach emphasizes correctness
and systematic property-based testing~\citep{Hughes2010} and the methods
discussed in \cref{section:moo} are a first step towards overcoming
implicit specifications and standard black-box heuristics in MOO.

\paragraph*{Classical models of uncertainty.}

Classical treatments of optimization under \emph{stochastic} uncertainty
include stochastic programming~\citep{StochProg} and robust
optimization~\citep{RobustOpt}.
Info-gap decision theory and dynamic adaptive planning offer alternative
frameworks~\citep{InfoGap, DAPP}.

\paragraph*{Functorial representations of uncertainty.}

The idea of representing uncertainty functorially and, in the aleatoric
case, monadically is well established in functional programming
\citep{giry1981,10.1017/S0956796805005721,ionescu2009,IONESCU_2016b,2017_Botta_Jansson_Ionescu}
and \cref{section:monadic} is a first step towards building a theory of
optimization under functorial uncertainty.
By only assuming \ensuremath{\Conid{U}} to be an uncertainty functor (see
\cref{subsection:ufgeneral}), the theory allows one to deal seamlessly
with different kinds of computational uncertainty: non-deterministic,
stochastic, fuzzy, etc. \citep{ionescu2009,IONESCU_2016b}.
At the semantic level, the theory supports accounting for different
attitudes with respect to risk under epistemic and aleatoric uncertainty
\citep{shepherd+al2018, shepherd2019} in an accountable way through
different measure functions and systematic property-based
testing~\citep{Hughes2010} of their monotonicity properties
\cref{subsection:monotonicity}.

\paragraph*{Complementary approaches.}

Other approaches towards dealing with uncertainty include probabilistic
programming languages \citep{ProbProg}. Abstract interpretation and
reachability provide sound but coarse
over-approximations~\citep{AbstractInterp, Reachability} and safe
reinforcement learning and constrained model predictive control tackle
uncertainty in online adaptive settings~\citep{SafeRL, MPC}. The work
presented here is complementary, targeting offline policy design with
formal guarantees.

\section{Discussion and Conclusion}
\label{section:conclusion}


\paragraph*{Summary of contributions.}

We have generalized well understood, generic \ensuremath{\Varid{min}} and \ensuremath{\Varid{argmin}}
algorithms to the case in which optimization is affected by value
uncertainty and by functorial uncertainty.
The resulting algorithms are easy to understand and to test and yet they
can be applied to tackle complex decision problems, as demonstrated in
\citep{Pusztai2023BayesOptMMI, jansson+2025}.
\DONE{Cite VCCS and other papers from the Chalmers plasma physics
  group}

Uncertainty functors represent epistemic uncertainty, aleatoric
uncertainty or more general kinds of uncertainty. In
\cref{section:monadic} we have build upon the framework originally
developed by Ionescu for modelling vulnerability with functional
programming and dependent types \citep{IONESCU_2016b} to derive a method
to account for these uncertainties in optimization problems.
In \cref{subsection:ufgeneral} we have discussed how uncertainty
functors relate to well understood abstractions in mathematics and
functional programming and in \cref{subsection:ufexamples} we have
provided examples of uncertainty functors that are relevant for
formalizing and understanding problems in physics, engineering and
(climate) modelling.

For these kind of problems and especially for decision making in the
context of climate impact research \citep{webster2000, webster2008,
  CARLINO202016593, shape2021, Botta2023MatterMost}, it is crucial to
avoid ad-hoc approaches and to apply consistent measure functions.
Where formal verification is not viable
\citep{ionescujansson:LIPIcs:2013:3899}, the methods should be amenable
to extensive property based testing via QuickCheck \citep{Hughes2010},
as demonstrated in \cref{subsection:ufexamples}. This requires precise
mathematical specifications.

The ones presented here extend those discussed in
\citep{2014_Botta_et_al,2017_Botta_Jansson_Ionescu,BREDE_BOTTA_2021}
for which we have derived verified implementations in Idris
\citep{idris2referenceguide}. Thus, our work can be seen as a first step
towards building generic, ``increasingly correct'' methods for
optimisation under value and functorial uncertainty following the
approach discussed in \citep{ionescujansson:LIPIcs:2013:3899}.

\paragraph*{Interpreting combined uncertainty in climate policy.}

In most engineering applications but also in (climate) integrated
assessment one typically has to do with decision problems which are
affected by both value uncertainty \citep{CARLINO202016593} and
by functorial uncertainty.
In the rest of this section, we sketch how the methods discussed in
\cref{section:moo,section:monadic} can be applied to these problems.

Climate policy offers a simple yet paradigmatic example.
In policy advice, the goal is to identify mitigation and adaptation
measures for the next few decades that satisfy two conflicting
objectives.
First, they must be likely to avoid unmanageable short- and long-term
damages from climate change (for example, on the global economy
\citep{TOL2024113922, KotzLevermannWenz2024, NealNewellPitman2025,
  BearparkHoganHsiang2025, Schoetz2025, Nordhaus12261}).
Second, they must keep the costs and societal risks of decarbonizing
fossil-fuel-based economies acceptable.
The problem (of selecting an optimal level of short-term investments
or costs) is affected by value uncertainty because there is a
trade-off between minimizing damages and short-term costs. It is
affected by functorial uncertainty because decisions about short-term
climate policies (for example, investments in mitigation and
adaptation policies from now until the year 2100) do not uniquely
determine damages from climate change.
This arises for at least three reasons:
\begin{itemize}
\item First, the relationship between short-term costs and actual
  short-term greenhouse gas (GHG) emissions reductions is uncertain:
  investments may turn out to be more or less effective, and decisions
  about such investments may not get implemented or may be implemented
  with delays.
\item Second, our knowledge of how GHG emissions impact the evolution
  of the climate system and this, in turn, economic damages, is far
  from perfect.
\item Third, decisions about climate policies until the year 2100 do
  not determine later GHG emissions uniquely, yet long-term damages
  depend on these later emissions.
\end{itemize}

This means that, given short-term costs (investments in adaptation and
mitigation options) \ensuremath{\Varid{c}\ \mathop{:}\ C_{\varphi}} in a \emph{feasible} range \ensuremath{C_{\varphi}}, the
associated \emph{possible} damages \ensuremath{\Varid{poss}\;\Varid{c}\ \mathop{:}\ \Conid{U}\;\Conid{D}} are characterized by
an uncertainty functor \ensuremath{\Conid{U}}, perhaps \ensuremath{\Conid{U}\mathrel{=}\Conid{I}} or \ensuremath{\Conid{U}\mathrel{=}\Conid{PDF}} as discussed in
\cref{subsection:ufexamples}.

\begin{figure}[h]
  \includegraphics[width=0.8\textwidth]{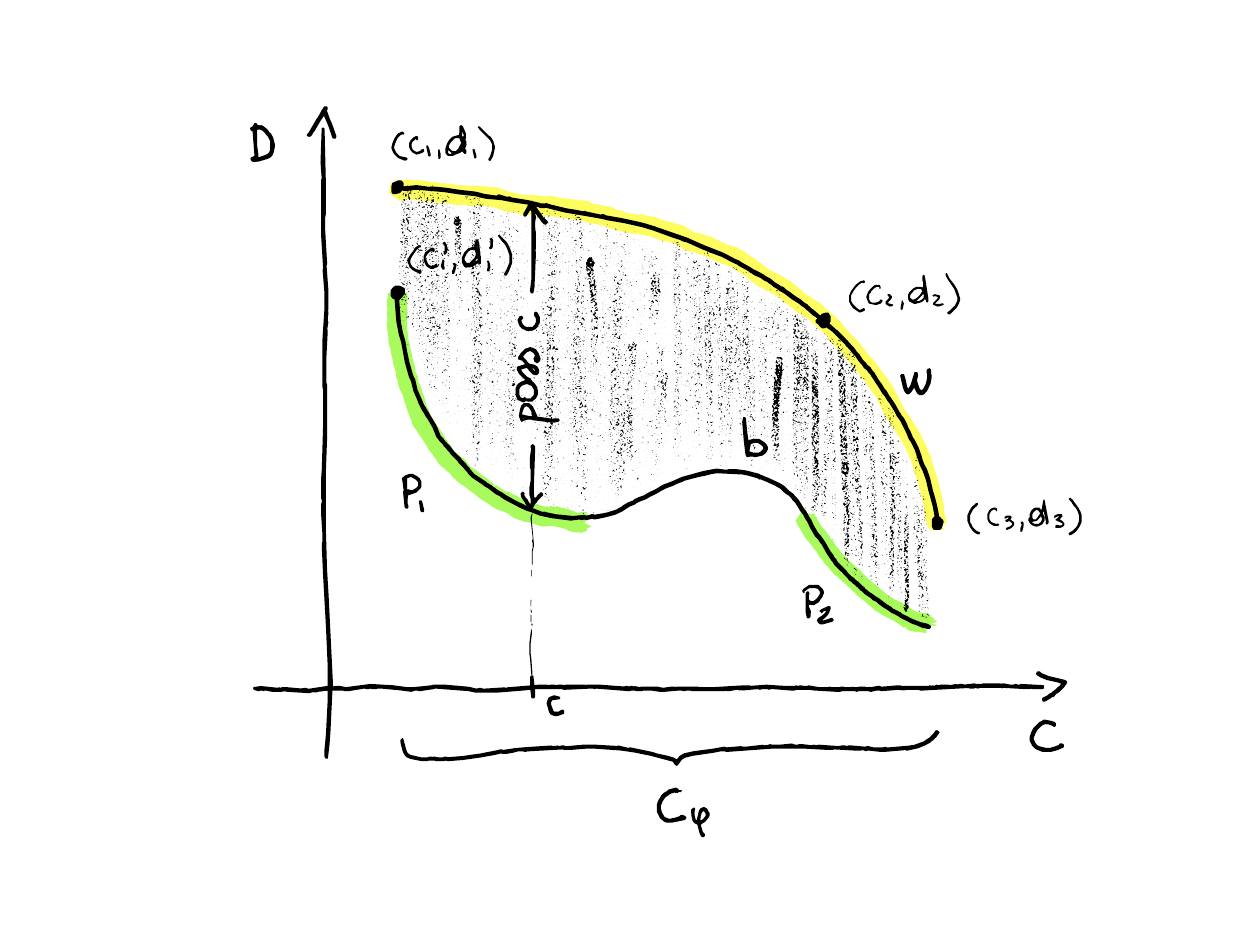}
  \caption{Pareto efficient manifolds for the worst-case (yellow) and
    for the best-case (green) measures: stable (1,3),
    unstable (1') and tipping (2) points.}
  \label{fig:four}
\end{figure}

As seen in \cref{section:monadic}, tackling decision problems affected
by functorial uncertainty requires selecting a measure function
\ensuremath{\mu}. This choice can have a crucial impact on the geometry of the
Pareto front of the resulting multi-objective optimization problem,
as illustrated in figure \cref{fig:four} for the best- and the
worst-case measure in a conceptual, stylized study.
The study is conceptual because it is not based on actual measurements
or on model data.
It is stylized because it encodes qualitative features that need to be
present in actual data for these to be relevant for decision making.
For example, damages must be relatively high for very low short-time
investments and vice-versa.
Furthermore, uncertainties for very low and for very high investments
must be significantly lower than uncertainties for average investments,
as sketched in \cref{fig:four}.
If this were not the case, the decision problem would be irrelevant:
there is little point in increasing short-time investments if there is
no confidence that very low (high) investments almost certainly yield
very high (low) damages!

The bottom line is that, in presence of value uncertainty and functorial
uncertainty, it is important to apply algorithms that are consistent,
well-understood and, ideally, generic and verified.
But it is also important to make sure that decision makers fully
understand the consequences of different attitudes towards risk (yellow
and green manifolds in figure \cref{fig:four}) and their implications for
the stability of optimal decisions (stable, unstable and tipping
points in figure \cref{fig:four}).
Crucially, this framework formalizes how stakeholders with identical
empirical data but different risk attitudes (e.g., worst-case vs.\
expected-value) can rationally arrive at contradictory optimal policies.
This insight suggests that societal disagreement often stems not from
ignoring scientific evidence, but from implicit differences in the
chosen measure of uncertainty.
Making these measures explicit allows for a more transparent and
rational democratic debate.

\DONE{Discuss how to apply methods for optimization under functorial
  uncertainty}

\DONE{Discuss the general case: value-uncertainty and functorial
  uncertainty together}

\section*{Acknowledgments}
%
The work presented in this paper heavily relies on free software, among
others on Haskell, GHC, git, vi, Emacs, \LaTeX\ and on the
FreeBSD and Debian GNU/Linux operating systems.
It is our pleasure to thank all developers of these excellent products.

\subsection*{Conflicts of Interest}
None.

\bibliographystyle{ACM-Reference-Format}
\bibliography{references}


\begin{thebibliography}{00}


\ifx \showCODEN    \undefined \def \showCODEN     #1{\unskip}     \fi
\ifx \showDOI      \undefined \def \showDOI       #1{#1}\fi
\ifx \showISBNx    \undefined \def \showISBNx     #1{\unskip}     \fi
\ifx \showISBNxiii \undefined \def \showISBNxiii  #1{\unskip}     \fi
\ifx \showISSN     \undefined \def \showISSN      #1{\unskip}     \fi
\ifx \showLCCN     \undefined \def \showLCCN      #1{\unskip}     \fi
\ifx \shownote     \undefined \def \shownote      #1{#1}          \fi
\ifx \showarticletitle \undefined \def \showarticletitle #1{#1}   \fi
\ifx \showURL      \undefined \def \showURL       {\relax}        \fi
\providecommand\bibfield[2]{#2}
\providecommand\bibinfo[2]{#2}
\providecommand\natexlab[1]{#1}
\providecommand\showeprint[2][]{arXiv:#2}

\bibitem[\protect\citeauthoryear{Althoff}{Althoff}{2010}]%
        {Reachability}
\bibfield{author}{\bibinfo{person}{Matthias Althoff}.}
  \bibinfo{year}{2010}\natexlab{}.
\newblock \showarticletitle{Reachability Analysis and its Application to the
  Safety Assessment of Autonomous Cars}.
\newblock \bibinfo{journal}{{\em PhD Thesis, Technical University of Munich\/}}
  (\bibinfo{year}{2010}).
\newblock


\bibitem[\protect\citeauthoryear{Bearpark, Hogan, and Hsiang}{Bearpark
  et~al\mbox{.}}{2025}]%
        {BearparkHoganHsiang2025}
\bibfield{author}{\bibinfo{person}{T. Bearpark}, \bibinfo{person}{D. Hogan},
  {and} \bibinfo{person}{S. Hsiang}.} \bibinfo{year}{2025}\natexlab{}.
\newblock \showarticletitle{Data anomalies and the economic commitment of
  climate change}.
\newblock \bibinfo{journal}{{\em Nature\/}}  \bibinfo{volume}{644}
  (\bibinfo{year}{2025}).
\newblock
Issue 8075.
\showDOI{%
\url{https://doi.org/10.1038/s41586-025-09320-4}}


\bibitem[\protect\citeauthoryear{Bellman}{Bellman}{1957}]%
        {bellman1957}
\bibfield{author}{\bibinfo{person}{Richard Bellman}.}
  \bibinfo{year}{1957}\natexlab{}.
\newblock \bibinfo{booktitle}{{\em Dynamic Programming}}.
\newblock \bibinfo{publisher}{Princeton University Press}.
\newblock


\bibitem[\protect\citeauthoryear{Ben-Haim}{Ben-Haim}{2006}]%
        {InfoGap}
\bibfield{author}{\bibinfo{person}{Yakov Ben-Haim}.}
  \bibinfo{year}{2006}\natexlab{}.
\newblock \showarticletitle{Info-Gap Decision Theory: Decisions Under Severe
  Uncertainty}.
\newblock \bibinfo{journal}{{\em Academic Press\/}} (\bibinfo{year}{2006}).
\newblock
\newblock
\shownote{Includes a retrospective essay and new material on hybrid
  uncertainties and robust-satisficing behavior.}


\bibitem[\protect\citeauthoryear{Ben-Tal, El~Ghaoui, and Nemirovski}{Ben-Tal
  et~al\mbox{.}}{2009}]%
        {RobustOpt}
\bibfield{author}{\bibinfo{person}{Aharon Ben-Tal}, \bibinfo{person}{Laurent
  El~Ghaoui}, {and} \bibinfo{person}{Arkadi Nemirovski}.}
  \bibinfo{year}{2009}\natexlab{}.
\newblock \showarticletitle{Robust Optimization}.
\newblock \bibinfo{journal}{{\em Mathematics of Operations Research\/}}
  \bibinfo{volume}{35}, \bibinfo{number}{3} (\bibinfo{year}{2009}),
  \bibinfo{pages}{1--39}.
\newblock


\bibitem[\protect\citeauthoryear{Bertsekas}{Bertsekas}{1995}]%
        {bertsekas1995}
\bibfield{author}{\bibinfo{person}{D.~P. Bertsekas}.}
  \bibinfo{year}{1995}\natexlab{}.
\newblock \bibinfo{booktitle}{{\em Dynamic Programming and Optimal Control}}.
\newblock \bibinfo{publisher}{Athena Scientific}, \bibinfo{address}{Belmont,
  Mass.}
\newblock


\bibitem[\protect\citeauthoryear{Bird and {de Moor}}{Bird and {de
  Moor}}{1997}]%
        {bird_algebra_1997}
\bibfield{author}{\bibinfo{person}{Richard Bird} {and} \bibinfo{person}{Oege
  {de Moor}}.} \bibinfo{year}{1997}\natexlab{}.
\newblock \bibinfo{booktitle}{{\em Algebra of Programming}}.
\newblock \bibinfo{publisher}{{Prentice-Hall,} Inc.}
\newblock


\bibitem[\protect\citeauthoryear{Birge and Louveaux}{Birge and
  Louveaux}{2011}]%
        {StochProg}
\bibfield{author}{\bibinfo{person}{John~R. Birge} {and}
  \bibinfo{person}{François Louveaux}.} \bibinfo{year}{2011}\natexlab{}.
\newblock \bibinfo{booktitle}{{\em Introduction to Stochastic Programming\/}
  (\bibinfo{edition}{2} ed.)}.
\newblock \bibinfo{publisher}{Springer}.
\newblock
\showISBNx{1461402360}


\bibitem[\protect\citeauthoryear{Botta, Brede, Crucifix, Ionescu, Jansson, Li,
  Mart{\'i}nez, and Richter}{Botta et~al\mbox{.}}{2023}]%
        {Botta2023MatterMost}
\bibfield{author}{\bibinfo{person}{Nicola Botta}, \bibinfo{person}{Nuria
  Brede}, \bibinfo{person}{Michel Crucifix}, \bibinfo{person}{Cezar Ionescu},
  \bibinfo{person}{Patrik Jansson}, \bibinfo{person}{Zheng Li},
  \bibinfo{person}{Marina Mart{\'i}nez}, {and} \bibinfo{person}{Tim Richter}.}
  \bibinfo{year}{2023}\natexlab{}.
\newblock \showarticletitle{Responsibility Under Uncertainty: Which Climate
  Decisions Matter Most?}
\newblock \bibinfo{journal}{{\em Environmental Modeling {\&} Assessment\/}}
  (\bibinfo{year}{2023}).
\newblock
\showISSN{1573-2967}
\showDOI{%
\url{https://doi.org/10.1007/s10666-022-09867-w}}


\bibitem[\protect\citeauthoryear{Botta, Jansson, and Ionescu}{Botta
  et~al\mbox{.}}{017b}]%
        {2017_Botta_Jansson_Ionescu}
\bibfield{author}{\bibinfo{person}{Nicola Botta}, \bibinfo{person}{Patrik
  Jansson}, {and} \bibinfo{person}{Cezar Ionescu}.}
  \bibinfo{year}{2017b}\natexlab{}.
\newblock \showarticletitle{Contributions to a computational theory of policy
  advice and avoidability}.
\newblock \bibinfo{journal}{{\em Journal of Functional Programming\/}}
  \bibinfo{volume}{27} (\bibinfo{year}{2017b}), \bibinfo{pages}{1--52}.
\newblock
\showISSN{0956-7968}
\showDOI{%
\url{https://doi.org/10.1017/S0956796817000156}}


\bibitem[\protect\citeauthoryear{Botta, Jansson, Ionescu, Christiansen, and
  Brady}{Botta et~al\mbox{.}}{017a}]%
        {2014_Botta_et_al}
\bibfield{author}{\bibinfo{person}{Nicola Botta}, \bibinfo{person}{Patrik
  Jansson}, \bibinfo{person}{Cezar Ionescu}, \bibinfo{person}{David~R.
  Christiansen}, {and} \bibinfo{person}{Edwin Brady}.}
  \bibinfo{year}{2017a}\natexlab{}.
\newblock \showarticletitle{Sequential decision problems, dependent types and
  generic solutions}.
\newblock \bibinfo{journal}{{\em Logical Methods in Computer Science\/}}
  \bibinfo{volume}{13}, \bibinfo{number}{1} (\bibinfo{year}{2017a}).
\newblock
\showDOI{%
\url{https://doi.org/10.23638/LMCS-13(1:7)2017}}


\bibitem[\protect\citeauthoryear{Botta, Jansson, and Richter}{Botta
  et~al\mbox{.}}{2025a}]%
        {MOOcodeRepo2025a}
\bibfield{author}{\bibinfo{person}{Nicola Botta}, \bibinfo{person}{Patrik
  Jansson}, {and} \bibinfo{person}{Tim Richter}.}
  \bibinfo{year}{2025}\natexlab{a}.
\newblock \bibinfo{title}{Associated code for the paper ``Optimization under
  uncertainty: understanding orders and testing programs with
  specifications''}.
\newblock   (\bibinfo{year}{2025}).
\newblock
\showURL{%
\url{https://gitlab.pik-potsdam.de/botta/papers/2025OUU}}
\newblock
\shownote{See
  \href{https://gitlab.pik-potsdam.de/botta/papers/2025OUU/haskell/Opt/MultiObj.lhs}{haskell/Opt/MultiObj.lhs}
  in the repository for the paper:
  \href{https://gitlab.pik-potsdam.de/botta/papers/2025OUU}{gitlab.pik-potsdam.de/botta/papers/2025OUU}.}


\bibitem[\protect\citeauthoryear{Botta, Jansson, and Richter}{Botta
  et~al\mbox{.}}{2025b}]%
        {MOOcodeRepo2025b}
\bibfield{author}{\bibinfo{person}{Nicola Botta}, \bibinfo{person}{Patrik
  Jansson}, {and} \bibinfo{person}{Tim Richter}.}
  \bibinfo{year}{2025}\natexlab{b}.
\newblock \bibinfo{title}{Associated code for the paper ``Optimization under
  uncertainty: understanding orders and testing programs with
  specifications''}.
\newblock   (\bibinfo{year}{2025}).
\newblock
\showURL{%
\url{https://gitlab.pik-potsdam.de/botta/papers/2025OUU}}
\newblock
\shownote{See
  \href{https://gitlab.pik-potsdam.de/botta/papers/2025OUU/agda/Functor/Basic.lagda}{agda/Functor/Basic.lagda}
  in the repository for the paper:
  \href{https://gitlab.pik-potsdam.de/botta/papers/2025OUU}{gitlab.pik-potsdam.de/botta/papers/2025OUU}.}


\bibitem[\protect\citeauthoryear{Brede and Botta}{Brede and Botta}{2021}]%
        {BREDE_BOTTA_2021}
\bibfield{author}{\bibinfo{person}{Nuria Brede} {and} \bibinfo{person}{Nicola
  Botta}.} \bibinfo{year}{2021}\natexlab{}.
\newblock \showarticletitle{On the correctness of monadic backward induction}.
\newblock \bibinfo{journal}{{\em Journal of Functional Programming\/}}
  \bibinfo{volume}{31} (\bibinfo{year}{2021}), \bibinfo{pages}{e26}.
\newblock
\showDOI{%
\url{https://doi.org/10.1017/S0956796821000228}}


\bibitem[\protect\citeauthoryear{Carlino, Giuliani, Tavoni, and
  Castelletti}{Carlino et~al\mbox{.}}{2020}]%
        {CARLINO202016593}
\bibfield{author}{\bibinfo{person}{Angelo Carlino}, \bibinfo{person}{Matteo
  Giuliani}, \bibinfo{person}{Massimo Tavoni}, {and} \bibinfo{person}{Andrea
  Castelletti}.} \bibinfo{year}{2020}\natexlab{}.
\newblock \showarticletitle{Multi-objective optimal control of a simple
  stochastic climate-economy model}.
\newblock \bibinfo{journal}{{\em IFAC-PapersOnLine\/}} \bibinfo{volume}{53},
  \bibinfo{number}{2} (\bibinfo{year}{2020}), \bibinfo{pages}{16593--16598}.
\newblock
\showISSN{2405-8963}
\showDOI{%
\url{https://doi.org/10.1016/j.ifacol.2020.12.786}}
\newblock
\shownote{21th IFAC World Congress.}


\bibitem[\protect\citeauthoryear{Cousot and Cousot}{Cousot and Cousot}{1977}]%
        {AbstractInterp}
\bibfield{author}{\bibinfo{person}{Patrick Cousot} {and}
  \bibinfo{person}{Radhia Cousot}.} \bibinfo{year}{1977}\natexlab{}.
\newblock \showarticletitle{Abstract Interpretation: A Unified Lattice Model
  for Static Analysis of Programs by construction or approximation of
  fixpoints}. In \bibinfo{booktitle}{{\em POPL}}. \bibinfo{pages}{238--252}.
\newblock


\bibitem[\protect\citeauthoryear{Deb, Pratap, Agarwal, and Meyarivan}{Deb
  et~al\mbox{.}}{2002}]%
        {NSGAII}
\bibfield{author}{\bibinfo{person}{K. Deb}, \bibinfo{person}{A. Pratap},
  \bibinfo{person}{S. Agarwal}, {and} \bibinfo{person}{T. Meyarivan}.}
  \bibinfo{year}{2002}\natexlab{}.
\newblock \showarticletitle{A Fast and Elitist Multiobjective Genetic
  Algorithm: {NSGA-II}}.
\newblock \bibinfo{journal}{{\em Trans. Evol. Comp\/}} \bibinfo{volume}{6},
  \bibinfo{number}{2} (\bibinfo{date}{April} \bibinfo{year}{2002}),
  \bibinfo{pages}{182--197}.
\newblock
\showISSN{1089-778X}
\showDOI{%
\url{https://doi.org/10.1109/4235.996017}}


\bibitem[\protect\citeauthoryear{Easterbrook}{Easterbrook}{2023}]%
        {Easterbrook_2023}
\bibfield{author}{\bibinfo{person}{Steve~M. Easterbrook}.}
  \bibinfo{year}{2023}\natexlab{}.
\newblock \bibinfo{booktitle}{{\em Computing the Climate: How We Know What We
  Know About Climate Change}}.
\newblock \bibinfo{publisher}{Cambridge University Press}.
\newblock


\bibitem[\protect\citeauthoryear{Emmerich and Deutz}{Emmerich and
  Deutz}{2018}]%
        {Emmerich+Deutz}
\bibfield{author}{\bibinfo{person}{Michael~TM Emmerich} {and}
  \bibinfo{person}{Andr{\'e}~H Deutz}.} \bibinfo{year}{2018}\natexlab{}.
\newblock \bibinfo{booktitle}{{\em A tutorial on multiobjective optimization:
  fundamentals and evolutionary methods}}. Vol.~\bibinfo{volume}{17}.
\newblock \bibinfo{publisher}{Springer}. 585--609 pages.
\newblock
\showDOI{%
\url{https://doi.org/10.1007/s11047-018-9685-y}}


\bibitem[\protect\citeauthoryear{Erwig and Kollmansberger}{Erwig and
  Kollmansberger}{2006}]%
        {10.1017/S0956796805005721}
\bibfield{author}{\bibinfo{person}{Martin Erwig} {and} \bibinfo{person}{Steve
  Kollmansberger}.} \bibinfo{year}{2006}\natexlab{}.
\newblock \showarticletitle{{FUNCTIONAL PEARLS}: Probabilistic Functional
  Programming in {Haskell}}.
\newblock \bibinfo{journal}{{\em Journal of Functional Programming\/}}
  \bibinfo{volume}{16}, \bibinfo{number}{1} (\bibinfo{date}{Jan.}
  \bibinfo{year}{2006}), \bibinfo{pages}{21--34}.
\newblock
\showISSN{0956-7968}
\showDOI{%
\url{https://doi.org/10.1017/S0956796805005721}}


\bibitem[\protect\citeauthoryear{Garc\'{\i}a and Fern\'{a}ndez}{Garc\'{\i}a and
  Fern\'{a}ndez}{2015}]%
        {SafeRL}
\bibfield{author}{\bibinfo{person}{Javier Garc\'{\i}a} {and}
  \bibinfo{person}{Fernando Fern\'{a}ndez}.} \bibinfo{year}{2015}\natexlab{}.
\newblock \showarticletitle{A Comprehensive Survey on Safe Reinforcement
  Learning}.
\newblock \bibinfo{journal}{{\em Journal of Machine Learning Research\/}}
  \bibinfo{volume}{16}, \bibinfo{number}{1} (\bibinfo{year}{2015}),
  \bibinfo{pages}{1437--1480}.
\newblock


\bibitem[\protect\citeauthoryear{Giry}{Giry}{1981}]%
        {giry1981}
\bibfield{author}{\bibinfo{person}{M. Giry}.} \bibinfo{year}{1981}\natexlab{}.
\newblock \showarticletitle{A categorial approach to probability theory}.
\newblock In \bibinfo{booktitle}{{\em Categorical Aspects of Topology and
  Analysis}}, \bibfield{editor}{\bibinfo{person}{B.~Banaschewski}} (Ed.).
  \bibinfo{series}{Lecture Notes in Mathematics}, Vol.~\bibinfo{volume}{915}.
  \bibinfo{publisher}{Springer}, \bibinfo{address}{Berlin},
  \bibinfo{pages}{68--85}.
\newblock


\bibitem[\protect\citeauthoryear{Gordon, Henzinger, Nori, and Rajamani}{Gordon
  et~al\mbox{.}}{2014}]%
        {ProbProg}
\bibfield{author}{\bibinfo{person}{{Andrew D.} Gordon},
  \bibinfo{person}{{Thomas A.} Henzinger}, \bibinfo{person}{{Aditya V.} Nori},
  {and} \bibinfo{person}{{Sriram K.} Rajamani}.}
  \bibinfo{year}{2014}\natexlab{}.
\newblock \showarticletitle{Probabilistic Programming}. In
  \bibinfo{booktitle}{{\em Proceedings of the on Future of Software
  Engineering}}. \bibinfo{publisher}{ACM}, \bibinfo{pages}{167--181}.
\newblock
\showDOI{%
\url{https://doi.org/10.1145/2593882.2593900}}


\bibitem[\protect\citeauthoryear{Haasnoot, Kwakkel, Walker, and ter
  Maat}{Haasnoot et~al\mbox{.}}{2013}]%
        {DAPP}
\bibfield{author}{\bibinfo{person}{Marjolijn Haasnoot}, \bibinfo{person}{Jan~H.
  Kwakkel}, \bibinfo{person}{Warren~E. Walker}, {and} \bibinfo{person}{Judith
  ter Maat}.} \bibinfo{year}{2013}\natexlab{}.
\newblock \showarticletitle{Dynamic Adaptive Policy Pathways: A Method for
  Crafting Robust Decisions for a Deeply Uncertain World}.
\newblock \bibinfo{journal}{{\em Global Environmental Change\/}}
  \bibinfo{volume}{23}, \bibinfo{number}{2} (\bibinfo{year}{2013}),
  \bibinfo{pages}{485--498}.
\newblock
\showDOI{%
\url{https://doi.org/10.1016/j.gloenvcha.2012.12.006}}


\bibitem[\protect\citeauthoryear{Helwegen, Wieners, Frank, and
  Dijkstra}{Helwegen et~al\mbox{.}}{2019}]%
        {esd-10-453-2019}
\bibfield{author}{\bibinfo{person}{Koen~G. Helwegen},
  \bibinfo{person}{Claudia~E. Wieners}, \bibinfo{person}{Jason~E. Frank}, {and}
  \bibinfo{person}{Henk~A. Dijkstra}.} \bibinfo{year}{2019}\natexlab{}.
\newblock \showarticletitle{{C}omplementing {$CO_2$} emission reduction by
  solar radiation management might strongly enhance future welfare}.
\newblock \bibinfo{journal}{{\em Earth System Dynamics\/}}
  \bibinfo{volume}{10}, \bibinfo{number}{3} (\bibinfo{year}{2019}),
  \bibinfo{pages}{453--472}.
\newblock
\showDOI{%
\url{https://doi.org/10.5194/esd-10-453-2019}}


\bibitem[\protect\citeauthoryear{Hoppe, Embreus, and Fülöp}{Hoppe
  et~al\mbox{.}}{2021}]%
        {HOPPE2021108098}
\bibfield{author}{\bibinfo{person}{Mathias Hoppe}, \bibinfo{person}{Ola
  Embreus}, {and} \bibinfo{person}{Tünde Fülöp}.}
  \bibinfo{year}{2021}\natexlab{}.
\newblock \showarticletitle{{DREAM}: A fluid-kinetic framework for tokamak
  disruption runaway electron simulations}.
\newblock \bibinfo{journal}{{\em Comput. Phys. Commun.\/}}
  \bibinfo{volume}{268} (\bibinfo{year}{2021}), \bibinfo{pages}{108098}.
\newblock
\showISSN{0010-4655}
\showDOI{%
\url{https://doi.org/10.1016/j.cpc.2021.108098}}


\bibitem[\protect\citeauthoryear{Hughes}{Hughes}{2010}]%
        {Hughes2010}
\bibfield{author}{\bibinfo{person}{John Hughes}.}
  \bibinfo{year}{2010}\natexlab{}.
\newblock \showarticletitle{Software Testing with {QuickCheck}}.
\newblock In \bibinfo{booktitle}{{\em CEFP 2009, Revised Selected Lectures}}.
  \bibinfo{publisher}{Springer}, \bibinfo{pages}{183--223}.
\newblock
\showISBNx{978-3-642-17685-2}
\showDOI{%
\url{https://doi.org/10.1007/978-3-642-17685-2_6}}


\bibitem[\protect\citeauthoryear{{Intergovernmental Panel on Climate
  Change}}{{Intergovernmental Panel on Climate Change}}{2018}]%
        {IPCC2018_SPM}
\bibfield{author}{\bibinfo{person}{{Intergovernmental Panel on Climate
  Change}}.} \bibinfo{year}{2018}\natexlab{}.
\newblock \showarticletitle{Summary for Policymakers}.
\newblock In \bibinfo{booktitle}{{\em Global Warming of 1.5\textdegree C.}},
  \bibfield{editor}{\bibinfo{person}{Valérie Masson-Delmotte} {et~al\mbox{.}}}
  (Eds.). \bibinfo{publisher}{Cambridge University Press},
  \bibinfo{pages}{3--24}.
\newblock
\showDOI{%
\url{https://doi.org/10.1017/9781009157940.001}}
\newblock
\shownote{See
  \href{https://doi.org/10.1017/9781009157940.001}{doi:10.1017/9781009157940.001}.}


\bibitem[\protect\citeauthoryear{Ionescu}{Ionescu}{2009}]%
        {ionescu2009}
\bibfield{author}{\bibinfo{person}{Cezar Ionescu}.}
  \bibinfo{year}{2009}\natexlab{}.
\newblock {\em \bibinfo{title}{Vulnerability Modelling and Monadic Dynamical
  Systems}}.
\newblock \bibinfo{thesistype}{Ph.D. Dissertation}. \bibinfo{school}{Freie
  Universit{\"a}t Berlin}.
\newblock
\showURL{%
\url{https://d-nb.info/1023491036/34}}


\bibitem[\protect\citeauthoryear{Ionescu}{Ionescu}{2016}]%
        {IONESCU_2016b}
\bibfield{author}{\bibinfo{person}{Cezar Ionescu}.}
  \bibinfo{year}{2016}\natexlab{}.
\newblock \showarticletitle{Vulnerability modelling with functional programming
  and dependent types}.
\newblock \bibinfo{journal}{{\em Mathematical Structures in Computer
  Science\/}} \bibinfo{volume}{26}, \bibinfo{number}{1} (\bibinfo{year}{2016}),
  \bibinfo{pages}{114–128}.
\newblock
\showDOI{%
\url{https://doi.org/10.1017/S0960129514000139}}


\bibitem[\protect\citeauthoryear{Ionescu and Jansson}{Ionescu and
  Jansson}{2013}]%
        {ionescujansson:LIPIcs:2013:3899}
\bibfield{author}{\bibinfo{person}{Cezar Ionescu} {and} \bibinfo{person}{Patrik
  Jansson}.} \bibinfo{year}{2013}\natexlab{}.
\newblock \showarticletitle{Testing versus proving in climate impact research}.
  In \bibinfo{booktitle}{{\em Proc. {TYPES} 2011}} {\em
  (\bibinfo{series}{Leibniz International Proceedings in Informatics
  (LIPIcs)})}, Vol.~\bibinfo{volume}{19}. \bibinfo{publisher}{Schloss Dagstuhl
  -- Leibniz-Zentrum f{\"u}r Informatik}, \bibinfo{address}{Dagstuhl, Germany},
  \bibinfo{pages}{41--54}.
\newblock
\showDOI{%
\url{https://doi.org/10.4230/LIPIcs.TYPES.2011.41}}


\bibitem[\protect\citeauthoryear{Jansson and Jansson}{Jansson and
  Jansson}{2023}]%
        {jansson_jansson_2023}
\bibfield{author}{\bibinfo{person}{Julia Jansson} {and} \bibinfo{person}{Patrik
  Jansson}.} \bibinfo{year}{2023}\natexlab{}.
\newblock \showarticletitle{Level-p-complexity of {Boolean} functions using
  thinning, memoization, and polynomials}.
\newblock \bibinfo{journal}{{\em Journal of Functional Programming\/}}
  \bibinfo{volume}{33} (\bibinfo{year}{2023}), \bibinfo{pages}{e13}.
\newblock
\showDOI{%
\url{https://doi.org/10.1017/S0956796823000102}}


\bibitem[\protect\citeauthoryear{Jansson, Smallbone, Botta, Sternvik, and
  Norozi}{Jansson et~al\mbox{.}}{2025}]%
        {jansson+2025}
\bibfield{author}{\bibinfo{person}{P. Jansson}, \bibinfo{person}{N. Smallbone},
  \bibinfo{person}{N. Botta}, \bibinfo{person}{E.~L. Sternvik}, {and}
  \bibinfo{person}{E. Norozi}.} \bibinfo{year}{2025}\natexlab{}.
\newblock \bibinfo{title}{Robust multi-objective optimization}.
\newblock \bibinfo{howpublished}{in preparation}.   (\bibinfo{year}{2025}).
\newblock


\bibitem[\protect\citeauthoryear{Kotz, Levermann, and Wenz}{Kotz
  et~al\mbox{.}}{2024}]%
        {KotzLevermannWenz2024}
\bibfield{author}{\bibinfo{person}{M. Kotz}, \bibinfo{person}{A. Levermann},
  {and} \bibinfo{person}{L. Wenz}.} \bibinfo{year}{2024}\natexlab{}.
\newblock \showarticletitle{The economic commitment of climate change}.
\newblock \bibinfo{journal}{{\em Nature\/}}  \bibinfo{volume}{628}
  (\bibinfo{year}{2024}).
\newblock
Issue 8008.
\showDOI{%
\url{https://doi.org/10.1038/s41586-024-07219-0}}


\bibitem[\protect\citeauthoryear{Mart\'{\i}nez~Montero, Crucifix, Couplet,
  Brede, and Botta}{Mart\'{\i}nez~Montero et~al\mbox{.}}{2022}]%
        {gmd-15-8059-2022}
\bibfield{author}{\bibinfo{person}{M. Mart\'{\i}nez~Montero},
  \bibinfo{person}{M. Crucifix}, \bibinfo{person}{V. Couplet},
  \bibinfo{person}{N. Brede}, {and} \bibinfo{person}{N. Botta}.}
  \bibinfo{year}{2022}\natexlab{}.
\newblock \showarticletitle{SURFER v2.0: a flexible and simple model linking
  anthropogenic CO$_2$ emissions and solar radiation modification to ocean
  acidification and sea level rise}.
\newblock \bibinfo{journal}{{\em Geoscientific Model Development\/}}
  \bibinfo{volume}{15}, \bibinfo{number}{21} (\bibinfo{year}{2022}),
  \bibinfo{pages}{8059--8084}.
\newblock
\showDOI{%
\url{https://doi.org/10.5194/gmd-15-8059-2022}}


\bibitem[\protect\citeauthoryear{Martínez~Montero, Brede, Couplet, Crucifix,
  Botta, and Wieners}{Martínez~Montero et~al\mbox{.}}{2024}]%
        {10.1093/oxfclm/kgae004}
\bibfield{author}{\bibinfo{person}{Marina Martínez~Montero},
  \bibinfo{person}{Nuria Brede}, \bibinfo{person}{Victor Couplet},
  \bibinfo{person}{Michel Crucifix}, \bibinfo{person}{Nicola Botta}, {and}
  \bibinfo{person}{Claudia Wieners}.} \bibinfo{year}{2024}\natexlab{}.
\newblock \showarticletitle{{Lost options commitment: how short-term policies
  affect long-term scope of action}}.
\newblock \bibinfo{journal}{{\em Oxford Open Climate Change\/}}
  \bibinfo{volume}{4}, \bibinfo{number}{1} (\bibinfo{date}{02}
  \bibinfo{year}{2024}), \bibinfo{pages}{kgae004}.
\newblock
\showISSN{2634-4068}
\showDOI{%
\url{https://doi.org/10.1093/oxfclm/kgae004}}
\showeprint{https://academic.oup.com/oocc/article-pdf/4/1/kgae004/56807120/kgae004.pdf}


\bibitem[\protect\citeauthoryear{Mayne, Rawlings, Rao, and Scokaert}{Mayne
  et~al\mbox{.}}{2000}]%
        {MPC}
\bibfield{author}{\bibinfo{person}{David~Q. Mayne}, \bibinfo{person}{James~B.
  Rawlings}, \bibinfo{person}{Christopher~V. Rao}, {and} \bibinfo{person}{Peter
  O.~M. Scokaert}.} \bibinfo{year}{2000}\natexlab{}.
\newblock \showarticletitle{Constrained Model Predictive Control: Stability and
  Optimality}.
\newblock \bibinfo{journal}{{\em Automatica\/}} \bibinfo{volume}{36},
  \bibinfo{number}{6} (\bibinfo{year}{2000}), \bibinfo{pages}{789--814}.
\newblock


\bibitem[\protect\citeauthoryear{McBride and Paterson}{McBride and
  Paterson}{2008}]%
        {mcbride2008applicative}
\bibfield{author}{\bibinfo{person}{Conor McBride} {and} \bibinfo{person}{Ross
  Paterson}.} \bibinfo{year}{2008}\natexlab{}.
\newblock \showarticletitle{Applicative Programming with Effects}.
\newblock \bibinfo{journal}{{\em Journal of Functional Programming\/}}
  \bibinfo{volume}{18}, \bibinfo{number}{1} (\bibinfo{date}{Jan.}
  \bibinfo{year}{2008}), \bibinfo{pages}{1--13}.
\newblock
\showISSN{0956-7968}
\showDOI{%
\url{https://doi.org/10.1017/S0956796807006326}}


\bibitem[\protect\citeauthoryear{Mitchell}{Mitchell}{1997}]%
        {10.5555/541177}
\bibfield{author}{\bibinfo{person}{Thomas~M. Mitchell}.}
  \bibinfo{year}{1997}\natexlab{}.
\newblock \bibinfo{booktitle}{{\em Machine Learning\/} (\bibinfo{edition}{1}
  ed.)}.
\newblock \bibinfo{publisher}{McGraw-Hill, Inc.}, \bibinfo{address}{USA}.
\newblock
\showISBNx{0070428077}


\bibitem[\protect\citeauthoryear{Moreno-Cruz and Keith}{Moreno-Cruz and
  Keith}{2012}]%
        {Moreno-Cruz+Keith2012}
\bibfield{author}{\bibinfo{person}{Juan Moreno-Cruz} {and}
  \bibinfo{person}{David Keith}.} \bibinfo{year}{2012}\natexlab{}.
\newblock \showarticletitle{Climate Policy under Uncertainty: A Case for
  Geoengineering}.
\newblock \bibinfo{journal}{{\em Climatic Change\/}} (\bibinfo{year}{2012}).
\newblock
\showDOI{%
\url{https://doi.org/10.1007/s10584-012-0487-4}}


\bibitem[\protect\citeauthoryear{Mu, Ko, and Jansson}{Mu et~al\mbox{.}}{2009}]%
        {mu2009algebra}
\bibfield{author}{\bibinfo{person}{Shin-Cheng Mu},
  \bibinfo{person}{Hsiang-Shang Ko}, {and} \bibinfo{person}{Patrik Jansson}.}
  \bibinfo{year}{2009}\natexlab{}.
\newblock \showarticletitle{Algebra of programming in {A}gda: dependent types
  for relational program derivation}.
\newblock \bibinfo{journal}{{\em Journal of Functional Programming\/}}
  \bibinfo{volume}{19}, \bibinfo{number}{5} (\bibinfo{year}{2009}),
  \bibinfo{pages}{545--579}.
\newblock


\bibitem[\protect\citeauthoryear{Neal, Newell, and Pitman}{Neal
  et~al\mbox{.}}{2025}]%
        {NealNewellPitman2025}
\bibfield{author}{\bibinfo{person}{T. Neal}, \bibinfo{person}{B.~R. Newell},
  {and} \bibinfo{person}{A. Pitman}.} \bibinfo{year}{2025}\natexlab{}.
\newblock \showarticletitle{Reconsidering the macroeconomic damage of severe
  warming}.
\newblock \bibinfo{journal}{{\em Environ. Res. Lett.\/}} \bibinfo{volume}{20},
  \bibinfo{number}{044029} (\bibinfo{year}{2025}).
\newblock
\showDOI{%
\url{https://doi.org/10.1088/1748-9326/adbd58}}


\bibitem[\protect\citeauthoryear{Nordhaus}{Nordhaus}{2018}]%
        {Nordhaus2018}
\bibfield{author}{\bibinfo{person}{William Nordhaus}.}
  \bibinfo{year}{2018}\natexlab{}.
\newblock \showarticletitle{{E}volution of modeling of the economics of global
  warming: changes in the {DICE} model, 1992--2017}.
\newblock \bibinfo{journal}{{\em Climatic Change\/}} \bibinfo{volume}{149},
  \bibinfo{number}{4} (\bibinfo{year}{2018}), \bibinfo{pages}{623--640}.
\newblock
\showDOI{%
\url{https://doi.org/10.1007/s10584-018-2218-y}}


\bibitem[\protect\citeauthoryear{Nordhaus}{Nordhaus}{2019}]%
        {Nordhaus12261}
\bibfield{author}{\bibinfo{person}{William Nordhaus}.}
  \bibinfo{year}{2019}\natexlab{}.
\newblock \showarticletitle{Economics of the disintegration of the Greenland
  ice sheet}.
\newblock \bibinfo{journal}{{\em Proceedings of the National Academy of
  Sciences\/}} \bibinfo{volume}{116}, \bibinfo{number}{25}
  (\bibinfo{year}{2019}), \bibinfo{pages}{12261--12269}.
\newblock
\showISSN{0027-8424}
\showDOI{%
\url{https://doi.org/10.1073/pnas.1814990116}}
\showeprint{https://www.pnas.org/content/116/25/12261.full.pdf}


\bibitem[\protect\citeauthoryear{Norell}{Norell}{2007}]%
        {norell2007thesis}
\bibfield{author}{\bibinfo{person}{Ulf Norell}.}
  \bibinfo{year}{2007}\natexlab{}.
\newblock {\em \bibinfo{title}{Towards a practical programming language based
  on dependent type theory}}.
\newblock \bibinfo{thesistype}{Ph.D. Dissertation}. \bibinfo{school}{Chalmers
  University of Technology}.
\newblock
\showDOI{%
\url{https://doi.org/10.1.1.436.7331}}


\bibitem[\protect\citeauthoryear{Pindyck}{Pindyck}{2017}]%
        {pindyck2017}
\bibfield{author}{\bibinfo{person}{Robert~S. Pindyck}.}
  \bibinfo{year}{2017}\natexlab{}.
\newblock \bibinfo{title}{{T}he {U}se and {M}isuse of {M}odels for {C}limate
  {P}olicy}.
\newblock \bibinfo{howpublished}{Review of Environmental Economics and Policy}.
    (\bibinfo{year}{2017}), \bibinfo{numpages}{100--114}~pages.
\newblock
Issue 1.
\showDOI{%
\url{https://doi.org/10.1093/reep/rew012}}


\bibitem[\protect\citeauthoryear{Pusztai, Ekmark, Bergström, Halldestam,
  Jansson, Hoppe, Vallhagen, and Fülöp}{Pusztai et~al\mbox{.}}{2023}]%
        {Pusztai2023BayesOptMMI}
\bibfield{author}{\bibinfo{person}{Istvan Pusztai}, \bibinfo{person}{Ida
  Ekmark}, \bibinfo{person}{Hannes Bergström}, \bibinfo{person}{Peter
  Halldestam}, \bibinfo{person}{Patrik Jansson}, \bibinfo{person}{Mathias
  Hoppe}, \bibinfo{person}{Oskar Vallhagen}, {and} \bibinfo{person}{Tünde
  Fülöp}.} \bibinfo{year}{2023}\natexlab{}.
\newblock \showarticletitle{Bayesian optimization of massive material injection
  for disruption mitigation in tokamaks}.
\newblock \bibinfo{journal}{{\em Journal of Plasma Physics\/}}
  \bibinfo{volume}{89}, \bibinfo{number}{2} (\bibinfo{year}{2023}),
  \bibinfo{pages}{905890204}.
\newblock
\showDOI{%
\url{https://doi.org/10.1017/S0022377823000193}}


\bibitem[\protect\citeauthoryear{Puterman}{Puterman}{2014}]%
        {puterman2014markov}
\bibfield{author}{\bibinfo{person}{Martin~L Puterman}.}
  \bibinfo{year}{2014}\natexlab{}.
\newblock \bibinfo{booktitle}{{\em Markov {D}ecision {P}rocesses: {D}iscrete
  {S}tochastic {D}ynamic {P}rogramming}}.
\newblock \bibinfo{publisher}{John Wiley \& Sons}.
\newblock


\bibitem[\protect\citeauthoryear{Richardson, Steffen, Lucht, Bendtsen, Cornell,
  Donges, Drüke, Fetzer, Bala, von Bloh, Feulner, Fiedler, Gerten, Gleeson,
  Hofmann, Huiskamp, Kummu, Mohan, Nogués-Bravo, Petri, Porkka, Rahmstorf,
  Schaphoff, Thonicke, Tobian, Virkki, Wang-Erlandsson, Weber, and
  Rockström}{Richardson et~al\mbox{.}}{2023}]%
        {doi:10.1126/sciadv.adh2458}
\bibfield{author}{\bibinfo{person}{Katherine Richardson}, \bibinfo{person}{Will
  Steffen}, \bibinfo{person}{Wolfgang Lucht}, \bibinfo{person}{Jørgen
  Bendtsen}, \bibinfo{person}{Sarah~E. Cornell}, \bibinfo{person}{Jonathan~F.
  Donges}, \bibinfo{person}{Markus Drüke}, \bibinfo{person}{Ingo Fetzer},
  \bibinfo{person}{Govindasamy Bala}, \bibinfo{person}{Werner von Bloh},
  \bibinfo{person}{Georg Feulner}, \bibinfo{person}{Stephanie Fiedler},
  \bibinfo{person}{Dieter Gerten}, \bibinfo{person}{Tom Gleeson},
  \bibinfo{person}{Matthias Hofmann}, \bibinfo{person}{Willem Huiskamp},
  \bibinfo{person}{Matti Kummu}, \bibinfo{person}{Chinchu Mohan},
  \bibinfo{person}{David Nogués-Bravo}, \bibinfo{person}{Stefan Petri},
  \bibinfo{person}{Miina Porkka}, \bibinfo{person}{Stefan Rahmstorf},
  \bibinfo{person}{Sibyll Schaphoff}, \bibinfo{person}{Kirsten Thonicke},
  \bibinfo{person}{Arne Tobian}, \bibinfo{person}{Vili Virkki},
  \bibinfo{person}{Lan Wang-Erlandsson}, \bibinfo{person}{Lisa Weber}, {and}
  \bibinfo{person}{Johan Rockström}.} \bibinfo{year}{2023}\natexlab{}.
\newblock \showarticletitle{Earth beyond six of nine planetary boundaries}.
\newblock \bibinfo{journal}{{\em Science Advances\/}} \bibinfo{volume}{9},
  \bibinfo{number}{37} (\bibinfo{year}{2023}), \bibinfo{pages}{eadh2458}.
\newblock
\showDOI{%
\url{https://doi.org/10.1126/sciadv.adh2458}}
\showeprint{https://www.science.org/doi/pdf/10.1126/sciadv.adh2458}


\bibitem[\protect\citeauthoryear{Rockström, Steffen, Noone, Persson, Chapin,
  Lambin, Lenton, Scheffer, Folke, Schellnhuber, Nykvist, De~Wit, Hughes,
  van~der Leeuw, Rodhe, Sörlin, Snyder, Costanza, Svedin, Falkenmark,
  Karlberg, Corell, Fabry, Hansen, Walker, Liverman, Richardson, Crutzen, and
  Foley}{Rockström et~al\mbox{.}}{2009}]%
        {ro06010m}
\bibfield{author}{\bibinfo{person}{J. Rockström}, \bibinfo{person}{W.
  Steffen}, \bibinfo{person}{K. Noone}, \bibinfo{person}{A. Persson},
  \bibinfo{person}{F.~S. Chapin}, \bibinfo{person}{E. Lambin},
  \bibinfo{person}{T.~M. Lenton}, \bibinfo{person}{M. Scheffer},
  \bibinfo{person}{C. Folke}, \bibinfo{person}{H. Schellnhuber},
  \bibinfo{person}{B. Nykvist}, \bibinfo{person}{C.~A. De~Wit},
  \bibinfo{person}{T. Hughes}, \bibinfo{person}{S. van~der Leeuw},
  \bibinfo{person}{H. Rodhe}, \bibinfo{person}{S. Sörlin},
  \bibinfo{person}{P.~K. Snyder}, \bibinfo{person}{R. Costanza},
  \bibinfo{person}{U. Svedin}, \bibinfo{person}{M. Falkenmark},
  \bibinfo{person}{L. Karlberg}, \bibinfo{person}{R.~W. Corell},
  \bibinfo{person}{V.~J. Fabry}, \bibinfo{person}{J. Hansen},
  \bibinfo{person}{B. Walker}, \bibinfo{person}{D. Liverman},
  \bibinfo{person}{K. Richardson}, \bibinfo{person}{P. Crutzen}, {and}
  \bibinfo{person}{J. Foley}.} \bibinfo{year}{2009}\natexlab{}.
\newblock \showarticletitle{Planetary boundaries: Exploring the safe operating
  space for humanity}.
\newblock \bibinfo{journal}{{\em Ecol. Soc.\/}} \bibinfo{volume}{14},
  \bibinfo{number}{2} (\bibinfo{year}{2009}), \bibinfo{pages}{32}.
\newblock


\bibitem[\protect\citeauthoryear{Schötz}{Schötz}{2025}]%
        {Schoetz2025}
\bibfield{author}{\bibinfo{person}{Christof Schötz}.}
  \bibinfo{year}{2025}\natexlab{}.
\newblock \showarticletitle{Spatial correlation in economic analysis of climate
  change}.
\newblock \bibinfo{journal}{{\em Nature\/}}  \bibinfo{volume}{644}
  (\bibinfo{year}{2025}).
\newblock
Issue 8076.
\showDOI{%
\url{https://doi.org/10.1038/s41586-025-09206-5}}


\bibitem[\protect\citeauthoryear{Sharpe, Mercure, Vinuales, Ives, Grubb,
  Pollitt, Knobloch, and Nijsse}{Sharpe et~al\mbox{.}}{2021}]%
        {shape2021}
\bibfield{author}{\bibinfo{person}{S. Sharpe}, \bibinfo{person}{J-F. Mercure},
  \bibinfo{person}{J. Vinuales}, \bibinfo{person}{M. Ives}, \bibinfo{person}{M.
  Grubb}, \bibinfo{person}{H. Pollitt}, \bibinfo{person}{F. Knobloch}, {and}
  \bibinfo{person}{F.J.M.M. Nijsse}.} \bibinfo{year}{2021}\natexlab{}.
\newblock \bibinfo{booktitle}{{\em {D}eciding how to decide: {R}isk-opportunity
  analysis as a generalisation of cost-benefit analysis}}.
\newblock \bibinfo{type}{{T}echnical {R}eport}. \bibinfo{institution}{UCL
  Institute for Innovation and Public Purpose, Working Paper Series (IIPP WP
  2021/03)}.
\newblock


\bibitem[\protect\citeauthoryear{Shepherd}{Shepherd}{2019}]%
        {shepherd2019}
\bibfield{author}{\bibinfo{person}{Theodore~G. Shepherd}.}
  \bibinfo{year}{2019}\natexlab{}.
\newblock \showarticletitle{{Storyline approach to the construction of regional
  climate change information}}.
\newblock \bibinfo{journal}{{\em Proc. R. Soc.\/}} \bibinfo{volume}{475},
  \bibinfo{number}{2225} (\bibinfo{year}{2019}).
\newblock
\showDOI{%
\url{https://doi.org/10.1098/rspa.2019.0013}}


\bibitem[\protect\citeauthoryear{Shepherd, Boyd, Calel, Chapman, Dessai,
  Dima-West, Fowler, James, Maraun, Martius, Senior, Sobel, Stainforth, Tett,
  Trenberth, van~den Hurk, Watkins, L., and Zenghelis}{Shepherd
  et~al\mbox{.}}{2018}]%
        {shepherd+al2018}
\bibfield{author}{\bibinfo{person}{T.~G. Shepherd}, \bibinfo{person}{E. Boyd},
  \bibinfo{person}{R.~A. Calel}, \bibinfo{person}{S.~C. Chapman},
  \bibinfo{person}{S. Dessai}, \bibinfo{person}{I.~M. Dima-West},
  \bibinfo{person}{H.~J. Fowler}, \bibinfo{person}{R. James},
  \bibinfo{person}{D. Maraun}, \bibinfo{person}{O. Martius},
  \bibinfo{person}{C.~A. Senior}, \bibinfo{person}{A.~H. Sobel},
  \bibinfo{person}{D.~A. Stainforth}, \bibinfo{person}{S.~F.~B. Tett},
  \bibinfo{person}{K.~E. Trenberth}, \bibinfo{person}{B.~J. J.~M. van~den
  Hurk}, \bibinfo{person}{N.~W. Watkins}, \bibinfo{person}{Wilby~R. L.}, {and}
  \bibinfo{person}{D.~A. Zenghelis}.} \bibinfo{year}{2018}\natexlab{}.
\newblock \showarticletitle{{Storylines: an alternative approach to
  representing uncertainty in physical aspects of climate change}}.
\newblock \bibinfo{journal}{{\em Climatic Change\/}} \bibinfo{number}{151}
  (\bibinfo{year}{2018}), \bibinfo{pages}{555–571}.
\newblock
\showDOI{%
\url{https://doi.org/10.1007/s10584-018-2317-9}}


\bibitem[\protect\citeauthoryear{Steffen, Richardson, Rockström, Cornell,
  Fetzer, Bennett, Biggs, Carpenter, de~Vries, de~Wit, Folke, Gerten, Heinke,
  Mace, Persson, Ramanathan, Reyers, and Sörlin}{Steffen
  et~al\mbox{.}}{2015}]%
        {doi:10.1126/science.1259855}
\bibfield{author}{\bibinfo{person}{Will Steffen}, \bibinfo{person}{Katherine
  Richardson}, \bibinfo{person}{Johan Rockström}, \bibinfo{person}{Sarah~E.
  Cornell}, \bibinfo{person}{Ingo Fetzer}, \bibinfo{person}{Elena~M. Bennett},
  \bibinfo{person}{Reinette Biggs}, \bibinfo{person}{Stephen~R. Carpenter},
  \bibinfo{person}{Wim de Vries}, \bibinfo{person}{Cynthia~A. de Wit},
  \bibinfo{person}{Carl Folke}, \bibinfo{person}{Dieter Gerten},
  \bibinfo{person}{Jens Heinke}, \bibinfo{person}{Georgina~M. Mace},
  \bibinfo{person}{Linn~M. Persson}, \bibinfo{person}{Veerabhadran Ramanathan},
  \bibinfo{person}{Belinda Reyers}, {and} \bibinfo{person}{Sverker Sörlin}.}
  \bibinfo{year}{2015}\natexlab{}.
\newblock \showarticletitle{Planetary boundaries: Guiding human development on
  a changing planet}.
\newblock \bibinfo{journal}{{\em Science\/}} \bibinfo{volume}{347},
  \bibinfo{number}{6223} (\bibinfo{year}{2015}), \bibinfo{pages}{1259855}.
\newblock
\showDOI{%
\url{https://doi.org/10.1126/science.1259855}}
\showeprint{https://www.science.org/doi/pdf/10.1126/science.1259855}


\bibitem[\protect\citeauthoryear{{The Idris Community}}{{The Idris
  Community}}{2020a}]%
        {idris2referenceguide}
\bibfield{author}{\bibinfo{person}{{The Idris Community}}.}
  \bibinfo{year}{2020}\natexlab{a}.
\newblock \bibinfo{title}{{I}dris2 {R}eference {G}uide}.
\newblock
  \bibinfo{howpublished}{\url{https://idris2.readthedocs.io/en/latest/reference/index.html}}.
    (\bibinfo{year}{2020}).
\newblock
\showURL{%
\url{https://idris2.readthedocs.io/en/latest/reference/index.html}}


\bibitem[\protect\citeauthoryear{{The Idris Community}}{{The Idris
  Community}}{2020b}]%
        {idrisreference}
\bibfield{author}{\bibinfo{person}{{The Idris Community}}.}
  \bibinfo{year}{2020}\natexlab{b}.
\newblock \bibinfo{title}{Language {Reference}}.
\newblock
  \bibinfo{howpublished}{\url{http://docs.idris-lang.org/en/latest/reference/index.html}}.
    (\bibinfo{year}{2020}).
\newblock
\showURL{%
\url{http://docs.idris-lang.org/en/latest/reference/index.html}}


\bibitem[\protect\citeauthoryear{{The Rocq Development Team}}{{The Rocq
  Development Team}}{2025}]%
        {CoqProofAssistant}
\bibfield{author}{\bibinfo{person}{{The Rocq Development Team}}.}
  \bibinfo{year}{2025}\natexlab{}.
\newblock \bibinfo{title}{Rocq Prover 9.0.0}.
\newblock   (\bibinfo{date}{March} \bibinfo{year}{2025}).
\newblock
\showDOI{%
\url{https://doi.org/10.5281/zenodo.1003420}}


\bibitem[\protect\citeauthoryear{Tol}{Tol}{2024}]%
        {TOL2024113922}
\bibfield{author}{\bibinfo{person}{Richard~S.J. Tol}.}
  \bibinfo{year}{2024}\natexlab{}.
\newblock \showarticletitle{A meta-analysis of the total economic impact of
  climate change}.
\newblock \bibinfo{journal}{{\em Energy Policy\/}}  \bibinfo{volume}{185}
  (\bibinfo{year}{2024}), \bibinfo{pages}{113922}.
\newblock
\showISSN{0301-4215}
\showDOI{%
\url{https://doi.org/10.1016/j.enpol.2023.113922}}


\bibitem[\protect\citeauthoryear{Varian}{Varian}{1992}]%
        {varian1992}
\bibfield{author}{\bibinfo{person}{Hal~R. Varian}.}
  \bibinfo{year}{1992}\natexlab{}.
\newblock \bibinfo{booktitle}{{\em {Microeconomic analysis}\/}
  (\bibinfo{edition}{third} ed.)}.
\newblock \bibinfo{publisher}{Norton}.
\newblock


\bibitem[\protect\citeauthoryear{Webster}{Webster}{2000}]%
        {webster2000}
\bibfield{author}{\bibinfo{person}{Mort~D. Webster}.}
  \bibinfo{year}{2000}\natexlab{}.
\newblock \bibinfo{booktitle}{{\em {T}he {C}urious {R}ole of "{L}earning" in
  {C}limate {P}olicy: {S}hould {W}e {W}ait for {M}ore {D}ata?}}
\newblock \bibinfo{type}{{T}echnical {R}eport}. \bibinfo{institution}{MIT Joint
  Program on the Science and Policy of Global Change, Report No. 67}.
\newblock


\bibitem[\protect\citeauthoryear{Webster}{Webster}{2008}]%
        {webster2008}
\bibfield{author}{\bibinfo{person}{Mort~D. Webster}.}
  \bibinfo{year}{2008}\natexlab{}.
\newblock \showarticletitle{{I}ncorporating {P}ath {D}ependency into
  {D}ecision-{A}nalytic {M}ethods: {A}n {A}pplication to {G}lobal
  {C}limate-{C}hange {P}olicy}.
\newblock \bibinfo{journal}{{\em Decision Analysis\/}} \bibinfo{volume}{5},
  \bibinfo{number}{2} (\bibinfo{year}{2008}), \bibinfo{pages}{60--75}.
\newblock


\bibitem[\protect\citeauthoryear{Zarnetske, Gurevitch, Franklin, Groffman,
  Harrison, Hellmann, Hoffman, Kothari, Robock, Tilmes, Visioni, Wu, Xia, and
  Yang}{Zarnetske et~al\mbox{.}}{2021}]%
        {Zarnetskee1921854118}
\bibfield{author}{\bibinfo{person}{Phoebe~L. Zarnetske},
  \bibinfo{person}{Jessica Gurevitch}, \bibinfo{person}{Janet Franklin},
  \bibinfo{person}{Peter~M. Groffman}, \bibinfo{person}{Cheryl~S. Harrison},
  \bibinfo{person}{Jessica~J. Hellmann}, \bibinfo{person}{Forrest~M. Hoffman},
  \bibinfo{person}{Shan Kothari}, \bibinfo{person}{Alan Robock},
  \bibinfo{person}{Simone Tilmes}, \bibinfo{person}{Daniele Visioni},
  \bibinfo{person}{Jin Wu}, \bibinfo{person}{Lili Xia}, {and}
  \bibinfo{person}{Cheng-En Yang}.} \bibinfo{year}{2021}\natexlab{}.
\newblock \showarticletitle{Potential ecological impacts of climate
  intervention by reflecting sunlight to cool Earth}.
\newblock \bibinfo{journal}{{\em Proceedings of the National Academy of
  Sciences\/}} \bibinfo{volume}{118}, \bibinfo{number}{15}
  (\bibinfo{year}{2021}).
\newblock
\showISSN{0027-8424}
\showDOI{%
\url{https://doi.org/10.1073/pnas.1921854118}}
\showeprint{https://www.pnas.org/content/118/15/e1921854118.full.pdf}


\bibitem[\protect\citeauthoryear{Zhang and Li}{Zhang and Li}{2007}]%
        {MOEAD}
\bibfield{author}{\bibinfo{person}{Qingfu Zhang} {and} \bibinfo{person}{Hui
  Li}.} \bibinfo{year}{2007}\natexlab{}.
\newblock \showarticletitle{{MOEA/D}: A Multiobjective Evolutionary Algorithm
  Based on Decomposition}.
\newblock \bibinfo{journal}{{\em IEEE Transactions on Evolutionary
  Computation\/}} \bibinfo{volume}{11}, \bibinfo{number}{6}
  (\bibinfo{year}{2007}), \bibinfo{pages}{712--731}.
\newblock
\showDOI{%
\url{https://doi.org/10.1109/TEVC.2007.892759}}


\end{thebibliography}
\label{lastpage01}


\end{document}